\documentclass[reqno,11pt]{amsart}

\usepackage{amsmath}
\usepackage{amsthm}
\usepackage{amssymb}
\usepackage{amscd}
\usepackage{xypic}
\usepackage{mathabx}
\usepackage{graphicx}
\usepackage{palatino}
\usepackage{bbm}
\usepackage{extarrows}
\usepackage{enumerate}

\usepackage[plainpages=false,colorlinks,hyperindex,pdfpagemode=None,bookmarksopen,linkcolor=red,citecolor=blue,urlcolor=blue]{hyperref}
\usepackage{pdflscape}
\usepackage{stmaryrd}

\usepackage{multirow}

\theoremstyle{plain}
\newtheorem{theorem}[subsubsection]{Theorem}

\newtheorem*{theorem*}{Theorem}
\newtheorem{proposition}[subsubsection]{Proposition}
\newtheorem*{proposition*}{Proposition}
\newtheorem{lemma}[subsubsection]{Lemma}
\newtheorem*{lemma*}{Lemma}

\newtheorem*{corollary*}{Corollary}

\theoremstyle{definition}
\newtheorem*{definition}{Definition}

\theoremstyle{remark}
\newtheorem{remark}[subsubsection]{Remark}

\newtheorem{example}[subsubsection]{Example}

\newcommand{\comment}[1] {  }

\DeclareFontFamily{OT1}{rsfs}{}
\DeclareFontShape{OT1}{rsfs}{n}{it}{<-> rsfs10}{}
\DeclareMathAlphabet{\mathscr}{OT1}{rsfs}{n}{it}

\newcommand{\Ad}{\mathrm{Ad}}

\newcommand{\CC}{\mathbb{C}}

\newcommand{\RR}{\mathbb{R}}

\newcommand{\B}{\mathcal{B}}

\newcommand{\Ind}{\operatorname{Ind}}
\newcommand{\Hom}{\operatorname{Hom}}
\newcommand{\End}{\operatorname{End}}
\newcommand{\Aut}{{\operatorname{Aut}}}

\newcommand{\Gm}{\mathbbm{G}_m}
\newcommand{\Ga}{\mathbbm{G}_a}

\newcommand{\GL}{\operatorname{GL}}
\newcommand{\Spin}{\operatorname{Spin}}
\newcommand{\Mat}{\operatorname{Mat}}

\newcommand{\PGL}{\operatorname{PGL}}
\newcommand{\SL}{\operatorname{SL}}
\newcommand{\Sp}{\operatorname{Sp}}

\newcommand{\SO}{{\operatorname{SO}}}

\newcommand{\tr}{\operatorname{tr}}
\newcommand{\Spec}{\operatorname{Spec}}

\newcommand{\diag}{{\operatorname{diag}}}

\newcommand{\ev}{\operatorname{ev}}

\newcommand{\reg}{{\operatorname{reg}}}

\newcommand{\Std}{{\operatorname{Std}}}

\newcommand{\codim}{{\operatorname{codim}}}

\renewcommand{\c}{\mathfrak{c}}

\newcommand{\g}{\mathfrak{g}}

\newcommand{\ua}{\underline{a}}
\newcommand{\ub}{\underline{b}}

\newcommand{\hor}{{\operatorname{hor}}}
\newcommand{\ver}{{\operatorname{ver}}}

\newcommand{\ssslash}{\mathbin{
  \mathchoice{/\mkern-10mu/\mkern-10mu/} 
    {/\mkern-10mu/\mkern-10mu/} 
    {/\mkern-8mu/\mkern-8mu/} 
    {/\mkern-8mu/\mkern-8mu/}}}

\newcommand{\bbbslash}{\mathbin{
  \mathchoice{\backslash\mkern-10mu\backslash\mkern-10mu\backslash} 
    {\backslash\mkern-10mu\backslash\mkern-10mu\backslash} 
    {\backslash\mkern-8mu\backslash\mkern-8mu\backslash} 
    {\backslash\mkern-8mu\backslash\mkern-8mu\backslash}}}

\makeatletter
\@namedef{subjclassname@2020}{\textup{2020} Mathematics Subject Classification}
\makeatother

\begin{document}
\title[Transfer operators and Hankel transforms]{Transfer operators and Hankel transforms: horospherical limits and quantization}
\author{Yiannis Sakellaridis}
\email{sakellar@jhu.edu}
\address{Department of Mathematics, Johns Hopkins University, Baltimore, MD 21218, USA.}


\begin{abstract}
Transfer operators are conjectural ``operators of functoriality,'' which transfer test measures and (relative) characters from one homogeneous space to another.  In previous work \cite{SaRankone, SaICM}, I computed transfer operators associated to spherical varieties of rank one, and gave an interpretation of them in terms of geometric quantization. In this paper, I study how these operators vary in the horospherical limits of these varieties, where they have a conceptual interpretation related to scattering theory. I also revisit Jacquet's Hankel transform for the Kuznetsov formula, which is related to the functional equation of the standard $L$-function of $\GL_n$, and provide an interpretation of it in terms of quantization.
\end{abstract}

\maketitle
\tableofcontents{}

\section{Introduction}

\subsection{Outline}
\subsubsection{} Let $G$ be a connected reductive group over a local field $F$, and, for a $G$-space $X$, consider the quotient stack $\mathfrak X = (X\times X)/G$, where $G$ acts diagonally. There is a notion of ``relative characters'' for $G$-representations associated to this quotient (generalizing characters on the adjoint quotient of the group), and many interesting questions that one can ask about them, such as 
\begin{quote}
 \item Are there relations between representations and characters associated to different spaces $\mathfrak X$, $\mathfrak Y$ as above?
\end{quote}

Such relations are often conjectured by the Langlands program and its generalization, the relative Langlands program.
Relative characters are, roughly, $G$-invariant generalized  functions on $X\times X$, or, equivalently, functionals on the Schwartz space $\mathcal S(\mathfrak X)$ \cite{SaStacks}, that are eigen- for the Bernstein center or the Harish-Chandra center. We would therefore like to answer such questions by describing ``transfer operators'' between Schwartz spaces
\[ \mathcal T: \mathcal S(\mathfrak Y) \to \mathcal S(\mathfrak X),\]
which pull back characters for $\mathfrak X$ to characters for $\mathfrak Y$. 

Such operators were described previously in \cite{SaRankone}, when $X$ is an affine homogeneous spherical variety of rank one, and $\mathfrak Y$ is the Kuznetsov quotient for the group $G^*=\SL_2$ or $\PGL_2$ (hence, the analogous to the space $X$ for $\mathfrak Y$ is the quotient of $G^*$ by a nontrivial unipotent subgroup, ``twisted'' by a nontrivial character $\psi$ of this subgroup). The surprising discovery was that, despite the non-abelian nature of the problem, the transfer operators $\mathcal T$ were given by explicit, abelian, Fourier transforms. A phenomenological ``quantization'' interpretation for this fact provided in \cite{SaICM}: namely, we can view $\mathcal S(\mathfrak Y)$ and $ \mathcal S(\mathfrak X)$ as ``quantizations'' of (roughly) the \emph{same} symplectic variety, and these Fourier transforms can be understood as direct analogs of the ``operators of change of Schr\"odinger model'' (albeit in a nonlinear setting). 

One goal of this paper is to provide another interpretation of these operators, in terms of the ``asymptotic cones'' (boundary degenerations) of $\mathfrak X$ and $\mathfrak Y$. In the case of $\mathfrak X$, this is the degeneration that one obtains when one lets the spherical variety $X$ degenerate to its horospherical limit, and in the case of $\mathfrak Y$, when one lets the character $\psi$ degenerate to the trivial character. In Section \ref{sec:degeneration}, we generalize to all rank-one (affine homogeneous) varieties a discovery of \cite[\S 4.3, 5]{SaTransfer1} for the special cases $X=\Gm\backslash \PGL_2$ and $X=\SL_2$, namely, 
that 
\emph{the transfer operators for $\mathfrak X$ and its boundary degeneration have exactly the same form} (in appropriate coordinates); see Theorem \ref{thmTempty}.

\subsubsection{} The argument for the study of transfer operators in the horospherical limit is no different than the one used in \cite{SaTransfer1}, up to knowing the \emph{scattering operators} in each case. Those are the operators that control the \emph{asymptotics} of generalized matrix coefficients on $X$, and the \emph{Plancherel formula} for the ``most continuous'' part of $L^2(X)$ \cite{SV, DHS}. Thus, we start in Section \ref{sec:scattering} by computing the scattering operators for all spherical varieties of rank one. The calculation should be very familiar to anyone with experience in the calculation of ``unramified/spherical (eigen)functions'' (for the Hecke algebra): Indeed, the scattering operators are the ``functional equations'' of the Casselman--Shalika method and its generalizations. Thus, our methods do not differ significanty from the ones employed in \cite{SaSph} to compute spherical functions. 

To extend this calculation to the ramified (principal series) representations, we need to be very pedantic with the definitions of Radon transforms (intertwining operators), to eliminate scalar ambiguities. Unfortunately, this makes that section quite technical, and the reader might be wise to just skim through it at first reading. On the flip side, we obtain a very rigid formula, encoded in Theorem \ref{thmscattering}, that relates scattering operators to the \emph{gamma factors} of the local functional equation of $L$-functions. These $\gamma$-factors will be our first encounter, in this paper, of the deep relations between local harmonic analysis and local $L$-functions.

\subsubsection{}

Of similar nature to the transfer operators are the so-called \emph{Hankel transforms}, which are the trace-formula-theoretic incarnations of the functional equations of local $L$-functions. For the purposes of the present paper, we do not need to recall general definitions or conjectures around those, as we will only be concerned with the standard $L$-function of $\GL_n$. By the work of Godement and Jacquet \cite{GJ}, the local functional equation for this $L$-function is afforded by the Fourier transform (depending on an additive character $\psi$)
\[ \mathcal F: \mathcal D(\Mat_n) \xrightarrow\sim \mathcal D(\Mat_n^*),\]
between half-densities on the vector space of $n\times n$ matrices, and on its dual, both viewed as $G=\GL_n\times^\Gm \GL_n$-spaces. The natural embeddings of $\GL_n$ in both $\Mat_n$ and $\Mat_n^*$ allow us to view $\mathcal F$ as a $G$-equivariant map between certain spaces of half-densities (or, by fixing a Haar measure, of measures) on $\GL_n$. It therefore acts by a scalar on characters of irreducible representations $\pi$ (at least, a scalar varying meromorphically, as the representation is twisted by characters of the determinant, and therefore defined for almost every $\pi$), and this scalar is, by definition, the gamma factor
\[ \gamma(\pi,\frac{1}{2},\psi)\]
of the standard (local) $L$-function of $\pi$.
Because of its equivariant nature, the Fourier transform descends to a \emph{Hankel transform} between spaces $(N,\psi)^2$-coinvariants (where $N$ is the upper triangular unipotent subgroup, and $\psi$ now also denotes a generic character of it)
\[ \mathcal H: \mathcal D(\Mat_n)_{(N,\psi)^2} \xrightarrow\sim \mathcal D(\Mat_n^*)_{(N,\psi)^2},\]
and we can now think of both sides as spaces $\mathcal D^-$, $\mathcal D^+$ of half-densities for the Kuznetsov quotient stack $\mathfrak Y$ of $\GL_n$. Jacquet \cite{Jacquet} has computed a formula for this transform -- see Theorem \ref{thmJacquet}.

In this paper, we will give a quantization interpretation for this formula, Theorem \ref{thmreform}, similar to the one given in \cite{SaICM} to the rank-one transfer operators. Namely, we will view the spaces $\mathcal D^-$, $\mathcal D^+$ as two different ``geometric quantizations'' of the \emph{same} cotangent stack $\mathfrak M$ (the latter being the two-sided Whittaker reduction of $T^*\Mat_n = T^* \Mat_n^*$), given by two different Lagrangian foliations on it. And we will show that Jacquet's Hankel transform is given by the integral along the leaves of these foliations, just as in the case of intertwining operators between Schr\"odinger models for the oscillator representation. This point of view also allows us to give a geometric reformulation of Jacquet's proof.

\subsubsection{}

I view the results of this paper as further evidence for the microlocal nature of conjectural ``operators of functoriality.'' This idea appeared already in my earlier work mentioned above, but the study of Jacquet's Hankel transform given here is the first time that it is being confirmed in higher rank. Similar ideas, but not in the context of trace formulas, have appeared in talks and unpublished notes of Vincent Lafforgue. 

Moreover, concepts such as ``geometric quantization of symplectic stacks'' are essentially unexplored, and presently very vague. This paper provides some examples and hints as to what they might mean. 

\subsection{Acknowledgments}

It is my pleasure to dedicate this paper to Toshiyuki Kobayashi, on the occasion of his 60th birthday. I met Professor Kobayashi in 2007, when I visited him with Joseph Bernstein in Kyoto, and I have since enjoyed the privilege of talking to him at various occasions in Japan, Israel, and elsewhere. I have always admired his originality and independence as a mathematical thinker, but also the breadth and depth of his knowledge beyond mathematics, which ranges from Japanese mountain vegetables to Greek mythology. 

I am deeply indebted to the anonymous referees for pointing out many errors in the original submission and for substantial suggestions on exposition. 

This work was supported by NSF grant DMS-2101700; the results of Sections \ref{sec:scattering}--\ref{sec:degeneration} were previously announced, without proof, in the first arXiv version of \cite{SaRankone}, but were not included in the published version.

\subsection{Notation}

We will be working over a local field $F$, which is non-Archi\-medean of characteristic zero in Sections \ref{sec:scattering}, \ref{sec:degeneration}, and Archimedean in Section \ref{sec:Hankel}. When no confusion arises, we will simply write $X$ for the $F$-points of a variety $X$ over $F$. 

Once we fix a Haar measure on $F$, every volume form $\omega$ on a smooth variety $X$ gives rise to a density (measure) $|\omega|$ on its $F$-points, in the standard way. 
We mostly follow the standard habit of denoting differentials and the associated densities by the same symbols, $dx$ etc., when the meaning is clear from the context; when not, we write $|\omega|$ for the density. One can also define associated half-densities, which will be denoted by $|\omega|^\frac{1}{2}$. 

We write $X\sslash G$ for the invariant-theoretic quotient $\Spec F[X]^G$ of a (typically affine) variety $X$ by the right action of a group $G$. When quotients are denoted by a single slash, $X/G$, unless we say otherwise, we will mean the stack quotient. However, knowledge of stacks is not required to read this paper: when this is not a variety, this symbol will mostly be a placeholder for an explicitly-defined object associated to the $G$-action on $X$, such as a space of orbital integrals. 

We will sometimes switch right actions to left actions, via the rule $g^{-1}\cdot x = x\cdot g$. The quotient of the product $X\times Y$ two right $G$-varieties by the diagonal action of $G$ will be denoted by $X\times^G Y$. 

The ``universal'' or ``abstract'' Cartan of a reductive group $G$ is defined as the quotient $A=B/N$ of \emph{any} Borel subgroup $B$ by its unipotent radical $N$; different choices for the Borel give a canonically isomorphic torus, which also comes equipped with a based root datum. Many constructions in this paper are ``universal'' in this sense: they can be described using a Borel subgroup, but a different choice leads to a canonically isomorphic construction. In those cases, we will feel free to use a Borel subgroup $B$, without commenting on the choice.

Finally, if $M$ is a Hamiltonian $G$-space (i.e., a symplectic $G$-variety equip\-ped with a moment map to $\mathfrak g^*$, for some group $G$), and $f\in \mathfrak g^*$ is a $G$-invariant element, the Hamiltonian reduction of $M$ at $f$, denoted $M\ssslash_f G$, is the symplectic ``space'' $(M \times_{\mathfrak g^*} \{f\})/G$. Often, this quotient does not make sense as a variety, i.e., the $G$-action is not free. In those cases, the proper way to think of $M\ssslash_f G$ is as a derived stack, i.e., we also need to understand the fiber product over $\mathfrak g^*$ as a derived fiber product. However, in this paper, we will usually restrict to an open subset which is a variety, or else explain some less sophisticated way of using this quotient.

\section{Scattering operators in rank one} \label{sec:scattering}

\subsection{Spherical varieties of rank one, and their asymptotic cones} \label{ssasymptoticcones}

\subsubsection{}
Let $X=H\backslash G$ be an affine, homogeneous spherical variety of rank one over a local field $F$, with $G$ and $H$ split reductive groups. In this section and the next, we will assume that $F$ is non-Archimedean in characteristic zero, because this is where a ``cleaner'' theory of asymptotics is available, by \cite[\S~5]{SV}. 

\begin{remark} There are no serious obstacles to extending the theory of asymptotics to positive characteristic, at least for the spaces of Table \eqref{thetable} below, other than some care that needs to be taken when defining these spaces and their boundary degenerations in small characteristics. Moreover, the calculations of scattering operators in the present chapter extend to the scattering operators for \emph{tempered} representations, developed in \cite{Carmona}, \\ \cite{DKS}, in the Archimedean case.
\end{remark}

We will assume that $X$ is contained in the following table, which contains all such varieties, up to the action of the ``center'' $\mathcal Z(X):=\Aut^G(X)$.

\vspace{0.5cm}~

~\hspace{-3.5cm}
{\small
$ \begin{array}{|c|c|c|c|c|c|}
\hline
&X & P(X) & \check G_X & \gamma & L_X \\
 & & & & &\\
\hline
\hline
A_1& \Gm \backslash\PGL_{2} &B &  \SL_{2} & \alpha & L(\Std,\frac{1}{2})^2 \\
\hline
A_n& \GL_n \backslash\PGL_{n+1} &P_{1,n-1,1} &  \SL_{2} & \alpha_1+\dots + \alpha_n& L(\Std,\frac{n}{2})^2 \\
\hline
B_n& \SO_{2n}\backslash\SO_{2n+1} & P_{\SO_{2n-1}} & \SL_2 & \alpha_1 + \dots + \alpha_n&  L(\Std, n-\frac{1}{2}) L(\Std,\frac{1}{2})\\
\hline
C_n& \Sp_{2n-2}\times \Sp_{2}\backslash \Sp_{2n} & P_{\SL_2 \times \Sp_{2(n-2)}} & \SL_{2} & \alpha_1+2\alpha_2+\dots + 2\alpha_{n-1}+ \alpha_n & L(\Std,n-\frac{1}{2}) L(\Std,n-\frac{3}{2}) \\ 
\hline
F_4& \Spin_9\backslash F_4 & P_{\Spin_7} & \SL_2 & \alpha_1+2\alpha_2+3\alpha_3+2\alpha_4 & L(\Std, \frac{11}{2}) L(\Std, \frac{5}{2}) \\
\hline
G_2& \SL_3\backslash G_2& P_{\SL_2} & \SL_2 & 2\alpha_1+ \alpha_2  & L(\Std, \frac{5}{2}) L(\Std, \frac{1}{2}) \\
\hline
\hline
D_2& \SL_2=\SO_{3}\backslash\SO_{4} & B & \PGL_2 & \alpha_1 + \alpha_2 &  L(\Ad,1) \\
\hline
D_n& \SO_{2n-1}\backslash\SO_{2n} & P_{\SO_{2n-2}} & \PGL_2 & 2\alpha_1+ \dots + 2\alpha_{n-2} + \alpha_{n-1}+\alpha_n &  L(\Ad,n-1) \\
\hline
D_4'' & \Spin_7\backslash\Spin_8 & P_{\Spin_6} & \PGL_2 & 2\alpha_1+ 2\alpha_{2} + \alpha_{3}+\alpha_4 & L(\Ad,3)\\
\hline
B_3''& G_2\backslash\Spin_7 & P_{\SL_3} & \PGL_2 & \alpha_1 + 2\alpha_2 + 3\alpha_3 &  L(\Ad,3)\\
\hline
\end{array}
$
}
\begin{equation}\label{thetable}\end{equation}

The table also shows the dual group $\check G_X$ of $X$, with its positive coroot $\gamma$ (the ``normalized spherical root'' of $X$) and the ``associated $L$-value.'' The dual group admits an embedding into the dual group $\check G$ of $G$ with positive coroot $\gamma$, and its Weyl group $W_X\simeq \mathbb Z/2$ has a generator $w$ which can also be considered as an element of the Weyl group of $G$. 

Associating $L$-values (i.e., special values of $L$-functions -- in this case, of local $L$-factors) to these varieties is motivated by number theory, but it turns out that these mysterious quantities control a great deal of harmonic analysis. In this paper, we will see how they control \emph{scattering operators} on the asymptotic cone of $X$, and \emph{transfer operators} for its relative trace formula. We postpone the discussion of these $L$-values until we encounter them in harmonic analysis.

\subsubsection{}

The purpose of this section is to recall and precisely calculate scattering operators, for all the varieties above.  We must first recall the \emph{asymptotic cone} (or \emph{boundary degeneration}, in the language of \cite{SV}) of a (quasiaffine, homogeneous) spherical variety $X$. This is a \emph{horospherical} variety $X_\emptyset$ (which, here, we will take to be homogeneous, by definition), which can be defined in several equivalent ways. One of them is by identifying $X_\emptyset$ with the open $G$-orbit in the normal bundle to a closed $G$-orbit in a ``wonderful'' (or rather, smooth toroidal) compactification of $X$ \cite[\S~2.4]{SV}. In the case at hand, where our varieties are of rank $1$, they all possess a canonical ``wonderful'' compactification $\bar X$, which is the union of $X$ and a projective orbit isomorphic to the flag variety $\B_X=P(X)^-\backslash G$, for some  parabolic $P(X)^-$. The opposite $P(X)$ of this parabolic (or rather, its conjugacy class) can be characterized as  the stabilizer of the open Borel orbit $X^\circ\subset X$, and it admits a quotient to a torus $A_X$, which is the quotient by which it acts on $X^\circ/U_{P(X)}$. In the rank-one cases of the table above, $A_X\simeq \Gm$, either via the spherical root $\gamma$ (when $\check G_X=\SL_2$) or via its square root (when $\check G_X=\PGL_2$). Let $q:P(X)\to A_X=\Gm$ denote the quotient map.

When $X$ is symmetric (as are almost all\footnote{In fact, as abstract varieties, all are symmetric; but the varieties $\SL_3\backslash G_2$ and $G_2\backslash \Spin_7$ are symmetric under a \emph{larger} group of automorphisms, namely, they are isomorphic to $\SO_6\backslash\SO_7$ and $\SO_7\backslash\SO_8$, respectively.} of the varieties in Table \eqref{thetable}), with $\theta=$ the involution associated to a point in the open orbit for a chosen Borel subgroup, $P(X)$ is known as the ``minimal $\theta$-split parabolic.'' The roots in the Levi of $P(X)$, for the varieties of Table \eqref{thetable}, are those that are orthogonal to the spherical root $\gamma$, and the boundary degeneration is an $A_X$-torsor over $P(X)^-\backslash G$.

The quotient $q:P(X)\to A_X=\Gm$, also defines a character $q^-$ for the opposite parabolic $P(X)^-$, by the natural identification of the abelianizations of $P(X)$ and $P(X)^-$. The boundary degeneration is an $A_X\times G$-variety which can be identified with 
\begin{equation}\label{Xempty}
 X_\emptyset = S\backslash G,
\end{equation}
where $S=\ker q^-$, and our convention is that $A_X$ acts on $X_\emptyset$ via $A_X\simeq P(X)^-/S$. For the cases of Table \eqref{thetable}, $P(X)$ is self-dual, i.e.,  $P(X)^-$ is conjugate to $P(X)$, so we could have written $X_\emptyset$ as an $A_X$-torsor over $P(X)\backslash G$. However, our presentation helps recall our conventions for the $A_X$-action, and the action of the universal Cartan $A$ of $G$: While $A_X$ is identified as a quotient of the abstract Cartan via $B\to P(X)\to A_X$, it acts on $X_\emptyset$ via $P(X)^-\to A_X$. 

Note that the identification \eqref{Xempty} is by no means canonical -- any translation by the action of $A_X$ leads to another such identification. But the boundary degeneration $X_\emptyset$ is rigid, by construction, and this has consequences for our geometric and harmonic-analytic calculations (which will not always be invariant under the action of $A_X$). 
In particular, by \cite[Theorem 5.1.1]{SV}, there is a canonical ``asymptotics'' morphism 
\begin{equation}\label{asymptotics} 
e_\emptyset^*: C^\infty(X)\to C^\infty(X_\emptyset),
\end{equation}
with the property that, in a suitable sense (that we will not review here), $\Phi$ and $e_\emptyset^*\Phi$ ``coincide in a neighborhood of the orbit at infinity.''

As another manifestation of the rigidity of $X_\emptyset$, in a subsequent subsection we will show that there is a distinguished $G$-orbit $X_\emptyset^R$ in the ``open Bruhat cell'' of $X_\emptyset\times X_\emptyset$, that is, whose image in $\B_X\times \B_X$ under the map induced from $X_\emptyset\to \B_X$ (the map taking a point to the normalizer of its stabilizer) belongs to the open $G$-orbit.

\subsubsection{} \label{sssaffinedegen} Before we do that, let us recall an alternative, equivalent definition of $X_\emptyset$: It is the open $G$-orbit in the special fiber of the ``affine degeneration'' of $X$ \cite[2.5]{SV}. This is a family $\mathcal X\to \overline{A_X}$ of affine varieties over a certain affine embedding of the torus $A_X$, which in our case can be identified with $\mathbb A^1$. The fibers over $A_X$ are isomorphic to $X$, while the fiber over $0\in \mathbb A^1$ is an affine horospherical variety, whose open $G$-orbit can be identified with $X_\emptyset$. 

\begin{example}
 Let $X=\SO_n\backslash \SO_{n+1}$, identified with the ``unit sphere'' in a maximally isotropic quadratic space $V$ of the appropriate discriminant (if we want $\SO_n$ to be split, too). Then, in one version of the affine degeneration, we can consider $V$ itself as a family of spaces containing $X$, with the  map to $\mathbb A^1$ being the quadratic form. Note, however, that the fibers over $t\ne \mathbb A^1$ are only isomorphic over the algebraic closure, and depend on the square class of $t$ over $F$. Therefore, it may be arithmetically preferable to define $\mathcal X = V\times_{\mathbb A^1} \mathbb A^1$, with $\mathbb A^1 \to \mathbb A^1$ the square map. This is the family over $\overline{A_X}$ mentioned above. 
\end{example}

\subsection{Scattering operators}

\subsubsection{} We generally normalize actions of various groups to be $L^2$-unitary. Specifically for the action of $A_X$ on functions on $X_\emptyset$, 
if $\delta_{P(X)}$ denotes the modular character of $P(X)$ (the inverse of the modular character of $P(X)^-$, when considered as a character of their common Levi subgroup), then 
this unitary action is given by
\begin{equation}\label{norm-functions} a\cdot \Phi(Sg) = \delta_{P(X)}^{\frac{1}{2}}(a) \Phi(Sag).
\end{equation}

For $\chi$ a character of $A_X$ in general position, consider the normalized (possibly degenerate) principal series representation 
\[I_{P(X)^-}^G(\chi) = \Ind_{P(X)^-}^G (\chi\delta_{P(X)}^{-\frac{1}{2}}).\]
Up to a choice of scalar, we have an isomorphism 
\begin{equation}\label{isomprincipal}
 C^\infty((A_X,\chi)\backslash X_\emptyset) \simeq I_{P(X)^-}^G(\chi),
\end{equation}
where the space on the left is the space of smooth $(A_X,\chi)$-eigenfunctions under the normalized action above.

The scattering maps that we would like to describe form a  meromorphic (in $\chi$) family of morphisms:
\[ \mathscr S_{w,\chi} : C^\infty((A_X,\chi^{-1})\backslash X_\emptyset) \to C^\infty((A_X,\chi)\backslash X_\emptyset),\]
such that\footnote{We rely throughout on the fact that, in rank one, ${^w\chi}=\chi^{-1}$; in general the scattering operators relate induced representations from $W_X$-conjugate characters.} $\mathscr S_{w,\chi^{-1}} \circ \mathscr S_{w,\chi} = I$, and characterized as follows:

For $\chi$ in general position, there is a unique up to scalar morphism $I_{P(X)^-}^G(\chi)\to C^\infty(X)$, and the scattering morphism $\mathscr S_{w,\chi}$ encodes its asymptotics, in the sense that the composition of this embedding with the asymptotics morphism \eqref{asymptotics} lives in the space of $\mathscr S_w$-invariant pairs of the direct sum
\[ C^\infty((A_X,\chi^{-1})\backslash X_\emptyset) \oplus C^\infty((A_X,\chi)\backslash X_\emptyset).\]

More precisely, the following commutative diagram characterizes the scattering morphisms, see \cite[\S~10.17]{DHS}:
\begin{equation}\label{fiberscatteringdiagram} \xymatrix{ && C^\infty((A_X,{\chi^{-1}})\backslash X_\emptyset^h) \ar@{-->}[dd] \ar[rr]^{\mathfrak M_{{\chi^{-1}}}^{-1}}  && C^\infty((A_X,{\chi^{-1}})\backslash X_\emptyset) \ar[dd]^{\mathscr S_{w,\chi}}  
\\ C_c^\infty(X) \ar[urr]^{\mathfrak N_{\chi^{-1}}}\ar[drr]_{\mathfrak N_{\chi}}&&&  \\ 
&& C^\infty((A_X,{\chi})\backslash X_\emptyset^h) \ar[rr]^{\mathfrak M_{\chi}^{-1}}  && C^\infty((A_X,{\chi})\backslash X_\emptyset),}
\end{equation}
where the notation is as follows: 
\begin{enumerate}
 \item 
The space $X_\emptyset^h$ is the space of \emph{generic horocycles} on $X$, or on $X_\emptyset$. It classifies pairs $(P,Y)$, where $P\in \B_X$ (that is, in the class of parabolics $P(X)$), with unipotent radical $U$, and $Y$ is a $U$-orbit in the open $P$-orbit of $X$, or of $X_\emptyset$; by \cite[Lemma 2.8.1]{SV}, $X$ and $X_\emptyset$ give canonically isomorphic spaces by this construction, and we will therefore also use the simplified notation $X^h$ for $X_\emptyset^h$. In the rank-one cases that we are considering, $X_\emptyset^h$ is $G$-isomorphic to $X_\emptyset$, but not canonically. Moreover, the $A_X$-action on $X_\emptyset$, defined above, gives rise to the $w$-twisted $A_X$-action on $X_\emptyset^h$ under such an isomorphism, where $w$ is the longest element of $W_X$ (in rank one, this is the inverse $A_X$-action). It is therefore best to think of $X_\emptyset^h$ as $S^+\backslash G$, where $S^+\subset P(X)$ is the kernel of the character $q$. 

 \item The operator $\mathfrak M_{\chi}$, which can be thought of as the ``standard intertwining operator'', is the operator which, in a region of convergence, takes a function in $C^\infty((A_X,{\chi})\backslash X_\emptyset)$ and integrates it over generic horocycles. \emph{Because there is no canonical measure on those horocycles, this operator depends on a choice of such measures}, and more canonically has image in the sections of a certain line bundle over $X_\emptyset^h$ (the line bundle dual to the line bundle whose fiber over a horocycle is the set of $U$-invariant measures on it --- see \cite[\S 15.2]{SV}). However, in the cases of Table \eqref{thetable} that we are interested in this paper (and, more generally, whenever $X$, hence also $X_\emptyset$, admits a $G$-invariant measure), such a choice can be made $G$-equivariantly, and it will not matter for the commutativity of the diagram --- the important point here being that horocycles in $X_\emptyset$ and $X$ are identified, and the choices of Haar measures must be made compatibly.

\item The operator $\mathfrak N_{\chi}$ is, similarly, the integral over the horocycles on $X$, followed by an averaging over horocycles in the same $B$-orbit, against the character $\chi^{-1}$ of $A_X$; that is, for a horocycle $Y$, considered both as a point in $X_\emptyset^h$ and as a subset of $X$, 
\begin{equation}\label{Nchi} \mathfrak N_\chi \Phi (Y) = \int_{A_X} \left(\int_{aY} \Phi(y) dy \right) \chi^{-1}\delta_{P(X)}^{-\frac{1}{2}}(a) da.
\end{equation}
In other words, this is the standard morphism to the principal series representation $C^\infty((A_X,\chi)\backslash X_\emptyset^h)$, given by an integral over the open Borel orbit.
Again, the Haar measure used on $A_X$ does not matter for the commutativity of the diagram.
\end{enumerate}

Note that, via the  noncanonical identification \eqref{isomprincipal}, the morphism $\mathscr S_{w,\chi}$ has to be (for almost all $\chi$) a multiple of the ``standard intertwining operator'' 
\[ \mathfrak R_\chi: I_{P(X)^-}^G(\chi^{-1}) \to I_{P(X)^-}^G(\chi)\]
given by the integral 
\[ \mathfrak R_\chi f(g) = \int_{U_{P(X)}^-} f(\tilde w ug) du,\] for some lift $\tilde w$ of $w$ to $G$.
Hence, our goal is to describe this constant of proportionality, but we must first give a careful definition of these intertwining operators, since they depend on the isomorphism \eqref{isomprincipal}, and the lift $\tilde w$, as well as the measure $u$ that appear in the definition of $\mathfrak R_\chi$. It turns out that there is a way to do define $\mathfrak R_\chi$ (which I will call ``spectral Radon transform'') that is independent of choices.

\subsection{The canonical Radon transform and the basic cases}\label{sscanonicalRadon}

\subsubsection{} The spectral Radon transforms $\mathfrak R_\chi$ will be obtained as Mellin transforms of a Radon transform $\mathfrak R: C_c^\infty(X_\emptyset)\to C^\infty(X_\emptyset)$ which, under isomorphisms $X_\emptyset \simeq S\backslash G$ as before, can be written 
\[ \mathfrak R \Phi(Sg) = \int_{U_{P(X)}^-} \Phi(S\tilde wug) du.\]
Interpreting such an integral without fixing such an isomorphism or a lift $\tilde w$ for the Weyl group element, we must describe 
a distinguished $G$-orbit $X_\emptyset^R \subset X_\emptyset\times X_\emptyset$, living over the open Bruhat cell of $\B_X\times \B_X$, such that the Radon transform 
is given by:
\begin{equation}\label{Radon-Gorbit}\mathfrak R(\Phi)(x) = \int_{(x,y)\in X_\emptyset^R} \Phi(y) dy.
\end{equation}

\subsubsection{Choice of measure for Radon transforms} \label{sssmeasureRadon}

This integral \eqref{Radon-Gorbit} also depends on fixing ($G$-equivariantly) measures $dy$ on the fibers of $X_\emptyset^R$ with respect to the first projection. The formulas for the scattering operators that we will present in this section are only true for one such choice of measure, because the spectral scattering operators $\mathscr S_{w,\chi}$ are completely canonical (and unitary on $L^2((A_X,{\chi})\backslash X_\emptyset)$, for $\chi$ unitary). The measure (and formulas) will depend on a nontrivial, unitary additive character $\psi: F\to \mathbb C^\times$, which we fix throughout, thereby fixing the corresponding self-dual Haar measure on $F$. 

It is enough to describe the measure for the standard intertwining operators for \emph{nondegenerate} principal series,\footnote{In keeping with standard conventions for boundary degenerations when $\chi$ stands for a character of the universal Cartan, i.e., the torus quotient of a Borel subgroup, we use $\chi$ to index the Radon transform whose image is in the induction of $\chi$ from the \emph{opposite} Borel subgroup. This ensures compatibility with the notation of \cite{SaTransfer1}.} 
\[\mathfrak R_{\chi,\tilde w}: I_{B^-}^G({^{w^{-1}}\chi}) \to I_{B^-}(\chi),\] defined for almost every character by meromorphic continuation of  
\begin{equation}\label{Radon-standard} \mathfrak R_{\chi,\tilde w}f(g) = \int_{(N^-\cap w^{-1} N^-w) \backslash N^-} f(\tilde w n g) dn,
\end{equation}
since the operator for degenerate principal series descends from those. This will be applied, in the course of our argument, to various spaces which are (noncanonically) quotients of $N^-\backslash G$ (i.e., horospherical), hence, being mindful of the noncanonical nature of such isomorphisms, we should think of $f$ above as a function on a space $Y$ which is isomorphic to $N^-\backslash G$, and of the representative $\tilde w$ of the Weyl element $w$ as determining a distinguished $G$-orbit in $Y\times Y$, so that $\mathfrak R_{\chi,\tilde w}$ is given by an integral analogous to \eqref{Radon-Gorbit}. Choosing a reduced decomposition $\tilde w = \tilde w_1 \cdots \tilde w_n$ into representatives for the \emph{simple} reflections, it is well-known (and immediate to check) that 
\begin{equation}\label{Radondecomp}\mathfrak R_{\tilde w} = \mathfrak R_{\tilde w_1} \circ \cdots \circ \mathfrak R_{\tilde w_n},
\end{equation}
 hence we only need to describe the measure for $w$ a simple reflection in the Weyl group. In this setting, we are reduced to the case of $\SL_2$, through the map from $\SL_2$ to the Levi of the parabolic corresponding to the simple reflection.

Hence, we take $G=\SL_2$, with a chosen $G$-orbit on $(N\backslash \SL_2)^2$. We can then fix an isomorphism between $\SL_2$ and the special linear group of a symplectic vector space $(V,\omega)$, with the chosen $G$-orbit equal to $V^R = \{ (v_1,v_2) | \omega(v_1, v_2) = 1\}$. (Any two such identifications differ by the scalar action of $\pm 1$, which does not affect Haar measures.) The canonical choice of measure, then, is the one described in \cite[\S~3.3]{SaTransfer1}, and it is as follows. By our isomorphism, the intertwining operator descends from the following Radon transform on $V$,
\begin{equation}\label{Radontransform} \mathfrak R_{\tilde w} \Phi(u) = \int_{\omega(u,v)=1} \Phi(v) dv,\end{equation}
with $dv$ the measure to be described.\footnote{The definition of Radon transform here is opposite to the convention of \cite{SaTransfer1}, where the integral was taken over the set $\omega(v,u)=1$. This change will only affect the formula \eqref{scattering-T-ext} below.} But the horocycle $\omega(u,v)=1$ (for fixed $u$), here, is canonically an $F$-torsor under $(x,v)\mapsto v+xu$, thus inheriting the Haar measure from $F$. 

This completes our description of the measure used to define Radon transforms. We complement it with a well-known formula for its inverse. It uses the \emph{gamma factors} of the local functional equation of $L$-functions:
\begin{equation}\label{FE}\gamma (\chi, s, \psi) L(\chi, s) = \epsilon(\chi, s,\psi) L(\chi^{-1}, 1-s).
\end{equation}
They are defined by the theory of local zeta integrals (Iwasawa--Tate), see \cite[2.1.4]{SaTransfer1} for a recollection. 

When $\chi$ is a character of a split torus $T$, and $r$ is a representation of the dual torus, we will use Langlands' notation for $L$- and $\gamma$-factors:
\begin{equation}\label{gammafactor}\gamma(\chi,r, s, \psi) = \prod_i \gamma(\chi\circ\lambda_i, s,\psi),
\end{equation}
where $\lambda_i$ ranges over the weights of $r$ (with multiplicities).

\begin{lemma}\label{lemmainverseRadon}
The operator $\mathfrak R_{\chi,\tilde w}$ is defined and invertible almost everywhere, with inverse 
\begin{equation}\label{inverseRadon}
 \mathfrak R_{\chi,\tilde w}^{-1} = \prod_{\alpha>0, w^{-1}\alpha<0} \gamma(\chi, \check\alpha, 0, \psi) \gamma(\chi, -\check\alpha, 0, \psi^{-1}) \mathfrak R_{{^{w^{-1}}\chi},\tilde w^{-1}}.
\end{equation}
\end{lemma}

\begin{proof}
 This reduces again, by the same argument, to $\SL_2$, where the formula can be proven by relating Radon transforms to Fourier transforms, as before \cite[(3.25)]{SaTransfer1}. However, there is a sign error in \emph{op.cit.}: The assertion that $\mathfrak F^*\circ \mathfrak F=1$ at the bottom of p.50 of \emph{op.cit}.\ is incorrect. The correct statement is that $\mathfrak F\circ \mathfrak F=1$, which implies that the characters $\psi^{-1}$ in (3.25) and (3.28) of \emph{op.cit}.\ should be $\psi$. However, the Radon transform of (3.25) corresponds to what we denote here by $\mathfrak R_{\chi, \tilde w}$ for a representative of the nontrivial Weyl element of $\SL_2$ with $\tilde w^2=-I$; this is $\chi(-I)\mathfrak R_{\chi, \tilde w^{-1}}$, and  using $\chi(-I) \gamma(\chi, -\check\alpha, 0, \psi) = \gamma(\chi, -\check\alpha, 0, \psi^{-1})$, we obtain our formula.
\end{proof}

\begin{remark}
 A technical detail: When $w^2=1$, the operator $\mathfrak R_{{^{w^{-1}}\chi},\tilde w^{-1}}$ is not necessarily equal to the operator $\mathfrak R_{{^{w^{-1}}\chi},\tilde w}$ -- which is why we stress the choice of representative $\tilde w$ in the notation. For example, for $\SL_2$, we cannot find $\tilde w$ with $\tilde w^2=1$. More abstractly, in terms of the chosen $G$-orbit on $Y\times Y$ (where $Y\simeq N^-\backslash G$), this orbit does not necessarily get preserved by switching the two copies of $Y$. However, this will be the case for the cases of Table \ref{thetable} -- even when $X=\SL_2$, the group acting will be $(\SL_2\times\SL_2)/\mu_2^\diag=\SO_4$, where the longest Weyl element is represented by an involution in the group.
\end{remark}

\subsubsection{Review of the basic cases} \label{sssreviewbasic}

Having fixed the measure for the Radon transform \eqref{Radon-Gorbit}, there remains to describe the ``canonical'' $G$-orbit $X_\emptyset^R \subset X_\emptyset\times X_\emptyset$.
Such $G$-orbits were described in \cite[\S 3]{SaTransfer1} for the basic cases of the two families of Table \eqref{thetable}, denoted by $A_1$ and $D_2$: $X=\SO_2\backslash\SO_3 =\Gm\backslash\PGL_2$ and $X=\SO_3\backslash \SO_4\simeq \SL_2$.
The word ``canonical'' is too strong, since there is nothing compelling about choosing these $G$-orbits over others. However, they are independent of choices of isomorphisms, and eventually their importance is that they are useful in obtaining exact formulas for the scattering operators.

\begin{remark}
 The varieties of the second group of Table \eqref{thetable} admit various forms over the field $F$ with $G$ and $H$ split, parametrized by $F^\times/(F^\times)^2$. Namely, since they are all isomorphic (as abstract varieties) to $\SO_{2n-1}\backslash \SO_{2n}$, the forms depend on the discriminant of the orthogonal complement of a $(2n-1)$-dimensional, maximally isotropic quadratic space $V$ inside of a split $(2n)$-dimensional quadratic space $U$.  We \emph{choose} to compute only for the cases where this discriminant is (square equivalent to) $1$, of which $X=\SL_2$ is the base case; minor modifications are needed to accommodate other discriminants. The same choice was made for the calculation of transfer operators in \cite{SaRankone}, but unfortunately I did not explicitly state that choice there.
\end{remark}

\subsubsection{} \label{ssstorus}

In the case of $X\simeq \Gm\backslash\PGL_2$, we should fix a $2$-dimensional symplectic vector space $V$, an isomorphism $G\simeq \PGL(V)$, and an isomorphism of $X\simeq \SO_2\backslash\SL_2$ with the space of quadratic forms of the form $q=xy$ for some choice of standard symplectic coordinates (i.e., such that the symplectic form  is $dx\wedge dy$).   Once the isomorphism $G\simeq \PGL(V)$ is fixed, such an identification for $X$ is unique up to the $G$-automorphism group $\mathbb Z/2$ of $X$, which sends the quadratic form $xy$ to the quadratic form $-xy$ (corresponding to the standard coordinates $(-y,x)$). The variety $X_\emptyset$ can, then, be identified with the variety of rank-one degenerate quadratic forms on $V$, which is the same as $V^{*\times}/\{\pm 1\}$ (where $V^*$ is the dual, and the exponent $~^\times$ denotes the complement of zero), since those quadratic forms, over the algebraic closure, are squares of linear forms. Via the identification $V\simeq V^*$ afforded by the symplectic form, we can also identify $X_\emptyset$ with $V^\times/\{\pm 1\}$. Evaluation of the quadratic form:
\[\ev:  X_\emptyset \times V \to \Ga\]
on $\Ga$ gives rise to our ``canonical'' $G$-orbit $X_\emptyset^R$, as the image of $\ev^{-1}(1)$ under $X_\emptyset\times V^\times \to X_\emptyset \times X_\emptyset$. Applying the automorphism $xy\mapsto -xy$ to $X$ also induces the analogous automorphism on $X_\emptyset$, fixing the distinguished orbit $X_\emptyset^R\subset X_\emptyset^2$. 

With this distinguished orbit fixed, the scattering operator 
\[ \mathscr S_{w,\chi}: C^\infty((A_X,\chi^{-1})\backslash X_\emptyset)\to C^\infty((A_X,\chi)\backslash X_\emptyset)\]
was calculated in \cite[Theorem 3.5.1]{SaTransfer1}
\begin{equation}\label{scattering-T}\mathscr S_{w,\chi} = \gamma(\chi, \frac{\check\alpha}{2}, \frac{1}{2}, \psi^{-1}) \gamma(\chi, \frac{\check\alpha}{2}, \frac{1}{2}, \psi) \gamma(\chi,-\check\alpha, 0,\psi) \cdot \mathfrak R_{\chi}.\end{equation}

Here, $\psi$ is any nontrivial unitary character of the additive group $F$, and the Radon transforms are computed with respect to
a measure that is directly proportional to the self-dual Haar measure on $F$ with respect to $\psi$. Note that the change-of-character formula for $\gamma$-factors,
\begin{equation}\label{changeofpsi}
 \gamma(\chi,s,\psi(a\bullet)) = \chi(a)|a|^{s-\frac{1}{2}}\gamma(\chi,s,\psi(\bullet))
\end{equation}
(see \cite[3.3.3]{Deligne-epsilon}, taking into account that the self-dual measure for $\psi(a\bullet)$ is $|a|^\frac{1}{2}$ the self-dual measure for $\psi$) implies that the factor 
\[\gamma(\chi, \frac{\check\alpha}{2}, \frac{1}{2}, \psi^{-1}) \gamma(\chi, \frac{\check\alpha}{2}, \frac{1}{2}, \psi) \gamma(\chi,-\check\alpha, 0,\psi)\] 
will get multiplied by $|a|^{-\frac{1}{2}}$, which is the inverse of the factor by which the measure defining $\mathfrak R_\chi$ will be multiplied, making the above expression for $\mathscr S_{w,\chi}$ independent of $\psi$.

\subsubsection{} \label{sssgroup}

In the case of $X=\SL_2=\SL(V)$, on the other hand, we can identify $X_\emptyset$ with the space of endomorphisms of $V$ of rank $1$, i.e., operators of the form $\tau: V/L\to L'$, where $L$ and $L'$ are lines in $V$. 
The distinguished $G^\diag$-orbit, in that case, is the set of pairs $(\tau_1,\tau_2)\in X_\emptyset\times X_\emptyset$ with $\tau_1+\tau_2\in \SL(V)$. With this distinguished orbit at hand, we define the spectral Radon transforms $\mathfrak R_\chi$, and a correction\footnote{See the proof of Lemma \ref{lemmainverseRadon}.}  of \cite[Theorem 3.5.1]{SaTransfer1} in this case gives the following formula for the  scattering operator:
\begin{equation}\label{scattering-G}
   \mathscr S_{w,\chi}= \gamma(\chi,\check\alpha, 0,\psi) \gamma(\chi,-\check\alpha,0,\psi)  \cdot \mathfrak R_\chi.
\end{equation}
Again, the product of gamma factors would be multiplied by $|a|^{-1}$, if we changed $\psi$ to $\psi(a\bullet)$, which is exactly inverse the factor by which the measure defining $\mathfrak R_\chi$ would change. (Notice here that the intertwining operator $\mathfrak R_\chi$ is given by an integral over a $2$-dimensional unipotent subgroup.)

\subsubsection{} \label{ssstorus-ext}

We will also need a twisted version of the scattering operator for $\Gm\backslash\PGL_2$ (\S~\ref{ssstorus}), where instead of functions we consider the induction of a character of $\Gm(F)$; equivalently, we can generalize the calculation of scattering operators to the variety $X=\Gm\backslash\GL_2$, where $\Gm$ is embedded as the subgroup $\begin{pmatrix} a \\ & 1\end{pmatrix}$. This variety can be identified with $\SL_2$, with the action of $\SO_4 \simeq \SL_2\times^{\mu_2}\SL_2$ restricted to a Levi subgroup isomorphic to $G=\SL_2\times^{\mu_2} \Gm\simeq\GL_2$, but we prefer to think, again, of $G$ as the general linear group of a $2$-dimensional symplectic vector space $V$, and of $X$ as the variety parametrizing standard bases $(v_1, v_2)$ (i.e., $\omega(v_1, v_2)=1$), with $G$-action 
\[ (v_1, v_2)\cdot g = ((\det g)^{-1} v_1 g, v_2 g).\]

In this case the torus $A_X$ is the full Cartan of $\GL_2$. The boundary degeneration $X_\emptyset$ can then be identified with the space of pairs $(v_1, v_2)$ of nonzero vectors with $\omega(v_1, v_2) = 0$ (i.e., colinear), with similar $G$-action, and there is a map $X_\emptyset\to\Gm$ taking such a pair to the constant $c$ such that $v_1=c v_2$. The fibers of this map can be identified with the nonzero vectors in $V$ via the projection to $v_2$, and for simplicity we will fix a Haar measure on $\Gm$ and fix the $G$-invariant measure on $X_\emptyset$ whose disintegration with respect to that specializes to the measure corresponding to the symplectic form on $V_1^\times=V^\times$, the fiber over $1\in\Gm$. (This choice will not affect the functional equations.)

We fix the distinguished $G$-orbits\footnote{Note that the distinguished $G$-orbit $X_\emptyset^R$ is not involutive in this case, and therefore it matters for the formulas in this paragraph that we define Radon transform with respect to the first, instead of the second projection. This will not matter, however, for the application of this case to the proof of Theorem \ref{thmscattering}, since we will only apply the formula for the scattering operator to even characters. Moreover, note that the intersection of $X_\emptyset^R$ with $V_1^\times\times V_1^\times$ coincides with the distinguished $\SL(V)$-orbit on $V\times V$ described in \cite[\S~3.3]{SaTransfer1}, hence equation (3.16) of \emph{op.cit}., relating Fourier transforms and Radon transforms, applies to this setting -- up to a change due to our conventions, explained in Footnote \ref{conventionRadon}.}
 $X^R\subset X \times X_\emptyset$ and $X_\emptyset^R\subset X_\emptyset\times X_\emptyset$ consisting of those quadruples $((v_1, v_2),(w_1,w_2))$ such that 
 \[\omega(v_1, w_2)=1= \omega(v_2,w_1).\] 
 The distinguished $G$-orbit $X^R$ identifies $X_\emptyset$ with the generic horocycle space $X^h$, taking the point $y\in X_\emptyset$ to the horocycle $\{x\in X| (x,y)\in X^R\}$; similarly, the $G$-orbit $X_\emptyset^R$ identifies $X_\emptyset$ with the generic horocycle space $X_\emptyset^h$ (intertwining the action of $a\in A_X$ with that of ${^wa}$, where $w$ is the nontrivial element of the Weyl group). These $G$-orbits correspond to each other, in the sense that these identifications are ``the same,'' i.e., compatible with the canonical identification $X^h=X_\emptyset^h$.\footnote{\label{footnotecompatibility}The reader who wishes to understand this should refer to the construction of the identification $X^h=X_\emptyset^h$ in \cite[\S~2.8]{SV}. A way to reformulate that identification is the following: On the affine degeneration $\mathcal X \to \mathbb A^1$ recalled in \S~\ref{sssaffinedegen}, which here we can take to be the family of pairs $(v_1, v_2) \in V^2$, with image $\omega(v_1,v_2)\in \mathbb A^1$, there is a $G$-commuting
  action of $A_X$ which matches its action on the special fiber $\mathcal X_0=X_\emptyset$. Here, if we identify $A_X = \Gm^2$ as the diagonal torus with respect to the upper triangular Borel subgroup, the action on the family is $\diag(a,b)\cdot (v_1,v_2) = (b^{-1}v_1,av_2)$. The canonical isomorphism $X^h\xrightarrow{\sim} X^h_\emptyset$ takes a generic horocycle $x\cdot N$ on $X$ to $\lim_{t\to 0} e^{\check\alpha}(t)\cdot (x\cdot N) \cdot e^{-\check\alpha}(t)$, where the action on the left is the action of $A_X$ on $\mathcal X$, and the action on the right is the action of $A_X$ on $X^h$, or equivalently the action of any section of the quotient $B\to A_X$ on $X$.} 
 
Using these orbits to define Radon transforms, as before, we have the following formula for the scattering operator $\mathscr S_{w,\chi}$.
\begin{equation}\label{scattering-T-ext}\mathscr S_{w,\chi} = \gamma(\chi, \check\epsilon_1,\frac{1}{2},\psi^{-1})\gamma(\chi,-\check\epsilon_2,\frac{1}{2},\psi) \gamma(\chi,-\check\alpha, 0,\psi) \cdot \mathfrak R_{\chi},\end{equation}
where $\check\epsilon_1$, $\check\epsilon_2$ are the standard coweights of $\GL_2$ (i.e., the weights of the standard representation of its dual). Note that, when the character is trivial on the center, the formula specializes to \eqref{scattering-T}. 

I sketch the proof, which is essentially the same as the proof of the special case \eqref{scattering-T} in \cite{SaTransfer1}: 

The adjoint of the operator $\mathfrak N_\chi$ of \eqref{fiberscatteringdiagram} is an operator 
\[C_c^\infty(X_\emptyset^h) \to C^\infty(X)\]
which factors through the dual of the principal series representation 
\[C^\infty((A_X,{\chi})\backslash X_\emptyset^h)\simeq I_{P(X)}(\chi).\] 
I claim that, when composed with evaluation at $(v_1, v_2)\in X$, this operator is given by the meromorphic continuation of the integral against an $A_X$-eigenfunction, whose restriction to $V_1^\times \subset X_\emptyset^h\simeq X_\emptyset$ is of the form 
\begin{equation}\label{eq:distribution} -a^{-1} v_1 + b v_2 \mapsto \chi_1(a) \chi_2(b),\end{equation}
for appropriate characters (to be described below) $\chi_1$, $\chi_2$. 
Indeed, it is easy to see that this is the integral over the preimage $\{(w_1,w_2)| \omega(v_1, w_2)=1= \omega(v_2,w_1)\}$ of $(v_1,v_2)$ in $X^R$, averaged against the appropriate character over the $A_X$-action.

The characters $\chi_1$, $\chi_2$ are computed by the equivariance properties of this operator, taking into account that, when $X_\emptyset$ is thought of as $X_\emptyset^h$, the $A_X$-action is given by $\diag(a,b)\cdot (w_1,w_2)= (a^{-1}w_1, b w_2)$ (with the conventions of Footnote \ref{footnotecompatibility}). One computes that $\chi_1 (z)= |z|^{\frac{1}{2}} \chi\circ e^{\check\epsilon_1}(z)$ and $\chi_2 (z)= |z|^{-\frac{1}{2}} \chi\circ e^{\check\epsilon_2}(z)$
(Compare\footnote{The minus sign in front of the coefficient $a^{-1}$, which does not appear in the cited formula, is because of the way that we have embedded $V_1^\times$ into $X_\emptyset$. When we mod out by the center $Z$ of $G$, we get $X_\emptyset/Z\simeq$ the space of rank-one quadratic forms $q$ on $V$, but here an element $v\in V_1^\times$ will map to $q=-(\omega(v,\bullet))^2$, while in \emph{op.cit}.\ it was mapping to the negative of that. In the notation of p.45 of \emph{op.cit}., our space $V_1$ corresponds to a twist $V^\alpha$ of $V$ corresponding to the quadratic extension $F(\sqrt{-1})$.} with the formula $\int \Phi(y,z) \tilde\chi(yz) dy dz$ on p.47 of \cite{SaTransfer1}, which is for the adjoint of the operator $\mathfrak N_{{^w\chi}}$, when the central character is trivial.) 

The integral against \eqref{eq:distribution} is a product of two Tate integrals, and the argument of \emph{op.cit}.\ generalizes to this case to show that, with $\mathfrak F^*$ denoting the $G$-equivariant extension of 
the symplectic Fourier transform
\[\mathfrak F^*\Phi(v) = \int_V \Phi(w) \psi^{-1}(\omega(w,v)) dw,\]  on $V_1$, we have 
\[ \mathfrak N_\chi^* \circ \mathfrak F^* = \gamma(\chi_1, 0, \psi^{-1})\gamma(\chi_2^{-1}, 0, \psi) \mathfrak N_{{^w\chi}} = \gamma(\chi, \check\epsilon_1,\frac{1}{2},\psi)\gamma(\chi,-\check\epsilon_2,\frac{1}{2},\psi^{-1}) \mathfrak N_{{^w\chi}}^*,\]
which as in \emph{op.cit}.\footnote{\label{conventionRadon}Since our convention for the definition of Radon transform in \eqref{Radontransform} is opposite to the one of \emph{op.cit}., in applying formula (3.16) of that reference one needs to change $\psi$ to $\psi^{-1}$. We then apply the identity $\gamma(\chi, \check\epsilon_1,\frac{1}{2},\psi)\gamma(\chi,-\check\epsilon_2,\frac{1}{2},\psi^{-1})\gamma(\chi,-\check\alpha, 0,\psi^{-1}) = \gamma(\chi, \check\epsilon_1,\frac{1}{2},\psi^{-1})\gamma(\chi,-\check\epsilon_2,\frac{1}{2},\psi) \gamma(\chi,-\check\alpha, 0,\psi)$ to obtain \eqref{Radontransform}.} implies \eqref{scattering-T-ext}.

\subsubsection{The other cases of Table \eqref{thetable}}

The calculation of scattering operators for the general case of a variety of Table \eqref{thetable} will be reduced to the basic cases by means of the following lemma. 
For every subset $I$ of the simple roots of $G$, we denote by $P_I$ the corresponding class of parabolics (or a representative), by $U_{P_I}$ its unipotent radical, and by $L_I$ its Levi quotient. 

\begin{lemma}\label{lemmaclosedorbits}
Let $Z\subset X$ be a closed $B$-orbit. 
\begin{itemize}
 \item For the cases of Table \eqref{thetable} with $\check G_X=\SL_2$, $Z$ is of rank zero, and for every simple root $\alpha$ such that $Y:=ZP_\alpha\ne Z$, $Y_2:=Y/U_{P_\alpha}$ is $L$-isomorphic to $\SO_2\backslash\SL_2 \simeq \Gm\backslash \PGL_2$ under some morphism $L\to \PGL_2$.
 \item For the cases of Table \eqref{thetable} with $\check G_X=\PGL_2$, and for all simple roots $\alpha$ we have $ZP_\alpha=Z$, except for two orthogonal simple roots $\alpha, \beta$ for which, setting $Y:=ZP_{\alpha\beta}$, $Y_2:=Y/U_{P_{\alpha\beta}}$ is $L$-isomorphic to $\SO_3\backslash\SO_4$ under some homomorphism $L\to \SO_4$.
\end{itemize}
\end{lemma}

As a matter of notation, what we denote by $Y/U_P$ here is the \emph{geometric} quotient of $Y$ by the $U_P$-action, not the stack quotient, which could be different because of nontrivial stabilizers.

\begin{proof}
 This is a combination of Lemmas 2.2.4 and 2.3.5 of \cite{SaRankone}.  
\end{proof}

In the remainder of this subsection, we will describe a ``canonical'' $G$-orbit $X_\emptyset^R \subset X_\emptyset\times X_\emptyset$, living over the open Bruhat cell in $\B_X\times \B_X$. This will give rise to a ``canonical'' Radon transform, which we will use in the next section to describe the scattering operators.

As mentioned before, a $G$-orbit on $X_\emptyset\times X_\emptyset$ over the open Bruhat cell in $\B_X\times \B_X$ is equivalent to a $G$-equivariant isomorphism $\iota: X_\emptyset \xrightarrow\sim X_\emptyset^h$, where $X_\emptyset^h$ is the space of generic horocycles on $X_\emptyset^h$, as above. Indeed, such an isomorphism defines the distinguished $G$-orbit 
\[X_\emptyset^R=\{(x,y) \in X_\emptyset\times X_\emptyset| x\in \iota(y)\},\]
and, vice versa, can be recovered from it. Note, however, that this isomorphism intertwines the action of $a\in A_X$ with the action of ${^wa}$, where $w\in W_X$ is the longest element. Similarly, this is equivalent to describing a distinguished $G$-orbit 
\[ X_\emptyset^{R,h} \subset X_\emptyset^h\times X_\emptyset^h,\]
which by the canonical isomorphism $X_\emptyset^h\simeq X^h$ of \cite[Lemma 2.8.1]{SV}, can be understood as a distinguished $G$-orbit $X^{R,h}\subset X^h\times X^h$. 
This is the orbit that we will describe. 

The orbit (and the Radon transform) will depend, a priori, on the choice of a closed $B$-orbit $Z$, as in Lemma \ref{lemmaclosedorbits}; a posteriori, by the calculation of scattering operators, it doesn't, a fact that can also be proved directly (but we will not do). Hence, fix $Z$, let $Y$ be as in Lemma \ref{lemmaclosedorbits}, and let $P$ be the parabolic $P_\alpha$, resp.\ $P_{\alpha\beta}$, appearing in the two cases of the lemma. Let $Y^\circ\subset Y$ be the open $B$-orbit; without fixing a Borel subgroup, we can consider $Y^\circ$ as a $G$-orbit on $X\times \B$ --- denoted by $\tilde Y^\circ$ to avoid confusion. Similarly to the definition of $X^h$, we can define 
\[ X^{h,Y} = \{(B, M)| B\in \B, \, M \mbox{ is a $U_B$-orbit with } (m,B)\in \tilde Y^\circ\mbox{ for any }m\in M\}.\]
Hence, these are not generic horocycles, but horocycles corresponding to the $B$-orbit $Y^\circ$.

Let $\tilde X^h$ be the base change $X^h\times_{\B_X}\B$ of $X^h$ to the full flag variety $\B$. Fixing a Borel subgroup $B \in \B$, we have noncanonical isomorphisms:
\begin{equation}\label{noncan1} \tilde X^h = A_X\times^B G,
\end{equation}
\begin{equation}\label{noncan2} X^{h,Y}= A_Y\times^B G,
\end{equation}
where $A_Y$ is the torus quotient by which $B$ acts on the geometric quotient $Y^\circ/U_B$. 

Let $w'$ be the Weyl group element $w_\alpha$ or $w_\alpha w_\beta$, respectively (where $w_\alpha$, $w_\beta$ denote simple reflections), for each of the two cases of Lemma \ref{lemmaclosedorbits}. A result of Knop \cite[\S~6]{KnOrbits} implies that there is an element $w_1$ of the Weyl group of $G$, such that \begin{equation}\label{opensmalliso}
Y^\circ \times^{N \cap  w_1 Nw_1^{-1}} \tilde w_1N\xrightarrow\sim X^\circ
\end{equation}
under the restriction of the action map to $\tilde w_1N\subset G$; here, $\tilde w_1$ is any lift of $w_1$, thought of as a double coset of $B\backslash G/B$, to $G$. Note that this implies that $\codim Y^\circ = \operatorname{length}(w_1)$, and $A_Y^{w_1} = A_X$.

Let $w_1\in W$ be such an element. It is known \cite{Brion}, \cite[\S~6.2]{SaSph} that the nontrivial element $w \in W_X$ is equal to $w_1^{-1} w' w_1$.

\begin{lemma}\label{lemmalength}
 In the setting above, the decomposition $w=w_1^{-1} w' w_1$ is reduced, i.e., we have $ \operatorname{length}(w) = 2  \operatorname{length}(w_1) +  \operatorname{length}(w')$.
\end{lemma}

\begin{proof}
One can make an abstract argument for that (for example, in the cases with $\check G_X=\SL_2$, $\gamma$ is a root of $G$, and any reduced decomposition of $w=w_\gamma$ has to correspond to a ``Brion path'' by considerations related to ranks of orbits), but it is also straightforward to see this by inspection of Table \eqref{thetable}: 

 {\small
$ \begin{array}{|c|c|c|c|c|}
\hline
&X = H\backslash G & \dim X & \dim Z & \operatorname{length}(w) \\
 & & & & \\
\hline
\hline
A_1& \Gm \backslash\PGL_{2} & 2 & 1 & 1 \\
\hline
A_n& \GL_n \backslash\PGL_{n+1} & 2n & n  & 2n-1\\
\hline
B_n& \SO_{2n}\backslash\SO_{2n+1} & 2n & n & 2n-1 \\
\hline
C_n& \Sp_{2n-2}\times \Sp_{2}\backslash \Sp_{2n} & 4n - 4 & 2n-2 & 4n-5 \\ 
\hline
F_4& \Spin_9\backslash F_4 & 16 & 8 & 15 \\
\hline
G_2& \SL_3\backslash G_2& 6 & 3 & 5 \\
\hline
\hline
D_2& \SL_2=\SO_{3}\backslash\SO_{4} & 3 & 2 & 2 \\
\hline
D_n& \SO_{2n-1}\backslash\SO_{2n} & 2n-1 & n & 2n-2  \\
\hline
D_4'' & \Spin_7\backslash\Spin_8 & 7 & 4 & 6\\
\hline
B_3''& G_2\backslash\Spin_7 & 7 & 4 & 6\\
\hline
\end{array}
$
}

Here, for $X=H\backslash G$, we computed $\dim Z$, the dimension of a closed Borel orbit, as $\dim B_G - \dim B_H  $, where by $B_*$ we denote the corresponding Borel subgroups. 
\end{proof}

By Lemma \ref{lemmaclosedorbits},
the $P$-variety $Y_2:=Y/U_P$ is isomorphic either to $\Gm\backslash\PGL_2$ or to $\SL_2$ under the action of the Levi $L$ of $P$. In each of the two cases, a distinguished $L$-orbit on $Y_2^h\times Y_2^h$, living over the open Bruhat cell, was described in \S~\ref{sssreviewbasic} above. This corresponds to a $G$-orbit on $X^{h,Y}\times X^{h,Y}$, living over the Bruhat cell corresponding to $w'$. We denote by $X^{R,h,Y}$ this $G$-orbit.

Now, the product $X^{h,Y}\times \tilde X^h$ lives over the product $\B \times \B$. It follows from \eqref{opensmalliso} that there is a \emph{distinguished} $G$-orbit $X'^h\subset X^{h,Y}\times \tilde X^h$ that lives over the Bruhat cell corresponding to $w_1$, that is, over the $G$-orbit of a pairs $(B,B'=\tilde w_1B\tilde w_1^{-1})$, where, as above, $\tilde w_1$ is any lift of $w_1\in B\backslash G/B$ to $G$. The fiber of this orbit over $(B,B')$ consists of all pairs of horocycles $((B, M) , (B', x N'))$ with $N'\subset B'$ the unipotent radical and $x\in M$.  Equation \eqref{opensmalliso} implies that if the horocycle $(B,M)$ belongs to $X^{h,Y}$, the horocycle $(B',xN')$ is generic, i.e., belongs to $\tilde X^h$.  

Consider the set of quadruples 
\[(x_1,y_1,y_2,x_2)\]
with $x_1, x_2\in \tilde X^h$, $(y_1,y_2)\in X^{R,h,Y}\subset (X^{h,Y})^2$, and $(y_i, x_i)\in X'^h$ for $i=1,2$. It lives over the set of quadruples 
\[ (B, \tilde w_1^{-1} B \tilde w_1,  \tilde w_1^{-1}\tilde w' B (\tilde w')^{-1}\tilde w_1 , \tilde w_1^{-1} \tilde w' \tilde w_1' B (\tilde w_1')^{-1} (\tilde w')^{-1}\tilde w_1)\]
of Borel subgroups, where the tilde denotes, again, lifts to $G$, and $\tilde w_1$, $\tilde w_1'$ are two possibly different lifts of $w_1$. However, because of Lemma \ref{lemmalength}, the element $\tilde w = \tilde w_1^{-1} \tilde w' \tilde w_1'$ is a lift of $w\in W_X$, and the quadruples of Borel subgroups as above form a single $G$-orbit.
\emph{The distinguished $G$-orbit $X^{R,h}\subset X^h\times X^h$, now, is the image of this set under the first and last projections, composed with the projection $\tilde X^h\to X^h$.}
The reader can immediately check that this is indeed a $G$-orbit, using the noncanonical isomorphisms \eqref{noncan1}, \eqref{noncan2}.

By the discussion above, this $G$-orbit corresponds to a distinguished $G$-orbit $X_\emptyset^R\subset X_\emptyset\times X_\emptyset$, that we use to define the Radon transform $\mathfrak R: C_c^\infty(X_\emptyset)\to C^\infty(X_\emptyset)$ by \eqref{Radon-Gorbit}. The spectral Radon transforms $\mathfrak R_\chi$ descend from it by Mellin transform. For the choice of  measure used to define $\mathfrak R$, see \S~\ref{sssmeasureRadon}.

\subsection{Formula for the scattering operators}

\subsubsection{}
The main result of this section is the following, where we fix a unitary additive character $\psi$:

\begin{theorem}\label{thmscattering}
 For the cases of Table \eqref{thetable}, in terms of the spectral Radon transforms 
 \[\mathfrak R_\chi: C^\infty((A_X,\chi^{-1})\backslash X_\emptyset) \to C^\infty((A_X,\chi)\backslash X_\emptyset)\]
 that descend from the canonical Radon $\mathfrak R$ described in \S~\ref{sscanonicalRadon},
 the scattering operator $\mathscr S_{w,\chi}$ for the nontrivial element $w\in W_X$ is given by 
 \begin{equation}
  \mathscr S_{w,\chi} = \mu_X(\chi)  \cdot \mathfrak R_{\chi},
 \end{equation}
where $\mu_X$ is given by the following formulas:
 \begin{itemize}
 \item  for the cases with $\check G_X=\SL_2$, with $L_X = L(\Std, s_1)L(\Std, s_2)$, 
\begin{equation} \label{scattering-typeT}
	\mu_X(\chi) = \gamma(\chi, \frac{\check\gamma}{2}, 1-s_1, \psi^{-1}) \gamma(\chi, \frac{\check\gamma}{2}, 1-s_2, \psi) \gamma(\chi,-\check\gamma, 0,\psi),\end{equation}
 \item for the cases with $\check G_X=\PGL_2$, with $L_X= L(\Ad, s_0)$,  
 \begin{equation} \label{scattering-typeG}
  \mu_X(\chi) = \gamma(\chi,\check\gamma, 1-s_0 ,\psi) \gamma(\chi,-\check\gamma,0,\psi).
\end{equation}
\end{itemize}
\end{theorem}

Here, $\gamma(\chi, \check\lambda, s, \psi)$ denotes the gamma factor \eqref{gammafactor} of the local functional equation for the abelian $L$-function associated to the composition of $\chi$ with the cocharacter $\check\lambda: \Gm\to A_X$. Notice that, by \eqref{changeofpsi}, if we replace $\psi$ by $\psi(a\bullet)$, for some $a\in F^\times$, the factor $\mu_X(\chi)$ changes by a factor of $|a|^{-s}$, where $s=s_1+s_2 - \frac{1}{2}$, in the first case, and $s=s_0$, in the second case. It so happens (see \cite[\S~1.2]{SaRankone}) that $2s = \dim X - 1 = $ the dimension of the unipotent radical of $P(X)$. Therefore, the measure used to define the Radon transform $\mathfrak R$, which is proportional to the self-dual measure with respect to $\psi$ to the power $\dim_{U_{P(X)}}$, changes by $|a|^s$, making the formula above for $\mathscr S_{w,\chi}$ independent of $\psi$. 

\subsubsection{}
The proof of the theorem will be given in a somewhat telegraphic fashion, because the arguments are essentially the same as the ones used to compute ``functional equations'' in \cite[Section 6]{SaSph}. The reader who wishes to read a detailed and explicit account of the arguments that follow is advised to look at that reference. The ``added value'' of the present work consists in the following:

\begin{itemize}
 \item We adopt the formalism of scattering operators, introduced in \cite{SV}. This adds an extra layer of complication; for example, \cite{SaSph} only considered the functional equations represented by the first vertical (dotted) arrow of \eqref{fiberscatteringdiagram}, while the scattering operators are given by the second vertical arrow. The benefit, however, is that the scattering operators are directly related to the asymptotics (and other constructions such as the ``unitary asymptotics'' of the Plancherel formula), a fact that we will use in the next section.
 \item We allow for general principal series representations, not only unramified ones. This is quite straightforward and does not change any of the arguments of \cite{SaSph}, once the ``basic cases'' of \S~\ref{sssreviewbasic} have been computed.
 \item We ``rigidify'' certain constructions to eliminate ambiguities up to automorphism groups. In the unramified case of \cite{SaSph}, integral structures provided such rigidifications, but, to consider ramified characters, this is not enough. Therefore, we replace various explicit integrals and principal series representations by (noncanonically) isomorphic constructions that live over the various horocycle spaces already introduced. This necessitated the pedantic discussion of the canonical Radon transforms of \S~\ref{sscanonicalRadon}, but is a more conceptual description of these constructions.
\end{itemize}

\begin{proof}[Proof of Theorem \ref{thmscattering}]
 The proof will follow the argument of \cite[\S~6.5]{SaSph}. Because of the isomorphism \eqref{opensmalliso}, the morphism $\mathfrak N_\chi$ of Diagram  \eqref{fiberscatteringdiagram}, given by an integral over the open $P(X)$-orbit $X^\circ$, can be expressed in terms of a similar integral over the smaller $B$-orbit $Y^\circ$. Namely, choosing a lift $\tilde w_1$ of $w_1$, we get a bijection $M\mapsto M \tilde w_1 N \tilde w_1^{-1}$ between horocycles in $Y^\circ$ for a given Borel subgroup $B$ and generic horocycles for the Borel subgroup $B'=\tilde w_1 B \tilde w_1^{-1}$.   
 Then, for $\Phi\in C_c^\infty(X)$, we can rewrite \eqref{Nchi} as 
 \begin{multline}\label{factorizeintegral} \mathfrak N_\chi \Phi(M \tilde w_1 N \tilde w_1^{-1}) = \int_{A_X} \int_{M\tilde w_1 N \tilde w_1^{-1}}   \Phi(x a ) dx \, \chi^{-1}\delta_{P(X)}^{-\frac{1}{2}}(a) da  = \\ \int_{(N\cap \tilde w_1^{-1} N \tilde w_1)\backslash N} \int_{A_Y} \int_{M}  \Phi(y a \tilde w_1 n \tilde w_1^{-1}) dy \, \cdot \\ {^{w_1}\left(\chi^{-1}\cdot \delta_{P(X)}^{-\frac{1}{2}}\right)}\cdot \delta_{(N\cap\tilde w_1^{-1} N \tilde w_1)\backslash N} (a) da \, dn,
 \end{multline}
where we have factored the invariant measure $dx$ on the horocycle $M\tilde w_1 N \tilde w_1^{-1}$ in terms of an invariant measure $dy$ on the horocycle $M$ and an invariant measure on $(N\cap \tilde w_1^{-1} N \tilde w_1)\backslash N$. (Recall that the choice of this measure is not canonical, but is canceled by the compatible choice of measure for the operator $\mathfrak M_{{\chi}}$ of Diagram \eqref{fiberscatteringdiagram}.)

\begin{remark}
We use the ``universal'' Cartan $A$ and its quotients $A_X$ and $A_Y$ in the notation above, with $\chi$, $\delta_{P(X)}$, etc.\ considered as characters of these universal quotients, and with $\delta_{(N\cap\tilde w_1^{-1} N \tilde w_1)\backslash N}$ considered as a character of $A$ via the quotient $B\to A$. For example, in the first integral, when $a\in A_X$ is represented by an element $t\in B'$, the character $\delta_{P(X)}$ really stands for the modular character of the representative $P(X)'$ of the class of $P(X)$ which contains $B'$. To prevent some natural confusion, I note here that the elements $a$ of the two integrals above can be represented by the \emph{same} element $t$ of a common torus $T\subset B\cap B'$, but the images of $t$ in the ``universal'' Cartan quotients $B\to A\to A_Y$ and $B'\to A\to A_X$ differ by the right action of $w_1$ on $A$, hence the character $T\to B'\to  A_X \to \CC^\times$, where the last arrow is $\chi^{-1}\delta_{P(X)}^{-\frac{1}{2}}$, is equal to the character $T\to B\to A_Y\to \CC^\times$, where the last arrow is ${^{w_1}\left(\chi^{-1}\cdot \delta_{P(X)}^{-\frac{1}{2}}\right)}$.
\end{remark}

 The inner two integrals of \eqref{factorizeintegral} represent a morphism that we can denote  
\[ \mathfrak N_{\chi}^Y: C_c^\infty(X)\to C^\infty((A_Y, {^{w_1}\chi})\backslash X^{h,Y})'.\]

There is some abuse of notation in the expressions above, due to the fact that the characters $^{w_1}\chi$, $\delta_{(N\cap\tilde w_1^{-1} N \tilde w_1)\backslash N}$ are characters of the universal Cartan $A$ of $G$ which do not, in general, factor through $A_Y$. Related to this is that the measures on the horocycles corresponding to the points of $X^{h,Y}$ cannot, in general, be fixed in a $G$-equivariant way.  
The rigorous interpretation is that the innermost integral of \eqref{factorizeintegral} is valued in a certain \emph{equivariant sheaf} over $X^{h,Y}$, whose fiber over the horocycle $M$ is dual to the invariant measures on this horocycle. The prime that appears in the notation $ C_c^\infty(X)\to C^\infty((A_Y, {^{w_1}\chi})\backslash X^{h,Y})'$ is supposed to remind us that the image of the map above is does not lie in scalar-valued functions, but in sections of this sheaf. This sheaf is induced from a certain character of $\ker(A\to A_Y)$, explicated, e.g., in \cite[\S~5.2]{SaSpc}, but we will not need to recall this calculation here. 

For clarity, we will explicate noncanonical isomorphisms with principal series representations, as we did for the asymptotic cone in \eqref{isomprincipal}. First of all, we can fix a ``standard'' split Cartan subgroup in a ``standard'' Borel subgroup $B$, which will determine the choice of lifts of the Weyl group elements, up to elements of that Cartan, so that it makes sense to write $w_1Bw_1^{-1}$. 
\begin{itemize}
 \item For the horospherical space $X^h =X_\emptyset^h$, we have $C^\infty((A_X,\chi)\backslash X^h)\simeq I_{P(X)}^G(\chi)$. Note that this is a subrepresentation of $I_B^G(\chi\delta_{L(X)}^{-\frac{1}{2}})$, where $\delta_{L(X)}$ is the modular character of the Borel of $L(X)$. 
 \item For the $Y$-horospherical variety $X^{h,Y}$, what we denoted above by $C^\infty((A_Y, {^{w_1}\chi})\backslash X^{h,Y})'$ is isomorphic to $I_{B}^G({^{w_1}(\chi\delta_{L(X)}^{-\frac{1}{2}})})$.
 \item The outer integral of \eqref{factorizeintegral} represents a standard intertwining operator 
 \[ I_B^G({^{w_1}(\chi\delta_{L(X)}^{-\frac{1}{2}})}) \to I_B^G(\chi\delta_{L(X)}^{-\frac{1}{2}}),\] 
 that takes the image of $\mathfrak N_{\chi}^Y$ into the subspace $I_{P(X)}^G(\chi)$. 
\end{itemize}

We can now expand the left-hand triangle of Diagram \eqref{fiberscatteringdiagram} to 
\begin{equation}\label{fiberscatteringdiagram2} \xymatrix{ &&  C^\infty((A_Y, {^{w_1}\chi^{-1}})\backslash X^{h,Y})'  \ar@{-->}[rr]^{\mathfrak T^Y_{{\chi^{-1}}}}\ar@{-->}[dd]  && C^\infty((A_X,{\chi^{-1}})\backslash X_\emptyset^h)  \ar@{-->}[dd]
\\ C_c^\infty(X) \ar[urr]^{\mathfrak N^Y_{\chi^{-1}}}\ar[urrrr]_(0.7){\mathfrak N_{\chi^{-1}}}\ar[drr]_{\mathfrak N^Y_\chi} \ar[drrrr]^(0.7){\mathfrak N_{\chi}}&&&  \\ 
&& C^\infty((A_Y, {^{w_1}\chi})\backslash X^{h,Y})' \ar@{-->}[rr]_{\mathfrak T^Y_{\chi}}  && C^\infty((A_X,{\chi})\backslash X_\emptyset^h),}
\end{equation}
where $\mathfrak T^Y_{\chi}$ is the transform represented by the outer integral of the right hand side of \eqref{factorizeintegral}, and we used horizontal dotted arrows to signify that the map is only defined in the image of the maps $\mathfrak N^Y_{\chi^{-1}}$, resp.\ $\mathfrak N^Y_{\chi}$.\footnote{The reader is advised to remember the isomorphism of  $C^\infty((A_X,\chi^{\pm 1})\backslash X^h)$ with the possibly degenerate principal series representation $I_{P(X)}^G(\chi^{\pm 1})$, and to think of these dotted arrows as actual maps into the non-degenerate principal series $I_B^G(\chi^{\pm 1}\delta_{L(X)}^{-\frac{1}{2}})$.} It is also some kind of Radon transform/standard intertwining operator, but between principal series representations realized on two different varieties. Note that this Radon transform can be expressed in terms of the canonical $G$-orbit in $X^{h,Y}\times \tilde X^h$, described at the end of \S~\ref{sscanonicalRadon}, by a formula analogous to \eqref{Radon-Gorbit}.\footnote{Note that the recipe of \S~\ref{sssmeasureRadon} for the choice of measure for this transform does not apply (because the source and target varieties are different), but this ambiguity cancels out with the ambiguity that forced us to denote the source of the map by $(\,\,)'$: the two choices of measures combine to a choice of measures for the horocycle $M\tilde w_1 N \tilde w_1^{-1}$ of the middle expression of \eqref{factorizeintegral}. The latter choice, we recall, can be performed $G$-equivariantly along all such ``generic'' horocycles, and will be canceled by the corresponding choice for the definition of $\mathfrak M_\chi$ in \eqref{fiberscatteringdiagram}.}

If we fix the isomorphisms \eqref{noncan1}, \eqref{noncan2}, this diagram becomes 
\begin{equation}\label{fiberscatteringdiagram3} \xymatrix{ &&  I_{B}^G({^{w_1}(\chi^{-1}\delta_{L(X)}^{-\frac{1}{2}})})  \ar@{-->}[rr]^{\mathfrak T^Y_{{\chi^{-1}}}}\ar@{-->}[dd]  && I_{P(X)}^G(\chi^{-1}) \ar@{-->}[dd]
\\ C_c^\infty(X) \ar[urr]^{\mathfrak N^Y_{{\chi^{-1}}}}\ar[urrrr]_(0.7){\mathfrak N_{\chi^{-1}}}\ar[drr]_{\mathfrak N^Y_{{\chi}}} \ar[drrrr]^(0.7){\mathfrak N_{\chi}}&&&  
\\ 
&& I_{B}^G({^{w_1}(\chi\delta_{L(X)}^{-\frac{1}{2}})}) \ar@{-->}[rr]_{\mathfrak T^Y_{\chi}}  && I_{P(X)}^G(\chi),}
\end{equation}
with the operators $\mathfrak T^Y_{{\chi}^{\pm 1}}$ translating to the standard intertwining operators \eqref{Radon-standard} corresponding to a lift $\tilde w_1^{-1}$ of the Weyl element $w_1^{-1}$ (up to the choice of measure, which however, will be made simultaneously for $\mathfrak T^Y_{\chi}$ and $\mathfrak T^Y_{\chi^{-1}}$, and therefore will not affect the right dotted arrows of the commutative diagram).

From now on we omit the character from the notation of these spectral Radon transforms. If we could compute the left vertical arrow of \eqref{fiberscatteringdiagram3} as a scalar multiple $c_{w'}(\chi) \mathfrak R_{\tilde w'}$ of a Radon transform, then the right vertical arrow would be equal to $c_{w'}(\chi) \mathfrak R_{\tilde w_1^{-1}} \circ \mathfrak R_{\tilde w'} \circ \mathfrak R_{\tilde w_1^{-1}}^{-1}$, which we could then express in terms of the transform $\mathfrak R_{\tilde w_1^{-1}}\circ  \mathfrak R_{\tilde w'} \circ \mathfrak R_{\tilde w_1} \overset{\eqref{Radondecomp}}{=} \mathfrak R_{\tilde w} $ using\footnote{\label{footnoteRadon}The indexing of Radon transform by Weyl elements here is as in \eqref{Radon-standard}, but with the integral taken over a quotient of $N$, not $N^-$. It should be clear from the diagrams which of the two conventions we are using, and hopefully this abuse of notation will make it easier to follow the arguments, rather than conjugating Weyl elements by the longest Weyl element every time we use the ``standard'' Borel $B$.} the Radon inversion formula \eqref{inverseRadon}.

Our next step is to relate the left vertical arrow of \eqref{fiberscatteringdiagram2} to the calculation of the ``basic cases'' of \S~\ref{sssreviewbasic}. The relevant basic case is the variety $Y_2:= Y/U_P$, where $P$ denotes the parabolic denoted by $P_\alpha$, resp.\ $P_{\alpha\beta}$, in the two cases of Lemma \ref{lemmaclosedorbits}. The operator $\mathfrak N^Y_{\chi}$ of \eqref{fiberscatteringdiagram2} can be further factored in the form 
\begin{multline}\label{NYfactor} \mathfrak N^Y_{\chi} \Phi(M) = \int_{A_Y} \int_{M_2} \int_{U_{P,y}\backslash U_P}  \Phi(y u a) du \, dy \, \cdot \\{^{w_1}\left(\chi^{-1}\cdot \delta_{P(X)}^{-\frac{1}{2}}\right)}\cdot \delta_{(N\cap\tilde w_1^{-1} N \tilde w_1)\backslash N} (a) da,\end{multline}
where we have written $M_2$ for the image of the horocycle $M$ modulo $U_P$. The innermost integral represents a $P$-equivariant map to a certain sheaf over $Y_2$ -- whose fiber over some point $\bar y$ is dual to the space of $U_P$-invariant measures on the fiber of $Y\to Y_2$ over $\bar y$. Again, by abuse of notation, we will denote the innermost integral as a surjection 
\[C_c^\infty(X) \twoheadrightarrow C_c^\infty(Y_2)',\] 
where the prime is supposed to remind us that these are not scalar-valued functions, but valued in the sheaf described above. As we vary $P$ over its conjugacy class, these maps combine to produce sections over a ``partial horocycle space'' associated to the pair $(Y,P)$ -- the analog of $X^{h,Y}$, with $B$ replaced by $P$. These can be identified with the induction of $C_c^\infty(Y_2)'$, that is, we have a $G$-morphism 
\[C_c^\infty(X) \to \Ind_P^G( C_c^\infty(Y_2)').\]
The two outer integrals of \eqref{NYfactor}, now,  represent the analog $\mathfrak N_{^{w_1}\chi}^{Y_2}$ of the operator $\mathfrak N_{\chi}$ for $Y_2$, with $\chi$ replaced by $^{w_1}(\chi^{-1}\cdot \delta_{P(X)}^{-\frac{1}{2}})$ as in \eqref{fiberscatteringdiagram3} (but we omit the modular character from the notation). We have a factorization of the morphisms $\mathfrak N^Y_{\chi^{\pm 1}}$,
\begin{equation}\label{fiberscatteringdiagram4} \xymatrix{ && \Ind_P^G( C_c^\infty((A_Y, {^{w_1}\chi^{-1}})\backslash Y_2^h)') =  C^\infty((A_Y, {^{w_1}\chi^{-1}})\backslash X^{h,Y})'  \ar@{-->}[dd]  
\\ 
C_c^\infty(X) \ar[r]  \ar[urr]^{\mathfrak N^Y_{\chi^{-1}}}\ar[drr]_{\mathfrak N^Y_\chi} & \Ind_P^G( C_c^\infty(Y_2)') \ar[ur] \ar[ur]_{\mathfrak N_{^{w_1}\chi^{-1}}^{Y_2}} \ar[dr]^{\mathfrak N_{^{w_1}\chi}^{Y_2}} &&  \\ 
&& \Ind_P^G( C_c^\infty((A_Y, {^{w_1}\chi})\backslash Y_2^h)') = C^\infty((A_Y, {^{w_1}\chi})\backslash X^{h,Y})'  
.}
\end{equation}

The point is now that, from the discussion of the basic cases in \S~\ref{sssreviewbasic} we essentially know the vertical arrow of the diagram above. 
 However, both the results of \S~\ref{sssreviewbasic} and the statement of the current theorem are about the scattering operators, i.e., the right vertical arrow of \eqref{fiberscatteringdiagram}, and to relate those, we need to compose with the appropriate inverse intertwining operators $\mathfrak M_{\chi^{\pm 1}}^{-1}$ of \eqref{fiberscatteringdiagram}. Adding those to the commutative square of Diagram \eqref{fiberscatteringdiagram2}, we obtain 
\begin{equation}\label{fiberscatteringdiagram5} \xymatrix{
\Ind_P^G( C_c^\infty((A_Y, {^{w_1}\chi^{-1}})\backslash Y_{2,\emptyset})')  \ar@/_8pc/[dddd]_{\mathscr S_{w',{^{w_1}\chi}}^{Y_2}} \ar@{<-}[d]^{\mathfrak M_{L,{^{w_1}\chi}^{- 1}}^{-1}}  \ar@{-->}[rr]^{?}  && C^\infty((A_X,{\chi^{-1}})\backslash X_\emptyset) \ar@/^6pc/[dddd]^{\mathscr S_{w,\chi}} \ar@{<-}[d]^{\mathfrak M_{\chi^{- 1}}^{-1}}
\\
\Ind_P^G( C_c^\infty((A_Y, {^{w_1}\chi^{-1}})\backslash Y_2^h)')   \ar@{-->}[dd]  \ar@{-->}[rr]^{\mathfrak T^Y_{{\chi^{-1}}}}  && C^\infty((A_X,{\chi^{-1}})\backslash X_\emptyset^h)  \ar@{-->}[dd]
\\ 
&&  \\ 
 \Ind_P^G( C_c^\infty((A_Y, {^{w_1}\chi})\backslash Y_2^h)')  \ar@{->}[d]^{\mathfrak M_{L,{^{w_1}\chi}}^{-1}} \ar@{-->}[rr]^{\mathfrak T^Y_{\chi}}  && C^\infty((A_X,{\chi})\backslash X_\emptyset^h) \ar@{->}[d]^{\mathfrak M_{\chi}^{-1}}
 \\
  \Ind_P^G( C_c^\infty((A_Y, {^{w_1}\chi})\backslash Y_{2,\emptyset})')  \ar@{-->}[rr]^{?}  && C^\infty((A_X,{\chi})\backslash X_\emptyset)
,}
\end{equation}
where ${\mathfrak M_{L,\chi^{\pm  1}}^{-1}}$ are the corresponding inverse Radon transforms for the Levi $L$ of $P$. 

The left curved arrow is the (induction from $P$ to $G$) of the scattering operator for $Y_2$, which was computed in \S~\ref{sssreviewbasic}.\footnote{The twist by certain sheaves, denoted by $~'$ throughout, plays a role in this calculation, as we will indicate in Example \ref{ex:A2} below.} Thus, if we calculate the arrows denoted by question marks, that will allow us to compute the scattering operator $\mathscr S_{w,\chi}$, as well.

To calculate the arrows $?$, we return to the diagram \eqref{fiberscatteringdiagram3}, and use the ``canonical'' $G$-orbit on $(X^h)^2$ (where we have fixed the isomorphism \eqref{noncan1} for $X^h$, so that, as explained above, the canonical $G$-orbit is represented by an involution $\tilde w\in G$) to identify $X^h \simeq X_\emptyset$ (in such a way that the canonical $G$-orbit corresponds to the distinguished ``tautological'' $G$-orbit on $X_\emptyset^h\times X_\emptyset$). We have now fixed an isomorphism between the bottom square of \eqref{fiberscatteringdiagram5}, and the square 
\begin{equation}\label{fiberscatteringdiagram6} \xymatrix{
 I_B^G({^{w_1}(\chi\delta_{L(X)}^{-\frac{1}{2}})}) \ar@{-->}[d]_{\mathfrak R_{\tilde w'}^{-1}} \ar@{-->}[rr]^{\mathfrak R_{\tilde w_1^{-1}}}  && I_{P(X)}^G(\chi)  \ar@{->}[d]^{\mathfrak R_{\tilde w}^{-1}}
 \\
 I_B^G({^{w_1}(\chi^{-1}\delta_{L(X)}^{-\frac{1}{2}})}) \ar@{-->}[rr]^{?}  && I_{P(X)}^G(\chi^{-1}),
}
\end{equation}
(and similarly for the top square of \eqref{fiberscatteringdiagram5}). 
By the standard property \eqref{Radondecomp} of Radon transforms, the arrow ? is the operator $\mathfrak R_{\tilde w_1}^{-1}$. 

To conclude, the scattering operator $\mathscr S_{w,\chi}$ can be computed from the scattering operator $\mathscr S_{w',{^{w_1}\chi}}^{Y_2}$ of the basic cases via the commutativity of the diagram 
\begin{equation}\label{fiberscatteringdiagram7} \xymatrix{
 I_B^G({^{w_1}(\chi\delta_{L(X)}^{-\frac{1}{2}})}) \ar[d]_{\mathscr S_{w',{^{w_1}\chi}}^{Y_2}} \ar@{-->}[rr]^{\mathfrak R_{\tilde w_1}^{-1}}  && I_{P(X)}^G(\chi) \ar[d]^{\mathscr S_{w,\chi}},
 \\
 I_B^G({^{w_1}(\chi^{-1}\delta_{L(X)}^{-\frac{1}{2}})}) \ar@{-->}[rr]^{\mathfrak R_{\tilde w_1}^{-1}}  && I_{P(X)}^G(\chi^{-1}),
}
\end{equation}
and the final result should be expressed in terms of the ``canonical'' Radon transform $I_{P(X)}^G(\chi)\to I_{P(X)}^G(\chi^{-1})$ which, under the isomorphisms that we fixed, is the standard Radon transform $\mathfrak R_{\tilde w}$. \footnote{Note that we are using, here, labeling conventions for the Radon transforms that correspond to induction from the ``standard'' Borel subgroup $B$, see Footnote \ref{footnoteRadon}; when applying formulas from \S~\ref{sscanonicalRadon}, they will need to be adjusted accordingly.}

The theorem now follows from the formula for $\mathscr S_{w',{^{w_1}\chi}}^{Y_2}$ recalled in \S~\ref{sssreviewbasic}, and formula \eqref{inverseRadon} for the inverse of Radon transform. The details of how the various gamma factors simplify to produce the final answer are essentially the same as in the unramified case, therefore I point the reader to \cite[\S~6.5]{SaSph}. In the three examples that follow, I will demonstrate the calculation for the cases $A_2$ and $C_3$, that the reader can compare with \cite[\S~6.6 and 6.17]{SaSph}, as well as the case of $F_4$, which was not included in \emph{op.\ cit.}

\end{proof}

\
\begin{example}\label{ex:A2}
	Let $X=\GL_2\backslash\SL_3$. For the sake of calculations, let us pick a point on $X$ and identify its stabilizer, $H=\GL_2$, with the subgroup of matrices of the form $\begin{pmatrix}
	\det g^{-1} \\ & g
	\end{pmatrix},$
	$g\in \GL_2$. Letting $B$ denote the upper triangular Borel subgroup, and labeling the simple roots $\alpha$, $\beta$ from top to bottom (so that the spherical root is $\gamma=\alpha+\beta$), the intersection $H\cap P_\alpha$ is equal to the upper triangular Borel subgroup of $H$. The open $B$-orbit in $H\cap P_\alpha\backslash P_\alpha$ will play the role of $Y^\circ$ in the argument above. We identify the Levi quotient $L_\alpha$ of $P_\alpha$ with $\GL_2$ through the top $2\times 2$ block of a matrix, and then the variety $Y_2$ is equal to $\Gm^2\backslash \GL_2\simeq \Gm\backslash\PGL_2$. .
	
	The inner integral of \eqref{NYfactor}, which corresponds to the map $C_c^\infty (X)\to \Ind_P^G(C_c^\infty(Y_2)')$ of \eqref{fiberscatteringdiagram4}, is an integral over the $(3,3)$ entry of a matrix (a copy of the additive group $\Ga$). The Haar measure on this $\Ga$ is \emph{not} invariant under the stabilizer $\Gm^2\subset\GL_2$ of a point on $Y_2$, but varies by the character $\eta:\diag(a,b)\mapsto |a^2 b|$. Thus, what we denote by $C_c^\infty(Y_2)'$ are sections of a $\GL_2$-equivariant sheaf over $Y_2$ -- the compact induction of this character.  The restriction of this sheaf to (the action of) $\SL_2$ is induced from the character $\diag(a,a^{-1})\mapsto |a|$ of $\Gm\subset \SL_2$.
	
	We can now calculate the scattering operator $\mathscr S_{w',{^{w_1}\chi}}^{Y_2}$ based on the case of the spherical variety $\Gm\backslash\GL_2$, discussed in \S~\ref{ssstorus-ext}. This scattering operator is a morphism between the $L_\alpha$-representations
	\[\Ind_{B_\alpha}^{L_\alpha}({^{w_\beta}\chi} \delta^\frac{1}{2}) = I_{B_\alpha^-}^{L_\alpha} ({^{w_\beta}\chi^{-1}})  \otimes |\det|^{\frac{3}{2}} \to  I_{B_\alpha^-}^{L_\alpha} ({^{w_\beta}\chi})  \otimes |\det|^{\frac{3}{2}}\]
	(taking into account that $\delta_{L(X)}^{-\frac{1}{2}}$ is trivial here), where $\delta=\delta_B$ is the modular character of the Borel of $G$, and we have denoted by $B_\alpha$, $B_\alpha^-$ the ``standard'' Borel subgroup of $L_\alpha$ and its opposite.

	Now, consider $L_\alpha\simeq\GL_2$ as a homogeneous space for $ \tilde L_\alpha:=\Gm^2\times^{\Gm}\GL_2$, and its derived group $\SL_2$ as a homogeneous space for $\Gm\times^{\mu_2}\SL_2\simeq\GL_2$.  
	The ``standard cocharacters'' $\check\epsilon_1$, $\check\epsilon_2$ with respect to the upper triangular Borel of $\GL_2$ can be written as $\frac{\check\mu + \check\alpha}{2}$, $\frac{\check\mu-\check\alpha}{2}$, respectively, in terms of the presentation $\Gm\times^{\mu_2}\SL_2$, where $\check\mu$ is the tautological cocharacter of $\Gm$. The character $\eta\otimes {^{w_\beta}\chi}|\det|^{\frac{3}{2}}$ descends to a character of the Cartan of $\tilde L_\alpha$, and restricts to a character $\chi_1$ of the Cartan of $\GL_2$; we compute:
	\[ e^{\check\epsilon_1}(\chi_1) (a) = e^{\frac{\check\alpha}{2}}({^{w_\beta}\chi})(a) |a|^{\frac{1}{2}}, \]
	\[ e^{\check\epsilon_2}(\chi_1) (a) = e^{-\frac{\check\alpha}{2}}({^{w_\beta}\chi})(a) |a|^{\frac{1}{2}}, \]
	where we have used the fact that ${^{w_\alpha}\chi}$ is trivial on the center of $L_\alpha'$, hence it makes sense to evaluate it on ${\frac{\check\alpha}{2}}$.

	Plugging this $\chi_1$ into \eqref{scattering-T-ext}, and using the fact that  $w_\beta\check\alpha = \check\gamma$,  we obtain
	\[\mathscr S_{w',{^{w_1}\chi}}^{Y_2} = \gamma(\chi,\frac{\check\gamma}{2},1,\psi^{-1})\gamma(\chi,\frac{\check\gamma}{2},0,\psi) \gamma(\chi,-\check\gamma, 0,\psi) \cdot \mathfrak R_{w_\alpha},\]
	where $\mathfrak R_{w_\alpha}$ is defined with respect to the ``canonical'' orbit introduced in \S~\ref{ssstorus-ext}, therefore 
	\begin{multline*}\mathscr S_{w,\chi}= \mathfrak R_{w_\beta}^{-1} \circ \mathscr S_{w',{^{w_\beta}\chi}}^{Y_2} \circ \mathfrak R_{w_\beta} \xlongequal{\eqref{inverseRadon}} 
	\gamma(\chi,\frac{\check\gamma}{2},1,\psi^{-1})\gamma(\chi,\frac{\check\gamma}{2},0,\psi) \gamma(\chi,-\check\gamma, 0,\psi) \\ \cdot  \gamma(\chi, \check\beta, 0, \psi^{-1}) \gamma(\chi, -\check\beta, 0, \psi) \mathfrak R_w.
	\end{multline*}
	Noting that in the cocharacter group of $A_X$ we have $\check\beta=\frac{\check\gamma}{2}$, this is 
	\begin{multline*}\gamma(\chi,\frac{\check\gamma}{2},1,\psi^{-1})\gamma(\chi,\frac{\check\gamma}{2},0,\psi) \gamma(\chi,-\check\gamma, 0,\psi) \\ \cdot  \gamma(\chi, \frac{\check\gamma}{2}, 0, \psi^{-1}) \gamma(\chi, -\frac{\check\gamma}{2}, 0, \psi) \mathfrak R_w \\
	= \gamma(\chi, \frac{\check\gamma}{2}, 0, \psi^{-1})\gamma(\chi,\frac{\check\gamma}{2},0,\psi) \gamma(\chi,-\check\gamma, 0,\psi) \mathfrak R_w.
	\end{multline*}

\end{example}

\begin{example}\label{ex:C3}
	Let $X=\Sp_2\times \Sp_4\backslash \Sp_6$. The spherical root is $\gamma=\alpha_1+2\alpha_2+\alpha_3$, with $\alpha_3$ denoting the long simple root, while the Levi $L(X)$ is the one that contains the simple roots $\alpha_1, \alpha_3$ (isomorphic to $\SL_2\times\Gm\times \SL_2$). The decomposition $w=w_1^{-1} w' w_1$ of the nontrivial element of $W_X$, as above, can be taken to be $w'=w_{\alpha_1}$ with 
	\[ w_1 = w_{\alpha_2}w_{\alpha_3}w_{\alpha_2}.\]
	
	Thus, a ``Brion path'' for the action of the Weyl group on the Borel orbits of maximal rank, defined by Knop \cite{KnOrbits}, is given by the following; I point the reader to \cite[\S~6]{SaSph} for details.
	
	\entrymodifiers={++[o][F-]}
	
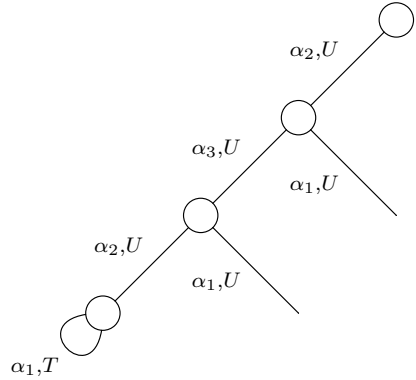
	\[\xymatrix{
		*{} &*{} &*{} & \ar@{-}[dl]_{\alpha_2,U}    & *{}& *{}& *{}
		\\
		*{} &*{} & \ar@{-}[dl]_{\alpha_3,U} \ar@{-}[dr]_{\alpha_1,U} &*{} &*{} 
		\\
		*{} &  \ar@{-}[dl]_{\alpha_2,U} \ar@{-}[dr]_{\alpha_1,U}&*{} &*{} &*{} 
		\\
		\ar@{-}@(d,l)[]^{\alpha_1,T}  & *{}& *{}& *{}
		\\
	}\]

	To calculate the right hand side of the equality $\mathscr S_{w,\chi}= \mathfrak R_{w_1}^{-1} \circ \mathscr S_{w',{^{w_1}\chi}}^{Y_2} \circ \mathfrak R_{w_1}$, in the notation of \eqref{fiberscatteringdiagram7}, we can invert the simple factors $\mathfrak R_{w_{\alpha_i}}$ of $\mathfrak R_{w_1}$ step-by-step, mimicking the inductive argument of \cite[\S~6]{SaSph}. 
	
	The two lower nodes correspond\footnote{In the sense that there is a closed orbit for the corresponding parabolic, whose quotient by the radical of the parabolic is isomorphic to this variety.} to the variety $\GL_2\backslash \PGL_3$, as in Example \ref{ex:A2} above.	These nodes will contribute a factor as in the previous example, except that we have to take into account another twist. If $P=L\ltimes U$ denotes the Siegel parabolic of $\Sp_6$, we will identify its Levi quotient $L$ with the general linear group of the Lagrangian fixed by $P$ under the left standard representation of $\Sp_6$. We consider the closed $P$-orbit on $X'$ corresponding to the bottom two nodes of the diagram above; the stabilizer of a point on that orbit is $P_1:= P\cap \Sp_2\times \Sp_4 = $ the product of Siegel parabolics of $\Sp_2\times \Sp_4$. The integral of a function in $C_c^\infty(X')$ over $U$-orbits on $P_1\backslash P$ lands in the induction, from $\GL_1\times \GL_2$ to $L=\GL_3$, of the character $|\mu|^2 \otimes |\det|$, where $\mu$ denotes the tautological character of $\GL_1$. In particular, restricting the action to the derived group $\SL_3$ of $L$, we are not inducing the trivial character of $\GL_2$, but the character $|\det|^{-1}$. This requires a modification of the calculation in the previous example, where $C_c^\infty(Y_2)'$ will not be induced from the character $\diag(a,b)\mapsto |a^2b|$ any more, but from the character $\diag(a,b)\mapsto |a^3 b|$. 
	This changes the factor contributed by the two lower nodes, which, following the argument of the previous example, can be computed\footnote{More details on this computation: The induced representation $\Ind_{B_\beta}^{L_\beta}({^{w_1}\chi} \delta^\frac{1}{2})$ can be computed to be equal to the normalized induction $I_{B_\alpha}^{L_\alpha}({^{w_1}\chi} |e^{\rho-\alpha_1-\alpha_2-\alpha_3}|)$, and when we plug this into \eqref{scattering-T-ext}, the scattering operator $\mathscr S_{w_{\alpha_1},{^{w_1}\chi}}^{Y_2}$ turns out to be 
		\[\gamma(\chi,\frac{\check\gamma}{2},\frac{3}{2},\psi^{-1})\gamma(\chi,\frac{\check\gamma}{2},-\frac{1}{2},\psi) \gamma(\chi,-\check\gamma, 0,\psi) \cdot \mathfrak R_{w_{\alpha_1}}.\]
	Similarly, the factor contributed by inverting the Radon transform $\mathfrak R_{w_{\alpha_2}}$ is, now, 	
	\begin{multline*} \gamma({^{w_{\alpha_3}w_{\alpha_2}}(\chi\delta_{L(X)}^{-\frac{1}{2}})},\check\alpha_2,0,\psi^{-1})\gamma({^{w_{\alpha_3}w_{\alpha_2}}(\chi\delta_{L(X)}^{-\frac{1}{2}})},-\check\alpha_2,0,\psi)\\
		= \gamma(\chi,\frac{\check\gamma}{2},\frac{1}{2},\psi^{-1})\gamma(\chi,-\frac{\check\gamma}{2},-\frac{1}{2},\psi),
	\end{multline*}
	and the product of the two simplifies to give the result.		
	} to be 
	\begin{equation}\label{gammaexample}\gamma(\chi, \frac{\check\gamma}{2}, \frac{1}{2}, \psi^{-1})\gamma(\chi,\frac{\check\gamma}{2},-\frac{1}{2},\psi) \gamma(\chi,-\check\gamma, 0,\psi).
	\end{equation}
	
	The inversion of the Radon transform $\mathfrak R_{w_{\alpha_3}}$ corresponding to the middle edge contributes, by \eqref{inverseRadon}, a factor of 
	\begin{multline*} \gamma({^{w_{\alpha_2}}(\chi\delta_{L(X)}^{-\frac{1}{2}})},\check\alpha_3,0,\psi^{-1})\gamma({^{w_{\alpha_2}}(\chi\delta_{L(X)}^{-\frac{1}{2}})},-\check\alpha_3,0,\psi) \\ = \gamma(\chi,\frac{\check\gamma}{2},-\frac{1}{2},\psi^{-1})\gamma(\chi,-\frac{\check\gamma}{2},\frac{1}{2},\psi),
\end{multline*}
	while the inversion of $\mathfrak R_{w_{\alpha_2}}$ corresponding to the top edge contributes
	\begin{multline*} \gamma({(\chi\delta_{L(X)}^{-\frac{1}{2}})},\check\alpha_2,0,\psi^{-1})\gamma({(\chi\delta_{L(X)}^{-\frac{1}{2}})},-\check\alpha_2,0,\psi) \\ = \gamma(\chi,\frac{\check\gamma}{2},-\frac{3}{2},\psi^{-1})\gamma(\chi,-\frac{\check\gamma}{2},\frac{3}{2},\psi).
	\end{multline*}

	Taken together, the factors above give the formula 
	\[\mathscr S_{w,\chi}=  \gamma(\chi, \frac{\check\gamma}{2}, -\frac{3}{2}, \psi^{-1})\gamma(\chi,\frac{\check\gamma}{2},-\frac{1}{2},\psi) \gamma(\chi,-\check\gamma, 0,\psi) \mathfrak R_w.\]
	for the scattering operator.
	
\end{example}

\begin{example}\label{ex:F4}
Let $X=\Spin_9\backslash F_4$. The spherical root is $\gamma=\alpha_1+2\alpha_2+3\alpha_3+\alpha_4$, with the roots labeled long first, short last, while the Levi $L(X)$ is the one that contains the simple roots $\alpha_1, \alpha_2,\alpha_3$ (isomorphic to $\Sp_6\times\Gm$). The decomposition $w=w_1^{-1} w' w_1$ of the nontrivial element of $W_X$, as above, can be taken to be $w'=w_{\alpha_4}$ with 
\[ w_1 = w_{\alpha_3}w_{\alpha_2}w_{\alpha_3}w_{\alpha_1}w_{\alpha_2}w_{\alpha_3}w_{\alpha_4}.\]

Thus, a ``Brion path'' for the action of the Weyl group on the Borel orbits of maximal rank, defined by Knop \cite{KnOrbits}, is given by the following.

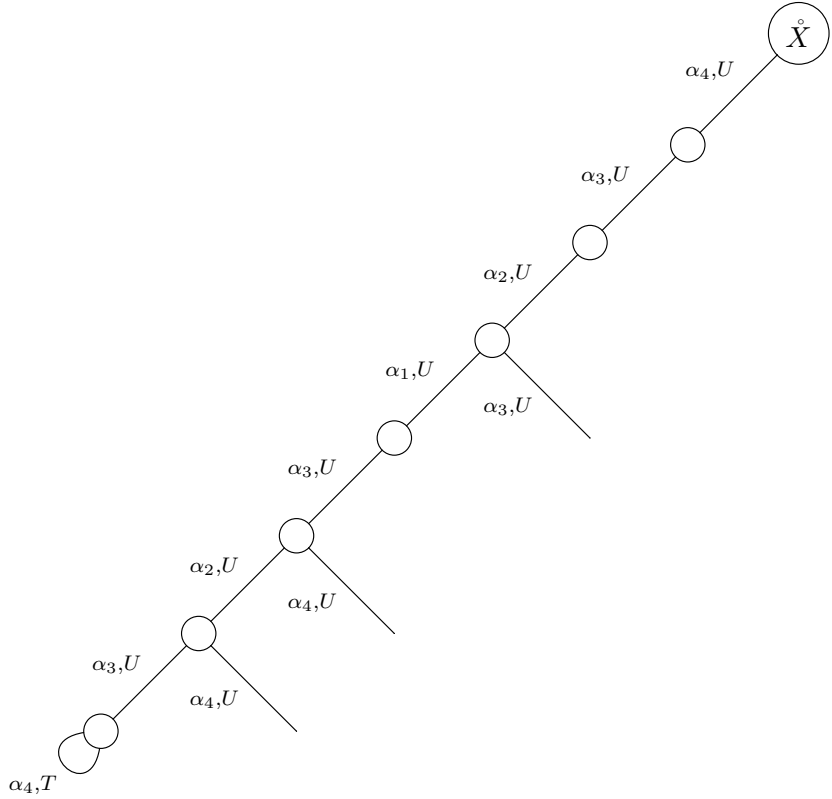
\entrymodifiers={++[o][F-]}

\[\xymatrix{
*{} &*{} &*{} & *{} & *{} &*{} &*{} & \mathring X \ar@{-}[dl]_{\alpha_4,U} & *{}
\\ 
*{} &*{} &*{} &*{} & *{} &   *{} & \ar@{-}[dl]_{\alpha_3,U} & *{}& *{}& *{}& *{}
\\ 
*{} &*{} &*{} & *{} &   *{} & \ar@{-}[dl]_{\alpha_2,U} & *{}& *{}& *{}& *{}
\\
*{} &*{} &*{} & *{} &   \ar@{-}[dl]_{\alpha_1,U} \ar@{-}[dr]_{\alpha_3,U}& *{}& *{}& *{}& *{}
\\
*{} &*{} &*{} & \ar@{-}[dl]_{\alpha_3,U}    & *{}& *{}& *{}
\\
*{} &*{} & \ar@{-}[dl]_{\alpha_2,U} \ar@{-}[dr]_{\alpha_4,U} &*{} &*{} 
\\
*{} &  \ar@{-}[dl]_{\alpha_3,U} \ar@{-}[dr]_{\alpha_4,U}&*{} &*{} &*{} 
\\
\ar@{-}@(d,l)[]^{\alpha_4,T}  & *{}& *{}& *{}
\\
}\]

Note that the three lower nodes correspond to the variety $\Sp_4\times \Sp_2\backslash \Sp_6$, as in Example \ref{ex:C3} above. 

The reader can now use the calculation of the previous example (there is no extra twist here, since the subgroup $\Sp_4\times \Sp_2$ has no nontrivial characters), together with the Radon inversion formula \eqref{inverseRadon}, to confirm the formula 
	\[\mathscr S_{w,\chi}=  \gamma(\chi, \frac{\check\gamma}{2}, -\frac{9}{2}, \psi^{-1})\gamma(\chi,\frac{\check\gamma}{2},-\frac{3}{2},\psi) \gamma(\chi,-\check\gamma, 0,\psi) \mathfrak R_w.\]

\end{example}

\begin{remark}
	The analog of the inductive hypothesis used in the proof of \cite[Proposition 6.5.1]{SaSph} for ``type $T$ spherical roots,'' here, is, with $\check\alpha$, $w$, $\check\gamma'$ as in that proposition, that the contribution of the nodes below a certain node in the diagram is the factor 
		\[ \gamma((\chi \delta_{L(X)}^\frac{1}{2})', \check\alpha, 0, \psi^{-1}) \gamma((\chi \delta_{L(X)}^\frac{1}{2})',-w\check\alpha,2-\langle\check\rho,\gamma'\rangle,\psi) \gamma(\chi,-\check\gamma, 0,\psi),\]
	where by $(\chi \delta_{L(X)}^\frac{1}{2})'$ we denote the conjugate of $\chi \delta_{L(X)}^\frac{1}{2}$ by the sequence of Weyl elements above the given node.
		
	For example, for the second-lowest node in Example \ref{ex:C3}, we have $\check\alpha=\check\alpha_2$, $\gamma'=w_{\alpha_2}\alpha_1=\alpha_1+\alpha_2$, and $w=$ the reflection corresponding to the root $\alpha_1+\alpha_2$, while the character $\chi \delta_{L(X)}^\frac{1}{2}$ has to be translated by the Weyl element $w_{\alpha_3} w_{\alpha_2}$. We then compute 
	\[\langle \check\alpha, w_{\alpha_3} w_{\alpha_2}\rho_{L(X)}\rangle = \frac{1}{2},\] 
	\[\langle -w\check\alpha, w_{\alpha_3} w_{\alpha_2}\rho_{L(X)} \rangle + 2 - \langle \check\rho, \gamma'\rangle = -\frac{1}{2},\]
	resulting in the gamma factors of \eqref{gammaexample}.
		
	Note also that the Brion path that was used for the calculation of the $F_4$-case above does not satisfy the assumptions of Proposition 6.4.1 of \emph{op.cit.} (at the node with edges labeled $\alpha_1$ and $\alpha_2$), but nonetheless satisfies this inductive hypothesis. 
\end{remark}

\section{Degeneration of transfer operators} \label{sec:degeneration}

In the previous section, we saw how the $L$-value associated to the spherical varieties of Table \eqref{thetable} controls, through its gamma factors, the scattering operators associated to the theory of asymptotics. 

On the other hand, in \cite{SaRankone} it was discovered that the $L$-value also has a different function: It controls certain ``transfer operators,'' which translate stable orbital integrals for the quotient $\mathfrak X = (X\times X)/G$ to orbital integrals for the associated Kuznetsov quotient $\mathfrak Y = (N,\psi)\backslash G^*/(N,\psi)$, where $G^*=\PGL_2$ or $\SL_2$, according as $\check G_X=\SL_2$ or $\PGL_2$, respectively. 

In this section, we will see how the two results are related via the \emph{degeneration of the transfer operators to the asymptotic cone}. 

\subsection{Transfer operators for rank-one spherical varieties} 

We keep assuming that the field $F$ is non-Archimedean. Here, we will need to work with measures, rather than functions, so we will use $\mathcal S(X)$ to denote the space of Schwartz measures (compactly supported, smooth) on the points of a variety $X$ over $F$. Note that the varieties of Table \eqref{thetable} all carry a $G$-invariant measure, so the translation from functions to measures is quite innocuous (up to the noncanonical choice of such a measure). But measures have natural pushforwards, and, using $\mathfrak X$ as a symbol for the quotient of $X\times X$ by the diagonal action of $G$, we will denote by $\mathcal S(\mathfrak X)$ the image of the pushforward map, from Schwartz measures on $X\times X$, to measures on the invariant-theoretic quotient $(X\times X)\sslash G:= \Spec F[X\times X]^G$, which in our cases is just an affine line. 

The meaning of $\mathcal S(\mathfrak Y)$ for the Kuznetsov quotient of $G^*$ is similar, but because of the twist by the Whittaker character $\psi$, we need to fix some conventions, as in \cite[\S~1.3]{SaRankone}. The result is a space of measures on the quotient $N\backslash G^*\sslash N$, which is again isomorphic to the affine line. 

The main theorem of \cite{SaRankone} is the following.

\begin{theorem}[{\cite[Theorem 1.3.1]{SaRankone}}] \label{thmtransfer}
For the varieties of Table \eqref{thetable}, and for appropriate identifications of the invariant-theoretic quotients $(X\times X)\sslash G$ and $N\backslash G^*\sslash N$ with the affine line (with coordinates that we denote by $\xi$ or $\zeta$, depending on whether $G^*=\PGL_2$ or $\SL_2$, respectively), the operator $\mathcal T$ described below gives rise to an injection
\begin{equation}\label{mapT} \mathcal T: \mathcal S(\mathfrak Y) \to \mathcal S(\mathfrak X).
\end{equation}

\begin{itemize}
 \item When $\check G_X=\SL_2$ with $L_X=L(\Std, s_1) L(\Std,s_2)$, $s_1\ge s_2$,
 \[    \mathcal T f(\xi) =  |\xi|^{s_1-\frac{1}{2}}  \left( |\bullet|^{\frac{1}{2}-s_1} \psi(\bullet) d\bullet\right) \star  \left( |\bullet|^{\frac{1}{2}-s_2} \psi(\bullet) d\bullet\right) \star f(\xi).\]
 \item When $\check G_X=\PGL_2$ with $L_X=L(\Ad, s_0)$,
\[    \mathcal T f(\zeta) =  |\zeta|^{s_0-1}  \left( |\bullet|^{1-s_0} \psi(\bullet) d\bullet\right) \star  f(\zeta).\]
\end{itemize}
 Here, $\star$ denotes multiplicative convolution on $F^\times$, $D\star f(x) = \int_{a\in F^\times} D(a) f(a^{-1} x)$.
\end{theorem}

The operator is bijective for a certain explicit enlargement of the space $\mathcal S(\mathfrak Y)$ of measures, that we will not recall here.
Conjecturally, the transfer operator also translates \emph{relative characters} of one quotient to relative characters of the other. This has been proven in several cases \cite{GW,SaTransfer1,SaTransfer2}. We will recall the notion of relative characters below, noting that in the case of $\SL_2 = \SO_3\backslash\SO_4$ they coincide with the usual \emph{stable} characters, while in the case of the Kuznetsov quotient $\mathfrak Y$ they are often called ``Bessel distributions.''

\subsection{Asymptotics of test measures}

To relate transfer operators to the scattering maps computed in the previous section, we 
recall \cite[Theorem 1.8]{DHS} that the asymptotics map \eqref{asymptotics} restricts to a morphism 
\begin{equation} e_\emptyset^*: \mathcal S(X)\to \mathcal S^+(X_\emptyset),
\end{equation}
 where $\mathcal S^+(X_\emptyset)$ denotes a certain enlargement of $\mathcal S(X_\emptyset)$, namely, a space of smooth measures on $X_\emptyset$, whose support has compact closure in an affine embedding (in this case, the ``affine closure'' $\Spec F[X_\emptyset]$). The ``restriction'' of the asymptotics map \eqref{asymptotics} from smooth functions to Schwartz measures makes sense, because 
by \cite[\S~4.2]{SV} an invariant measure on $X$ canonically induces an invariant measure on $X_\emptyset$. 

The spectral scattering maps $\mathscr S_{w,\chi}$ studied in the previous section are actually the Mellin transforms of a scattering operator $\mathfrak S_w$, an involution on $\mathcal S^+(X_\emptyset)$ which we think of as an action of the Weyl group $W_X\simeq \mathbb Z/2$.  
This involution is $(A_X,w_\gamma)$-equivariant, that is, it intertwines the action of $a\in A_X$ with the action of ${^{w_\gamma}a}=a^{-1}$, when this action is normalized to be unitary. Since we are working with measures, here, the unitary action analogous to \eqref{norm-functions} is 
\begin{equation}\label{norm-measures} a\cdot f(Sg) = \delta_{P(X)}^{-\frac{1}{2}}(a) f(Sag).
\end{equation}

Because of the equivariance property of $\mathfrak S_w$, it descends to an operator from the $(A_X,\chi^{-1})$-coinvariants to the $(A_X,\chi)$-coinvariants of $\mathcal S^+(X_\emptyset)$; for $\chi$ in general position, these can be identified with the corresponding coinvariants of the standard Schwartz space:

\begin{lemma}\label{lemmachicoinvariants}
For an open dense set of $\chi \in \widehat{A_X}_\CC$ (the complex Lie group of characters of $A_X$), the inclusion $\mathcal S(X_\emptyset)\hookrightarrow \mathcal S^+(X_\emptyset)$ induces an isomorphism on $(A_X,\chi)$-coinvariants, hence identifying those with $\mathcal S((A_X,\chi)\backslash X_\emptyset)$, the space of smooth measures on $A_X\backslash X_\emptyset$, valued in the sheaf whose sections are $(A_X,\chi)$-equivariant functions (for the normalized action \eqref{norm-functions} on functions) on $X_\emptyset$. Moreover, the ``twisted pushforward maps''
\begin{equation}\label{Mellin} \mathcal S^+(X_\emptyset)\to \mathcal S((A_X,\chi)\backslash X_\emptyset),\end{equation}
are meromorphic in $\chi$.
\end{lemma}

\begin{proof}
 This follows directly from the description of $\mathcal S^+(X_\emptyset)$ as a ``fractional ideal'' in the space of rational sections of $\chi\mapsto \mathcal S((A_X,\chi)\backslash X_\emptyset)$, in \cite[(1.18)]{DHS}.
\end{proof}

We will denote the map \eqref{Mellin} by $f\mapsto \check f(\chi)$, and think of it as a Mellin transform.

The spectral scattering morphisms of the previous section are the meromorphic family of operators descending from $\mathfrak S_w$ through this map:
\[\mathscr S_{w_\gamma,\chi} : \mathcal S((A_X,\chi^{-1})\backslash X_\emptyset) \to \mathcal S((A_X,\chi)\backslash X_\emptyset).\]
Note that, up to a choice of measure, $\mathcal S((A_X,\chi)\backslash X_\emptyset)$ is what was denoted before by $C^\infty((A_X,\chi)\backslash X_\emptyset)$.

Similar maps exist for the Whittaker model $\mathcal S(N,\psi\backslash G^*)$, with $e_\emptyset^*$ there mapping to a space $\mathcal S^+(N\backslash G^*)$ of measures on the space $N\backslash G^*$ (the ``degenerate Whittaker model,'' with the trivial character on $N$). In what follows, we will also be denoting by $Y$ the ``quotient'' $(N,\psi)\backslash G^*$ (i.e., the space $N\backslash G^*$ endowed with a sheaf defined by the nondegenerate character $\psi$), and by $Y_\emptyset$ its degeneration (the space $N\backslash G$ with the trivial sheaf). The image of the Schwartz space $\mathcal S(Y):= \mathcal S(N,\psi\backslash G^*)$ under the asymptotics map will be denoted by $\mathcal S^+(N\backslash G^*)$:
 \[ e_\emptyset^*: \mathcal S(N,\psi\backslash G^*) \to \mathcal S^+(N\backslash G^*).\]

Now, let us repeat the construction of test measures for the relative trace formula, for the asymptotic cones of $X$ and $Y$. Denote by $\mathcal S^+(X_\emptyset\times X_\emptyset/G) $ the pushforward of $\mathcal S^+(X_\emptyset)\otimes \mathcal S^+(X_\emptyset)$ to the invariant-theoretic quotient $(X_\emptyset\times X_\emptyset)\sslash G$. We adopt the same convention for $(N\backslash G^* \times N\backslash G^*)/G^*= N\backslash G^*/N$, letting $\mathcal S^+(N\backslash G^*/N)$ be the image of the pushforward map from $\mathcal S^+(N\backslash G^*)\otimes \mathcal S^+(N\backslash G^*)$ to $N\backslash G^*\sslash N = (N\backslash G^*) \times (N\backslash G^*)\sslash G^*$.

Composing the asymptotics maps with these pushforwards, we obtain maps
 \begin{equation}\label{diagramasymp}
 \xymatrix{\mathcal S(N,\psi\backslash G^*)\otimes \mathcal S(N,\psi^{-1}\backslash G^*)  
\ar[rr]^{(e_\emptyset^*\otimes e_\emptyset^*)_{G^*}} && \mathcal S^+(N\backslash G^*/N) \ar@{-->}[d]^{\mathcal T_\emptyset} 
\\
\mathcal S(X\times X) \ar[rr]^{(e_\emptyset^*\otimes e_\emptyset^*)_G} && \mathcal S^+(X_\emptyset\times X_\emptyset/G)},  
\end{equation}
with the map $\mathcal T_\emptyset$ to be introduced. The main result of this section will be a description of a canonical ``transfer operator'' $\mathcal T_\emptyset$, characterized by compatibility with relative characters; therefore, let us first discuss those. 

\subsection{Relative characters}

\subsubsection{}
Let $\pi$ be an admissible representation of $G$, with $\tilde\pi$ its contragredient. A \emph{relative character} on the quotient $\mathfrak X=(X\times X)/G$ for $\pi$ is a functional that factors 
\[ J_\pi: \mathcal S(X\times X)\to \pi\otimes\tilde\pi\to \CC,
\]
with the first map $G\times G$-equivariant and the second the defining pairing between $\pi$ and $\tilde\pi$. 

A source of relative characters is the Plancherel formula for $L^2(X)$: Once we choose a Plancherel measure $\mu_X$, as well as a $G$-invariant measure $dx$ in order to embed $\mathcal S(X)\hookrightarrow L^2(X)$, we have a decomposition of the bilinear pairing $\left< f_1, f_2\right> = \int_X \frac{f_1 f_2}{dx}$ as
\[ \left< f_1, f_2\right> = \int_{\hat G} J_\pi(f_1 \otimes f_2) \mu_X(\pi),\]
uniquely determined by this formula for $\mu_X$-almost every $\pi$ in the unitary dual $\hat G$. 

On the other hand, we can build relative characters by pullback via the asymptotics maps. Of interest to us here will be the relative characters that we can pull back from the space $\mathcal S^+(X_\emptyset\times X_\emptyset/G)$.
 Namely, there is a ``canonical'' open embedding $A_X\hookrightarrow X_\emptyset\times X_\emptyset\sslash G$, such that 
\begin{itemize}
 \item the identity maps to the image of the distinguished $G$-orbit $X_\emptyset^R\subset X_\emptyset\times X_\emptyset$, described in \S~\ref{sscanonicalRadon}. 
\item the map is equivariant with respect to the $A_X$-action descending from the action on the first, or equivalently the second, copy of $X_\emptyset$. 
\end{itemize}

Consider the meromorphic family of functionals $I_\chi$ obtained as pullbacks of the composition of maps
\[\xymatrix{ I_\chi: \mathcal S(X\times X)\ar[rr]^{(e_\emptyset^*\otimes e_\emptyset^*)_G} && \mathcal S^+(X_\emptyset\times X_\emptyset/G) \ar[rr]^{\,\,\,\,\,\,\,\,\,\,\,\, \int \chi^{-1}\delta_{P(X)}^{-\frac{1}{2}}} && \CC,}\]
 where the last arrow denotes the integral against the pullback of the character $\chi^{-1}\delta_{P(X)}^{-\frac{1}{2}}$ from $A_X$ to (a dense open subset of) $X_\emptyset\times X_\emptyset$. A priori, this arrow converges only on the subspace of such measures that are supported on $A_X\subset (X_\emptyset\times X_\emptyset)\sslash G$, but it is not hard to make sense of it for almost every $\chi$:
 
 \begin{lemma}\label{lemmarelchar}
The functionals $I_\chi$ converge when $|\chi(\varpi^{\check\gamma})| \ll 1$, and extend meromorphically to all $\chi \in \widehat{A_X}_\CC$ (the complex Lie group of characters of $A_X$). For an open dense set of $\chi$'s, the functional $I_\chi$ is a relative character for the normalized principal series representation $\pi_\chi=I_{P(X)}(\chi)$, that is, it factors through a morphism
 \[ \mathcal S(X\times X) \to \pi_\chi\otimes \widetilde{\pi_\chi} \xrightarrow{\left<\,\, , \,\, \right>} \CC.\]
\end{lemma}

\begin{proof}
 The elements of $\mathcal S^+(X_\emptyset\times X_\emptyset)$ are of moderate growth, and their support has compact closure in the affine closure of $X_\emptyset\times X_\emptyset$. Therefore,  the elements of $\mathcal S^+(X_\emptyset\times X_\emptyset/G)$ are of moderate growth, and their support has compact closure in $\mathbb A^1(F)$, when we identify $A_X\simeq \Gm$ via the character $\gamma$ or $\frac{\gamma}{2}$. Convergence for $|\chi(\varpi^{\check\gamma})| \ll 1$ follows. 
 
 Before we sketch the proof of meromorphic continuation, let us explain why they are relative characters. The composition 
 \[\xymatrix{I_\chi^\emptyset: \mathcal S^+(X_\emptyset)\otimes \mathcal S^+(X_\emptyset)\ar[rr] && \mathcal S^+(X_\emptyset\times X_\emptyset/G) \ar[rr]^{\,\,\,\,\,\,\,\,\,\,\,\, \int \chi^{-1}\delta_{P(X)}^{-\frac{1}{2}}} && \CC,}\]
is, by construction, $(A_X,\chi^{-1})$-equivariant with respect to the normalized action of $A_X$ on \emph{either} $\mathcal S^+(X_\emptyset)$-factor, and therefore the map factors through the $(A_X,\chi)$-coinvariants of each factor. By Lemma \ref{lemmachicoinvariants}, for $\chi$ in general position, these coinvariants are equal to $\mathcal S((A_X,\chi)\backslash X_\emptyset)$, which is isomorphic to $I_{P(X)^-}(\chi)\simeq \pi_{\chi^{-1}}$. Hence, $I_\chi$ factors through $\pi_{\chi^{-1}}\otimes\pi_{\chi^{-1}}$.  Since $\chi$ is $w$-conjugate to its inverse, for generic $\chi$ we have that $\pi_\chi \simeq \pi_{\chi^{-1}} = \widetilde{\pi_\chi}$.
 
Finally, for the claim of meromorphicity, it is not hard to express the map $I_\chi^\emptyset$, on the tensor product $\mathcal S((A_X,\chi)\backslash X_\emptyset)\otimes \mathcal S((A_X,\chi)\backslash X_\emptyset)$ of these coinvariant spaces, in terms of standard intertwining operators and the standard pairing between $\pi_{\chi}$ and $\pi_{\chi^{-1}}$. More precisely, if we use the canonical Radon transform introduced in \S~\ref{sscanonicalRadon} (translated to measures, by multiplying both sides by a $G$-invariant measure):
\[ \mathfrak R_\chi: \mathcal S((A_X,\chi^{-1})\backslash X_\emptyset) \to \mathcal S((A_X,\chi)\backslash X_\emptyset),\]
then the map $I_\chi^\emptyset$ factors through 
 \begin{multline}\label{Iinftyfactor} \mathcal S((A_X,\chi)\backslash X_\emptyset) \otimes \mathcal S((A_X,\chi)\backslash X_\emptyset)
  \xrightarrow{I\times \mathfrak R_{\chi^{-1}}} \\  \mathcal S((A_X,\chi)\backslash X_\emptyset) \otimes \mathcal S((A_X,\chi^{-1})\backslash X_\emptyset)  \xrightarrow{\left<\bullet\right>} \CC,
\end{multline}
where the arrow labeled $\left<\bullet\right>$ is the integral over the diagonal of $A_X\backslash X_\emptyset$, against an invariant measure that is prescribed by the measure used for the Radon transform. (I leave the details of this measure to the reader.)
\end{proof}

The factorization \eqref{Iinftyfactor} will be very useful in what follows, so let us record it as 
\begin{equation}
 I_\chi^\emptyset(f_1\otimes f_2) = \left< \check f_1(\chi), \mathfrak R_{\chi^{-1}}\check f_2(\chi)\right>.
\end{equation}
Of course, we could have applied Radon transform to the first, instead of the second factor.

 Similarly, consider the Kuznetsov quotient for $G^*$.  Identifying the Cartan $A^*$ of $G^*$ with the torus of diagonal elements through the upper triangular Borel, the embedding 
 \[A^*\to \begin{pmatrix} & -1 \\ 1\end{pmatrix} A^*\subset G^*\] 
 descends to an embedding of $A^*$ in $N\backslash G^*\sslash N$. 
 We similarly have a relative character $J_\chi$ on $\mathcal S(N,\psi\backslash G^*)\otimes \mathcal S(N,\psi^{-1}\backslash G^*)$, obtained as the composition  
 $$\xymatrix{ \mathcal S(N,\psi\backslash G^*)\otimes \mathcal S(N,\psi^{-1}\backslash G^*) \ar[rr]^{\,\,\,\,\,\,\,\,\,\, e_\emptyset^*\otimes e_\emptyset^*} && \mathcal S^+(N\backslash G^*/N) \ar[rr]^{\,\,\,\,\,\,\,\,\int \chi^{-1}\delta_{B^*}^{-\frac{1}{2}}} && \CC}.$$
 Notice that here we are integrating here against the character $\chi^{-1}\delta_{B^*}^{-\frac{1}{2}}$ of $A^*\subset N\backslash G^*\sslash N$, where $B^*\subset G^*$ is the Borel subgroup of $G^*$; this makes $J_\chi$ a relative character for the normalized principal series $I_{B^*}^{G^*}(\chi)$. The analogous statements of Lemma \ref{lemmarelchar} all hold in this setting.
 
 We identify $A^*$ with the Cartan $A_X$ of $X$ via the identification of the dual groups with $\SL_2$ or $\PGL_2$, i.e., so that the positive root of $G^*$ corresponds to the spherical root $\gamma$. 
 
 \subsubsection{}
 
 We will now discuss the relation of the relative characters $I_\chi$, $J_\chi$ to the Plancherel formulas for $L^2(X)$, resp.\ $L^2((N,\psi)\backslash G^*)$. 
 Recall \cite[Theorem 7.3.1]{SV} that the spaces $L^2(X)$, $L^2(Y)$ have discrete and continuous spectra, with the continuous spectra naturally parametrized by unitary characters of $A_X$, modulo inversion (i.e., modulo the action of $W_X$). More precisely, using the index $\emptyset$ for the orthogonal complement of the subspace spanned by relative discrete series (the images of irreducible subrepresentations $\pi\hookrightarrow L^2(X)$ or $L^2(Y)$), there are Plancherel decompositions
\begin{equation}\label{Plancherel-cont} L^2(X)_\emptyset\mbox{ or }L^2(Y)_\emptyset = \int_{\widehat{A_X}/W_X} \mathcal H_\chi d\chi,
\end{equation}
where $\widehat{A_X}$ denotes the unitary dual of $A_X$, and the unitary representation $\mathcal H_\chi$ can be identified with ``the'' $(A_X,\chi)$-equivariant (for the normalized action) Hilbert space completion of $\mathcal S(X_\emptyset)$ (resp.\ of $\mathcal S(Y_\emptyset)$).

To be precise, $L^2(X)$ is not quite a completion of $\mathcal S(X)$ (it is, rather, a completion of a space of half-densities), but, fixing a Haar measure on $X$ (and hence a Haar half-density), we can consider it to be so. We can similarly fix a Haar measure on $N\backslash G^*$, and take $d\chi$, in the decompositions above, to be such that it pulls back to a fixed Haar measure on $\widehat{A_X}$ under the finite map $\widehat{A_X}\to \widehat{A_X}/W_X$. The Plancherel decompositions \eqref{Plancherel-cont}, then, give rise to \emph{relative characters $J_\chi^Z$ (where $Z=X$ or $Y$)}, which are the pullbacks of the Hermitian forms of $\mathcal H_\chi$ to $\mathcal S(Z)$ -- but our convention will be to consider them as bilinear forms, i.e., as functionals on $\mathcal S(X)\otimes \mathcal S(X)$, resp.\ on $\mathcal S(N,\psi\backslash G^*)\otimes \mathcal S(N,\psi^{-1}\backslash G^*)$.

 \begin{proposition}
 Let $\mu(\chi)$ be the function introduced in Theorem \ref{thmscattering}. 
 
  The product $I_\chi \mu_X(\chi)$ is a family of relative characters on $X$ that is invariant under the $W_X$-action $\chi\mapsto \chi^{-1}$. 
  
For suitable choices of invariant measures, the most continuous part of the Plancherel decomposition of $L^2(X)$ (corresponding to the canonical subspace $L^2(X)_\emptyset\subset L^2(X)$) reads:
 \begin{equation}\label{PlancherelX} \left<f_1, \overline{f_2}\right>_\emptyset  =  \int_{\widehat{A_X}/W_X} I_\chi(f_1\otimes f_2) \mu_X(\chi) d\chi,\end{equation}
 where $d\chi$ is such that it pulls back to a Haar measure on $\widehat{A_X}$.
 
Similarly, for the Whittaker space $Y=(N,\psi)\backslash G^*$, the family of relative characters $\gamma(\chi,-\check\gamma, 0,\psi) J_\chi$ is invariant under $\chi\mapsto \chi^{-1}$, and the Plancherel formula for the most continuous part of $L^2(N,\psi\backslash G^*)$ reads:
 \begin{equation}\label{PlancherelW} \left<f_1, \overline{f_2}\right>_\emptyset  =  \int_{\widehat{A_X}/W_X} J_\chi(f_1\otimes f_2) \gamma(\chi,-\check\gamma, 0,\psi) d\chi.
 \end{equation}
 \end{proposition}

\begin{proof}
 The proof, which is based on Theorem \ref{thmscattering} and the Plancherel formula of \cite[Theorem 7.3.1]{SV}, is identical to that of \cite[Theorem 3.6.3]{SaTransfer1}. The Whittaker case is already included in \cite[Theorem 3.6.3]{SaTransfer1} (in the case of $G^*=\SL_2$, but the case of $G^*=\PGL_2$ is identical, since the scattering operators only depend on the pullback to $\SL_2$).
\end{proof}

The relative characters $I_\chi(f_1\otimes f_2) \mu_X(\chi)$ and $J_\chi(f_1\otimes f_2) \gamma(\chi,-\check\gamma, 0,\psi)$ of the proposition, which appear in the Plancherel formulas of two different spaces with the \emph{same} Plancherel measure, or any multiple of this pair by the same scalar, will be called ``matching'' relative characters.

\subsection{Asymptotic transfer operators}

Putting everything together, we can now prove the main result of this section, which is a generalization of the results of \cite[\S 4.3, 5]{SaTransfer1} to all varieties of Table \eqref{thetable}. 

Let us first fix coordinates: We identify $A_X=A^*$ with $\Gm$ via the character $\gamma$, when $\check G_X=\SL_2$, resp.\ $\frac{\gamma}{2}$, when $\check G_X=\PGL_2$. We previously embedded $A_X \subset X_\emptyset \times X_\emptyset \sslash G$ sending the identity to the image of the distinguished $G$-orbit, but here we will work with the \emph{opposite} embedding $\Gm\simeq A_X \hookrightarrow X_\emptyset\times X_\emptyset\sslash G$, sending $-1$ to the image of the distinguished $G$-orbit.\footnote{Such a choice was also made in \cite[\S~5]{SaTransfer1}, and should also have been made in \S~4.3 of \emph{op.cit.}, due to the correction noted in the proof of Lemma \ref{lemmainverseRadon}, and the subsequent appearance of the character $\psi$ in both factors of \eqref{scattering-typeG}. Of course, these choices are only made for compatibility of the final formulas with those of \cite{SaRankone}, which also depend on non-canonical choices of coordinates.}

Having identified an open subset of $(X_\emptyset\times X_\emptyset)\sslash G$ and of $N\backslash G^*\sslash N$ with , we will use the following coordinate for this space:
\begin{itemize}
 \item when $\check G_X=\SL_2$, an element of $A^*$ will be denoted $\begin{pmatrix} \xi \\ & 1\end{pmatrix}$;
 \item when $\check G_X=\PGL_2$, an element of $A^*$ will be denoted $\begin{pmatrix} \zeta \\ & \zeta^{-1}\end{pmatrix}$.
\end{itemize}
We will denote the smallest $A_X=A^*$-stable subspace of $\mathcal S^+(N\backslash G^*/N)$ (resp.\ of $\mathcal S^+(X_\emptyset\times X_\emptyset/G)$) containing the image of the map  $e_\emptyset^*\otimes e_\emptyset^*$ by $\mathcal S^+(N\backslash G^*/N)'$, resp.\ $\mathcal S^+(X_\emptyset\times X_\emptyset/G)'$.

\begin{theorem}\label{thmTempty}
 There is a unique $A_X$-equivariant operator 
 \[\mathcal T_\emptyset: \mathcal S^+(N\backslash G^*/N)' \to \mathcal S^+(X_\emptyset \times X_\emptyset/G)',\]
 such that, for all almost all $\chi\in \widehat{A_X}$, the pullbacks of 
 \begin{equation}\label{diagramdegen}
 \xymatrix{\mathcal S(N,\psi\backslash G^*)\otimes \mathcal S(N,\psi^{-1}\backslash G^*)
\ar[rr]^{e_\emptyset^*\otimes e_\emptyset^*} && \mathcal S^+(N\backslash G^*/N)' \ar[d]^{\mathcal T_\emptyset} \\
\mathcal S(X\times X) \ar[rr]^{e_\emptyset^*\otimes e_\emptyset^*} && \mathcal S^+(X_\emptyset\times X_\emptyset/G)' \ar[rrr]^{~\,\,\,\,\,\,\,\,\,\, \int \chi^{-1}\delta_{P(X)}^{-\frac{1}{2}}} &&& \CC,}
\end{equation}
are matching relative characters for $X$ and $(N,\psi)\backslash G^*$.

 Moreover, in the coordinates fixed above, the operator is given by the following formula: 
\begin{itemize}
 \item When $\check G_X=\SL_2$ with $L_X=L(\Std, s_1) L(\Std,s_2)$, $s_1\ge s_2$,
 $$    \mathcal T_\emptyset f(\xi) =  |\xi|^{s_1-\frac{1}{2}}  \left( |\bullet|^{\frac{1}{2}-s_1} \psi(\bullet) d\bullet\right) \star  \left( |\bullet|^{\frac{1}{2}-s_2} \psi(\bullet) d\bullet\right) \star f(\xi).$$
 \item When $\check G_X=\PGL_2$ with $L_X=L(\Ad, s_0)$,
$$    \mathcal T_\emptyset f(\zeta) =  |\zeta|^{s_0-1}  \left( |\bullet|^{1-s_0} \psi(\bullet) d\bullet\right) \star  f(\zeta).$$
\end{itemize}
 
\end{theorem}

The term ``$A_X$-equivariant'', here, refers to the normalized action of $A_X$ on $\mathcal S^+(X_\emptyset \times X_\emptyset/G)$ that descends from \eqref{norm-measures}, and, similarly, its analogously normalized action (but using the modular character $\delta_{B^*}$ instead of $\delta_{P(X)}$) on $\mathcal S^+(N\backslash G^*/N)$. 
The factor $|\xi|^{s_1-\frac{1}{2}}$, resp.\ $|\zeta|^{s_0-1}$, in the formula for $\mathcal T_\emptyset$ is due to the difference between the characters $\delta_{B^*}^{-\frac{1}{2}}$ and $\delta_{P(X)}^{-\frac{1}{2}}$ in the definition of the relative characters $I_\chi$ and $J_\chi$; in terms of the torus $A_X$, this factor can be written $|e^{\rho_{P(X)} - \frac{\gamma}{2}}| = |e^{\rho_{P(X)} - \rho_{B^*}}| = \delta_{P(X)}^{\frac{1}{2}}\delta_{B^*}^{-\frac{1}{2}}$.

It would have been more natural to work with half-densities instead of measures, in order to avoid these factors; the downside would be that the pushforward maps from half-densities on $X\times X$ to half-densities on $X\times X\sslash G$ are not completely canonical. However, there seems to be a distinguished way to define such pushforwards, that was used in \cite{SaICM} to give a quantization interpretation of the transfer operators. Although we will discuss that picture in the next subsection, we will not dwell any further, in this article, on the reformulation of these results using half-densities.

\begin{proof}
 Knowing the scattering operators by Theorem \ref{thmscattering}, the proof is then identical to that of \cite[Theorem 4.3.1]{SaTransfer1}.
\end{proof}

The remarkable feature of the formulas of Theorem \ref{thmTempty} is that the transfer operators for the boundary are given by \emph{exactly} the same formulas as for the original spaces, Theorem \ref{thmtransfer}. 

I expect that \eqref{diagramdegen} descends to a commutative diagram 
\begin{equation} \xymatrix{
\mathcal S(N,\psi\backslash G^*/N,\psi) \ar[rr]^{e_\emptyset^*\otimes e_\emptyset^*}\ar[d]^{\mathcal T} && \mathcal S^+(N\backslash G^*/N) \ar[d]^{\mathcal T_\emptyset} \\
\mathcal S(X\times X/G) \ar[rr]^{e_\emptyset^*\otimes e_\emptyset^*} && \mathcal S^+(X_\emptyset\times X_\emptyset/G) },
\end{equation}
where $\mathcal T$ is the transfer operator of Theorem \ref{thmtransfer}. This would imply that the relative characters under $\mathcal T$ satisfy:
\begin{equation}
 \mathcal T^*(\mu_X(\chi) I_\chi) = \gamma(\chi,-\check\gamma, 0,\psi) J_\chi.
\end{equation}
This was proven for the basic cases $A_1$ and $D_2$ in \cite{SaTransfer1}.

\section{Hankel transforms for the standard $L$-function of $\GL_n$} \label{sec:Hankel}

In this section, we change our setting, to discuss a close analog of transfer operators, the \emph{Hankel transforms} which realize the functional equation of $L$-functions at the level of trace formulas. Our goal is to give an interpretation of a theorem of Jacquet, and of its proof, from the point of view of quantization. 

Here, we will take $F=\mathbb R$, in order to use the language of geometric quantization (which involves line bundles and connections). The results can then be transferred formally to any local field, and indeed the theorem of Jacquet holds over any local field.\footnote{In fact, Jacquet wrote \cite{Jacquet} for non-Archimedean fields, but Jacquet's proof of the Theorem \ref{thmJacquet} that we are discussing here holds over Archimedean fields, as well.} For example, the flat sections of the connection $\nabla = \nabla^0 - i\hbar dx$ on the trivial (complex) line bundle over $\mathbb R$, where $\nabla^0$ is the usual connection, are the multiples of the exponential $e^{i\hbar x}$, and this can be replaced by an additive character when $F$ is any other local field. Similarly, Jacquet's formula (see Theorem \ref{thmJacquet} below) has a ``natural'' meaning (and is correct) over any local field. The reformulation of Jacquet's proof that we provide also makes sense over any such field, due to the theory of the Weil representation. However, since we do not develop a general theory for how to understand geometric quantization over general local fields, there is little benefit to complicating the presentation by including other local fields $F$; it should be straightforward for the interested reader to fill in the translations and verify that every step makes sense over any $F$.

\begin{remark}
 An important note on proofs and notation: The arguments that we present are mere reformulations of Jacquet's arguments. There is therefore no point in commenting again on convergence issues and other technical details that have been dealt with in \cite{Jacquet}. For the same reason, we will take the freedom to be a bit vague with some of our notation; in particular, we will often denote by $\mathcal D(X)$ an unspecified space of half-densities on a variety $X$. The rigorous mathematical interpretation of these spaces is that we \emph{start} with well-defined spaces, such as the space $\mathcal D(V)$ of \emph{Schwartz} half-densities on a vector space, and then the other spaces $\mathcal D(X)$ are images of this space under various integral constructions that make sense, as proven by Jacquet. Finally, we will sometimes write $\mathcal D(X)$ for a space of half-densities defined on an open dense subset $X^\bullet$ of $X$; it will be clear (and we will usually say) what this subset is.
\end{remark}

\subsection{The theorem of Jacquet}

In this section, we denote by $G$ the group $G=\GL_n$, and we will work with the Kuznetsov formula for $G$, with respect to the standard character $\psi: (x_{ij})_{ij} \mapsto \psi(\sum_{i=1}^{n-1} x_{i, i+1})$ of the upper triangular unipotent subgroup $N\subset G$, where, again, we use $\psi$ both for a fixed nontrivial unitary character of $F$, and for this character of $N$. Since $F=\RR$, $\psi$ has the form $\psi(x) = e^{i\hbar x}$, for some nonzero constant $\hbar$; we will be writing $dx$ for the self-dual Haar measure with respect to $\psi$ (or $|dx|$, when we want to distinguish the measure from the differential form). 

As in the case of $G^*$ in the previous section, we will identify the dense subset $N\backslash G_B/N \subset N\backslash G\sslash N$, where $G_B$ denotes the open Bruhat cell, with the Cartan $A$ of diagonal matrices, but this time (for compatibility with the literature), we will use the identification that descends from $A\to w A$, where $w$ is the antidiagonal permutation matrix (i.e., its entries are all $1$, rather than the antidiagonal of $(1,-1)$ that we used previously in this paper). Fortunately, these noncanonical choices will not play a role, once we reformulate our theorem in terms of geometric quantization -- our reformulation of \S~\ref{ssJacquetreform} will not depend on choices of representatives for the orbits.

It is essential here to work with half-densities, so we let $\mathcal D(X)$ be the space of Schwartz half-densities on the $F$-points of a smooth variety $X$. Those are products of Schwartz functions by half-densities of the form $|\omega|^\frac{1}{2}$, where $\omega$ is a nowhere-vanishing polynomial volume form on $X$. (In the examples at hand, one can always find such a form; in the general case, one would have to be more careful with the definitions.)

We are particularly interested in the case where $X=\Mat_n$, the space of $n\times n$-matrices, or $\Mat_n^*$, its linear dual. In calculations that follow, we will be identifying $\Mat_n^*$, as a space, with $\Mat_n$, via the (symmetric) trace pairing $(A,B)=  \tr(AB)$; note, however, that, defining the right $G\times G$-action on $\Mat_n$ as 
\begin{equation}\label{actionM} A\cdot (g_1, g_2):= g_1^{-1}A g_2,
\end{equation}
 the dual action on $\Mat_n^*$ is 
\begin{equation}\label{actionMstar}B\cdot(g_1,g_2) = g_2^{-1} B g_1.
\end{equation}

The two natural embeddings $G\hookrightarrow\Mat_n$ and $G\hookrightarrow\Mat_n^*$ differ by $g\mapsto g^{-1}$, once we apply the identification $\Mat_n^*=\Mat_n$. These embeddings allow us to restrict half-densities to $G$, and then we define twisted pushforwards to the double quotient space $N\backslash G\sslash N$, or rather to its open subset identified, as above, with the Cartan $A$: If $\varphi = \Phi (dg)^\frac{1}{2}$ is a half-density on $G$, where $dg$ is a Haar measure on $G$ and $\Phi$ is a function, its twisted push-forward to $A\hookrightarrow N\backslash G\sslash N$ is defined as the product 
\begin{equation}\label{pushfdensities} \delta^\frac{1}{2}(a) O_a(\Phi) (da)^\frac{1}{2},
\end{equation}
 where $\delta$ is the modular character of the Borel subgroup $B\supset N$,  $da$ is a Haar measure on $A$, and $O_a$ is the Kuznetsov orbital integral
\begin{equation}\label{orbital-Kuz} O_a(\Phi) = \int_{N\times N} \Phi(n_1 w a n_2) \psi^{-1}(n_1 n_2) dn_1 dn_2,
\end{equation}
where, again, $w$ is the antidiagonal permutation matrix. The choice of Haar measures used here will not affect any of the results of this section.

\begin{remark}\label{remarkhalfdensities}
 We pause to emphasize the importance of the definition \eqref{pushfdensities} for the pushforward of a half-density. Canonically, only measures admit pushforwards. The pushforward of a measure of the form $\Phi dg$ is easily computed to be equal to $\delta(a) O^0_a(\Phi) (da)$ (for compatible choices of Haar measures), where now $O_a^0$ is the orbital integral \eqref{orbital-Kuz}, but with the trivial character instead of $\psi$. (Or, better, use $O_a$ and think of it as a twisted pushforward.) Less canonically, the function $a\mapsto  O_a(\Phi)$ can be thought of as (twisted) ``pushforward of functions'' -- it is the natural definition that we obtain from fixing the same Haar measure $|d\underline{n}|$ on $N\times N$ for all its orbits (in the open Bruhat cell), and integrating; it is ambiguous up to a constant, that depends on the choice of Haar measure. Similarly, the definition \eqref{pushfdensities} corresponds to fixing a Haar measure on $N\times N$, taking its square root $|d\underline{n}|^\frac{1}{2}$ (a Haar half-density), and multiplying $\varphi$. The product $\varphi \cdot |d\underline{n}|^\frac{1}{2}$ is a ``measure along the fibers of the map to $A$, multiplied by a half-density in the transverse direction,'' and it makes sense to push it forward to a half-density on $A$. We should note, however, that, although natural, it is not completely clear why we should take the same Haar measure for every $N\times N$-orbit; the full justification of this choice will come with the proof of the main theorem, see Proposition \ref{proposition-main}.
\end{remark}

For reasons that have to do with the Godement--Jacquet method of representing $L$-functions, we will adopt the notations of \cite{SaHanoi}, and denote the images of $\mathcal D(\Mat_n)$, $\mathcal D(\Mat_n^*)$ under these pushforward maps by \\ $\mathcal D^-_{L(\Std,\frac{1}{2})}(\mathfrak Y)$, $\mathcal D^-_{L(\Std^\vee,\frac{1}{2})}(\mathfrak Y)$, respectively, where $\mathfrak Y$ stands as a symbol for the twisted Kuznetsov quotient $(N,\psi)\backslash G/(N,\psi)$. These are understood as spaces of half-densities on $A$, and they are genuinely different; for example, when $n=1$, the embedding $A\hookrightarrow \Mat_1$ attaches the point $0$ to $A\simeq F^\times$, while the embedding $A\hookrightarrow \Mat_1^*$ attaches the point $\infty$.

The theorem that follows is due to Jacquet \cite[Theorem 1]{Jacquet}; I present it in the reformulation of \cite[Theorem 9.1]{SaHanoi}. Thinking of half-densities for $\mathfrak Y$ as half-densities on the torus $A$ of diagonal element, we also identify the latter with the ``universal'' Cartan of $G$, via the quotient $B\to A$ of the upper-triangular Borel subgroup. Thus, it makes sense to write $e^\alpha: A\to \Gm$ for a root (where we use exponential notation for the character of $A$, reserving the additive notation $\alpha$ for its differential), and we also denote by $e^{\check\epsilon_i}:\Gm\to A$ the cocharacter into the $i$-th entry of the diagonal. The Hankel transform of the theorem that follows will be expressed in terms of  operators of the form ``multiplication by a function $\psi(e^\alpha)$ on $A$'' as well as \emph{multiplicative Fourier convolutions} along the cocharacters $e^{-\check\epsilon_i}$. Those will be denoted by $\mathcal F_{-\check\epsilon_i,\frac{1}{2}}$, and are given, explicitly, by 
\[\mathcal F_{-\check\epsilon_i,\frac{1}{2}} \varphi(\diag(a_1, \dots, a_n)) = \int_{F^\times} \varphi(\diag(a_1, \dots, a_i x, \dots, a_n)) |x|^\frac{1}{2} \psi^{-1}(x) d^\times x.\]
\emph{(Caution: Compared to the notation of \cite{SaHanoi}, we have changed $\psi$ to $\psi^{-1}$ here, and we will make some corresponding changes below.)}

\begin{theorem}\label{thmJacquet}
Let $G=\GL_n$. Consider the diagram
 \[\xymatrix{
 \mathcal D(\Mat_n) \ar[d]\ar[r]^{\mathcal F}& \mathcal D(\Mat_n^*)\ar[d]\\
 \mathcal D^-_{L(\Std,\frac{1}{2})}(\mathfrak Y) \ar@{-->}[r]^{\mathcal H_{\Std}} & \mathcal D^-_{L(\Std^\vee,\frac{1}{2})}(\mathfrak Y)
 }\]
where $\mathcal F$ denotes the equivariant Fourier transform: 
\[ \mathcal F(\varphi) (y) = \left(\int_{\Mat_n} \varphi(x) \psi(-\left<x, y\right>) dx^\frac{1}{2}\right) dy^\frac{1}{2}\]
(for dual Haar measures $dx$, $dy$ on $\Mat_n$ and $\Mat_n^*$ with respect to the character $\psi$).

There is a linear isomorphism $\mathcal H_{\Std}$ as above, 
making the diagram commute. Moreover, $\mathcal H_{\Std}$ is given by the following formula:
 
 \begin{equation}\label{HStd}
  \mathcal H_{\Std} = \mathcal F_{-\check\epsilon_1,\frac{1}{2}} \circ \psi(e^{-\alpha_1}) \circ \mathcal F_{-\check\epsilon_2,\frac{1}{2}} \circ \cdots \circ \psi(e^{-\alpha_{n-1}}) \circ \mathcal F_{-\check\epsilon_n, \frac{1}{2}}.
 \end{equation}
 \end{theorem}

Explicitly, and denoting diagonal matrices simply as $n$-tuples, \eqref{HStd} reads:
 \begin{multline}\label{HStd-expl}
  \mathcal H_{\Std} f (b_1, \dots, b_n)  = \\
 \int f(b_1 p_1, \dots, b_n p_n) \psi(-\sum_{i=1}^n p_i + \sum_{i=1}^{n-1} \frac{b_{i+1}}{p_i b_i}) |p_1 \cdots p_n|^\frac{1}{2} d^\times p_1 \cdots d^\times p_n.
 \end{multline}
 
\begin{remark}
 Opposite to the convention in \cite{SaHanoi} and \cite{Jacquet} (where $\psi$ is denoted by $\theta^{-1}$), we have changed the convention of Fourier transform, using the character $\psi^{-1}$ instead of $\psi$. This resulted in some sign changes, but will simplify some expressions in the rest of the section.
\end{remark}

\subsection{Cotangent reformulation of Jacquet's theorem} \label{ssJacquetreform}

\subsubsection{}
The reformulation of Jacquet's theorem that we will present is directly analogous to the interpretation of transfer operators given in \cite{SaICM}: Namely, we will interpret both of the  spaces $\mathcal D^-_{L(\Std,\frac{1}{2})}(N,\psi\backslash G/N,\psi)$ and $\mathcal D^-_{L(\Std^\vee,\frac{1}{2})}(N,\psi\backslash G/N,\psi)$ as geometric quantizations corresponding to two different Lagrangian foliations on the \emph{same} symplectic space, and then 
the operator of Formula \eqref{HStd} represents integrals over the leaves of such a foliation.

Let us start by introducing the notions of foliations and integrals over the leaves; our definitions are more restrictive than usual, but good enough for our purposes. 

\begin{definition}\label{deffoliation}
\begin{enumerate}
 \item We will call \emph{foliation} a smooth morphism of smooth varieties $\mathscr F: X\to Y$, and \emph{leaves} its fibers; we will also be denoting $Y$ by $X/\mathscr F$. If $X$ is symplectic, a foliation is \emph{Lagrangian} if its leaves are Lagrangian subvarieties. 

\item Given a foliation as above, let $(L,\nabla)$ be a smooth vector bundle on (the real points of) $X$, with a connection that is flat with trivial monodromies along the leaves of $\mathscr F$. 
The \emph{space of parallel half-densities along $\mathscr F$}, $\mathcal D_{\mathscr F}(X,L)$, is defined as follows: Its elements are half-densities on $Y$, valued in a line bundle $L_\mathscr F$, whose sections are those sections of $L$ that are $\nabla$-flat along the leaves of $Y$. 

\item Suppose, now, that $\mathscr F_1$, $\mathscr F_2$ are two Lagrangian foliations on a symplectic space $X$, and $(L,\nabla)$ is a smooth vector bundle on $X$, with a connection whose curvature is equal to $i\hbar \omega$ -- in particular, flat along the leaves of Lagrangian foliations. Assume, also, that the restriction of each foliation to the leaves of the other is also a foliation; in particular, the intersections of two leaves are smooth, and their tangent spaces coincide with the intersections of the tangent spaces of the leaves. If $Z$ is a leaf, say of $\mathscr F_1$, we will be writing $Z/\mathscr F_2 \subset Y_2$ for its space of leaves for the restriction of the $\mathscr F_2$-foliation. 

In this setting, the \emph{integral} of an element  $\varphi \in \mathcal D_{\mathscr F_1}(X,L)$ along the leaves of $\mathscr F_2$ is the element of $\mathcal D_{\mathscr F_2}(X,L)$ whose ``value'' at $x\in X$ is given by 
\begin{equation}\label{intfoliation} \int_{\mathscr F_{2,x}/\mathscr F_1} T_{x,z}\varphi(z) |\omega|^{\frac{1}{2}\dim (\mathscr F_{2,x}/\mathscr F_1)}, 
\end{equation}
provided that the integral converges. Here, $\mathscr F_{2,x}$ denotes the $\mathscr F_2$-leaf containing $x$, and $T_{x,z}$ denotes parallel transport from the point $z$ to the point $x$. The integral makes sense as a parallel half-density along $\mathscr F_2$, exactly as in the linear case \cite[\S~3.4]{WWLi} (see Remark \ref{remarkint} below for an explanation). In particular, our convention is that, for a volume form $\Omega$ on a smooth variety, $|\Omega|$ denotes the corresponding positive measure on the $\mathbb R$-points of the variety, that is induced from the self-dual measure on $F=\mathbb R$ with respect to the character $x\mapsto e^{i\hbar x}$.
\end{enumerate}

\end{definition}

\begin{remark}\label{remarkint}
 We explain why the outcome of integral \eqref{intfoliation} represents a half-density on $X/\mathscr F_2$ (valued in the line bundle $L_{\mathscr F_2}$): 
 
 Let $x\in X$, set $M=T_xX$, and let $\ell_1, \ell_2 \subset M$ be the tangent spaces of the leaves of $\mathscr F_1$, $\mathscr F_2$ passing through $x$. The symplectic form $\omega$ restricts to a symplectic form on the quotient $(\ell_1+\ell_2)/(\ell_1\cap \ell_2)$, giving rise to a volume form $\omega \wedge\cdots \wedge \omega$ ($\frac{1}{2}\dim (\ell_1+\ell_2)/(\ell_1\cap \ell_2)$ times). We have a canonical isomorphism of lines 
 \begin{equation}\label{isomlines}
 \det (M/\ell_1)\otimes \det ((\ell_1+\ell_2)/(\ell_1\cap\ell_2)) = \det (M/\ell_2) \otimes \det (\ell_2/(\ell_1\cap\ell_2))^{\otimes 2}.
 \end{equation}
 
 Taking the ``square root of the absolute value'' of the dual spaces, we deduce that the product of a half-density along $X/\mathscr F_1$ by $|\omega|^{\frac{1}{2}\dim (\mathscr F_{2,x}/\mathscr F_1)}:= |\omega\wedge\cdots\wedge\omega|^\frac{1}{2}$ can be canonically understood as a density (measure) on  $\mathscr F_{2,x}/\mathscr F_1$ times a half-density on $X/\mathscr F_2$. This is how the integral \eqref{intfoliation} gives rise to a half-density on $X/\mathscr F_2$, valued in the appropriate line bundle. For further discussion of these integrals along Lagrangians, see \S~\ref{sssLagrangians} below.
\end{remark}

\subsubsection{}

For the case at hand, let us adopt a basis-independent point of view, and write $V$ for an $n$-dimensional vector space, $V^*$ for its linear dual. The role of $\Mat_n$ will be played by $\End(V)\simeq V^*\otimes V$, and we identify its linear dual with $\End(V^*)\simeq V\otimes V^*$. The cotangent spaces of both are identified with the space $M:= \End(V) \oplus \End(V^*)$, and the foliations $M\to \End(V)$ and $M\to\End(V^*)$ will be called the ``vertical'' and ``horizontal'' foliation, respectively. Breaking the symmetry,\footnote{Only the product of the symplectic form by an chosen nonzero constant $\hbar$ matters for what follows; we are not really breaking the symmetry -- just choosing opposite constants for the two quotients.} we need to fix a sign for the symplectic form on this vector space, and we choose it to be 
\begin{equation}\label{sympM} \omega_M = \sum_j dT_j \wedge dT_j^*,
\end{equation}
 where $(T_j)_j, (T_j^*)_j$ are dual bases for $\End(V)$, $\End(V^*)$. 

The group $G = \GL(V)\simeq \GL(V^*)$ admits an open embedding into both $\End(V)$ and $\End(V^*)$, giving rise to two \emph{distinct} embeddings of $T^*G$ into $M$. 

Now we will define a smooth line bundle $L$ on (the real points of) $M$, equipped with a connection $\nabla$, whose curvature is $i\hbar\omega_M$, for some nonzero constant $\hbar$. Let $\psi$ be the additive character $x\mapsto e^{i\hbar x}$. In order not to privilege one Lagrangian over the other, we will define $L$ in a symmetric fashion: Smooth sections of $L$ are pairs $(\Phi,\Phi^*)$ consisting of $\Phi,\Phi^*\in C^\infty(M)$, with the property that 
\begin{equation}\label{phiphistar} \Phi^*(A,B) = \psi(-\langle A, B\rangle) \Phi(A,B).
\end{equation}
 If $\nabla^0$ denotes the standard flat connection on the trivial line bundle on $M$, then the connection $\nabla$ is defined as 
\begin{equation}\label{nabla}
\nabla= (\nabla^0 - i\hbar  \sum_j  T_j^* dT_j, \nabla^0 + i \hbar \sum_j T_j dT_j^*), 
\end{equation}
where again we have used a dual basis. In particular, a function $f$ on $\End(V)$ pulls back to a section $\tilde f = (\Phi,\Phi^*)$ of $L$ on $M$, given by 
\[\Phi(A,B)=f(A), \,\, \Phi^*(A,B)=f(A)\psi(-\left<A,B\right>),\] with the property that $\nabla_Z \tilde f =0$ for any vector field $Z$ which is parallel to the leaves of the ``vertical'' Lagrangian foliation $M\to \End(V)$. Similarly, a function $f$ on $\End(V^*)$ pulls back to 
\[(\Phi(A,B)= f(B) \psi(\left<A,B\right>), \,\, \Phi^*(A,B)=f(B)),\]
which is $\nabla$-flat along the leaves of the ``horizontal'' foliation $M\to \End(V^*)$. 

\subsubsection{}

We now introduce the \emph{Kuznetsov cotangent space} $T^*\mathfrak Y$. We first give a rather clumsy presentation of the Whittaker model, by fixing a maximal unipotent subgroup $N\subset G$, and a homomorphism $f : N\to \Ga$ which is a generic (i.e., nonzero on every simple root space). The Whittaker ``space'' $Y$ is then the homogeneous space $N\backslash G$, but we will think of every point $Ng$ of it as the set of triples $(Ng, N', f')$, where $N'$ is the maximal unipotent subgroup $g^{-1}Ng$, and $f'$ is the homomorphism $N'\to \Ga$ obtained from $f$ via conjugation by $g$. There are nicer presentations of the Whittaker model, which do not depend on the choice of a base point; for example, if $\dim V=2$, we can let $Y=$ the set of pairs $(v, \omega)$, where $v\in V^\times$ and $\omega$ is a non-zero alternating form on $V$, and we can endow the stabilizer $N\subset G$ of such a point with the homomorphism $f$ such that $v^* - n\cdot v^*= f(n) v$, for every $v^*\in V$ with $\omega(v, v^*)=1$. 
(Also, if the group were adjoint, the set of pairs $(N,f)$ consisting of a maximal unipotent subgroup and a generic morphism to $\Ga$ would constitute the Whittaker model.) However, for reasons of conciseness, we will fix a base point. 

We may, and will, identify $f$ with its differential, as well, which is an element of $\mathfrak n^*$. We will then symbolically write $f+n^\perp$ for the preimage of $f$ under the restriction map $\mathfrak g^*\to\mathfrak n^*$, and define the Whittaker cotangent space $T^*Y$ as the space $T^*Y=(f+n^\perp)\times^N G$. We define the opposite Whittaker cotangent space as $T^*Y^- = (-f+n^\perp)\times^N G$. 

The inclusion $(f+n^\perp)\hookrightarrow \g^*$ gives rise to a $G$-equivariant map $T^*Y\to \g^*$, and $T^*Y$ is naturally a Hamiltonian space with this moment map; same for $T^*Y^-$. 

Finally, the Kuznetsov cotangent space $T^*\mathfrak Y$ is defined as the Hamiltonian reduction of $T^*Y^- \times T^*Y$ with respect to the diagonal $G$-action:\footnote{Our notation stands for Hamiltonian reduction at $(f,f)$ with respect to a \emph{left} and a \emph{right} action of $N$; when we convert this to a right $(N\times N)$-action, it will correspond to Hamiltonian reduction at $(-f,f)$.}
\begin{equation}\label{Kuzcotang} T^*\mathfrak Y=T^*Y^-\times_{\mathfrak g^*, \pm}^G T^*Y = N\bbbslash_f T^*G \ssslash_f N,
\end{equation}
where $\pm$ denotes that the images under the moment map should be \emph{opposite}. Of course, the map $-f\mapsto f$ induces an equivariant map $T^*Y^-\to T^*Y$ which inverts the moment map, so we could also present $T^*\mathfrak Y$ as $T^*Y\times_{\mathfrak g^*}^G T^*Y$, but in order to avoid sign confusions in the calculations that follow, it is better to think of the presentation \eqref{Kuzcotang}.

For completeness, we use the last presentation to provide a standard description of this space, although we will not use it further in this paper. Recall the \emph{group scheme of regular centralizers} $J_G \to \c^*$, where $\c^*:= \g^*\sslash G$; it comes equipped with an isomorphism of its pullback to the regular locus $\g^{*,\reg}$ with the inertia group scheme for the adjoint action of $G$. A well-known result of Kostant says that $T^*Y$ is a $J_G$-torsor over $\g^{*,\reg}$. As a corollary, $T^*\mathfrak Y$ is \emph{canonically} isomorphic to $J_G$, with the isomorphism induced by the (action, projection) map $J_G\times_{\c^*} T^*Y \to T^*Y\times_{\mathfrak g^*} T^*Y$.

\subsubsection{} 

Next, we study the double Hamiltonian reduction 
\[\mathfrak M = N\bbbslash_f M \ssslash_f N = M_{-f,f}/(N\times N),\] 
where $M_{-f,f}=M \times_{\mathfrak n^*\times \mathfrak n^*} (\{-f\}\times\{f\})$ (with respect to the right action of $N\times N$). The two embeddings of $T^*G$ into $M$ induce two distinct embeddings of $J_G\simeq T^*\mathfrak Y$ into it. 

The coordinate-dependent embedding $A \hookrightarrow N\backslash G\sslash N$ that we defined when we identified $G$ with $\GL_n$ is more canonically an $A$-torsor $A_1 \hookrightarrow N\backslash \GL(V)\sslash N$, which further maps into the quotient $N\bbslash \End(V)\sslash N$. The coordinate ring of $N\bbslash \End(V)\sslash N$ is spanned by the semiinvariants (highest weight vectors) for $B\times B$ on $F[\End(V)]$, and similarly for the other spaces, which shows that $A_1$, the locus where these semiinvariants are nonvanishing, is open in each of those spaces, and its preimage is the open $B\times B$-orbit (Bruhat cell), which we will denote by $G_B\subset G$. We have the analogous torsor $A_2\subset N\bbslash \End(V^*)\sslash N$.  
In particular, the restriction of the maps
\[ N\bbbslash_f M \ssslash_f N \to N\bbslash \End(V)\sslash N,\]
\[ N\bbbslash_f M \ssslash_f N \to N\bbslash \End(V^*)\sslash N,\]
to $A_1$, resp.\ $A_2$, coincides with the maps 
\[ N\bbbslash_f T^*G_B \ssslash_f N \to A_1\subset N\bbslash \End(V)\sslash N,\]
\[ N\bbbslash_f T^*G_B \ssslash_f N \to A_2\subset N\bbslash \End(V^*)\sslash N,\]
respectively, but note again that these two refer to two \emph{different} embeddings of $T^*G_B$ into $M$. We will denote by $M^\circ$ the intersection of the two embeddings, which is an open subset of $M$. Explicitly, in coordinates identifying $M$ with pairs $(A,B)$ of $n\times n$ matrices, 
\[M^\circ = NwAN \times NAwN,\] 
with this presentation corresponding to how we think of the GIT quotients as $A$-torsors. From the freeness of the $N\times N$-action on $G_B$ it is immediate that $(T^*G_B)_{-f,f}$ is an $N\times N$-torsor over $T^*A$ (for each of the two embeddings above), and from this it is easy to deduce that $M^\circ_{-f,f}$ is an $N\times N$-torsor over $A_1\times A_2$.

To summarize, we have an open, $N\times N$-equivariant subset $M^\circ\subset M$, and the double Hamiltonian reduction $J^\circ:= N\bbbslash_f M^\circ \ssslash_f N$ admits two Lagrangian fibrations, $J^\circ\to A_1$ and $J^\circ\to A_2$, identifying $J^\circ$ with the product $A_1\times A_2$.  The symplectic space $J^\circ$ can be identified as an open subset of the Kuznetsov cotangent space $J_G=T^*\mathfrak Y$, but this identification depends on which of the two embeddings $T^*G\hookrightarrow M$ we choose.

\begin{lemma}\label{lemmalagrangian}
The fibers of the composition  $M^\circ_{-f,f} \to J^\circ\to A$ ($A=A_1$ or $A_2$) are Lagrangian subspaces. In particular, the connection $\nabla$ on the line bundle $L$ restricts to a flat connection on those fibers. The monodromy of this flat connection is trivial.
\end{lemma}

\begin{proof}
 Both statements are true already for the map $(T^*G_B)_{-f,f}\to A$. The analog of the first statement is true for any Hamiltonian space $X$ for a group $H$, endowed with a Lagrangian fibration $X \to Y$ with free $H$-action on $Y$, and the Hamiltonian reduction at an $H$-fixed point of $\mathfrak h^*$. (Take, here, $X=T^*G_B$ and $Y=G_B$.) Triviality of the monodromy follows from the sequence\footnote{I thank the anonymous referee for pointing out a mistake in a previous version of this proof and suggesting this argument.}
 \[(T^*G_B)_{-f,f} \to T^*A \to A,\]
 with the first arrow being an $N\times N$-torsor and the second an $\mathfrak a^*$-torsor (hence both contractible).
\end{proof}

Note, also that the fibers of $M^\circ_{-f,f}\to J^\circ$ are $N\times N$-torsors. The line bundle, with its connection, descend to $J^\circ$, and will also be denoted by $(L,\nabla)$. We will call the foliations $\mathscr G_1: J^\circ \to A_1$ and $\mathscr G_2: J^\circ\to A_2$ the ``vertical'' and ``horizontal'' foliations, respectively, but I emphasize that these are \emph{different} foliations than the ``vertical'' and ``horizontal'' foliations on $M$: their preimages in $M^\circ_{-f,f} \to J^\circ$ \emph{do not} coincide with leaves of the foliations on $M$.

\subsubsection{}

We are now ready to reformulate Jacquet's Theorem \ref{thmJacquet}. 
We will think of $\mathcal D(\End(V))$ and $\mathcal D(\End(V^*))$ as geometric quantizations of $M$, given by the data $(L,\nabla)$. This means that we consider functions on $\End(V)$ (resp.\ on $\End(V^*)$) as sections of $L$ which are flat along the leaves of the ``vertical'' (resp.\ ``horizontal'') foliation, as explained before, and we will be writing $\mathcal D(\End(V))= \mathcal D_\hor(M, L)$ (resp.\ $\mathcal D(\End(V^*))= \mathcal D_\ver(M, L)$), thinking of the elements as (Schwartz) half densities on the space of leaves for the corresponding foliation, valued in the line bundle of flat sections along the leaves. 

Fourier transform, then, becomes the isomorphism
\[\mathcal F: \mathcal D_\hor(M, L) \to \mathcal D_\ver(M, L)\]
given by the standard intertwiners of ``integration along the leaves of $\mathscr F_\hor$,'' as in Definition \ref{deffoliation}. Note that the inverse Fourier transform is also given by integration along the leaves, this time of $\mathscr F_\ver$.

\subsubsection{} \label{sssdefpushf}

The pushforwards \eqref{pushfdensities} of half-densities (restricted to $M^\circ$) can now be seen as maps  
\begin{equation}\label{intM}
\begin{aligned} \mathcal D_\hor(M, L) \to \mathcal D_\hor(J^\circ, L),  \\
 \mathcal D_\ver(M, L) \to \mathcal D_\ver(J^\circ, L).
\end{aligned}
\end{equation}

We can understand these maps as integrals along the Lagrangian leaves of the map $M^\circ_{-f,f}\to A$ (where $A=A_1$, resp.\ $A_2$), but some care is required in understanding those integrals, since the expression \eqref{intfoliation} does not make sense in the absence of a Lagrangian foliation on the entire space $M$. 

Rather, what we should do is pick a Haar measure on $N\times N$, and use it (or rather, its square root half-density) to integrate the given half-densities along the fibers of the map $M^\circ_{-f,f} \to J^\circ$, which are $N\times N$-orbits.

In more detail:
Consider the Lagrangian leaves of the foliation $\mathscr G_\ver: J^\circ\to A_1$. By Lemma \ref{lemmalagrangian}, these are quotients by the (free) $N\times N$-action of Lagrangian subvarieties of $M^\circ$; we will denote by $\mathscr G_{\ver,a}$ the leaf of $\mathscr G_\ver$ over $a\in A_1$, and by $\widetilde{\mathscr G}_{\ver,a}$ its preimage in $M^\circ$. 
The intersection of every $\widetilde{\mathscr G}_{\ver,a}$ with the leaves of the ``vertical'' foliation $\mathscr F_\ver: M\to\End(V)$ is easily seen to be of dimension $n$ ($=\dim A$). The orbits of $N\times N$ provide sections for each quotient 
\[ \widetilde{\mathscr G}_{\ver,a} \to \widetilde{\mathscr G}_{\ver,a}/\mathscr F_{\ver}.\]

The elements $\varphi \in \mathcal D_\hor(M, L)$ are half-densities in the transverse direction to $\mathscr F_\ver$ (valued in the line bundle $L_\ver$ of vertically parallel sections of $L$). The analog of \eqref{intfoliation}, here, is an integral of $\varphi$ over the quotient $\widetilde{\mathscr G}_{\ver,a}/\mathscr F_\ver$. For such an integral to make sense, \emph{we do not multiply $\varphi$ by a power of $|\omega|$, as in \eqref{intfoliation}, but by our fixed Haar half-density $|d\underline{n}|^\frac{1}{2}$ along $N\times N$-orbits}. Using an isomorphism analogous to \eqref{isomlines}, the product $\varphi \cdot |d\underline{n}|^\frac{1}{2}$ can be written as the product of the Haar measure $|d\underline{n}|$ on $\widetilde{\mathscr G}_{\ver,a}/\mathscr F_\ver$ by a half-density along $J^\circ/\mathscr G_\ver$ and by a section of $L$. 

The reader can check, using Remark \ref{remarkhalfdensities}, that this definition of the pushforward \eqref{intM} corresponds to the one that we gave in coordinates in \eqref{pushfdensities}. The definition depends on the choice of Haar measure on $N\times N$, but this choice applies to both sides of Theorem \ref{thmJacquet}, and will not affect the Hankel transforms. As we pointed out in that remark, although it is natural, this definition of pushforward is not entirely justified yet -- one could imagine varying the $N\times N$-Haar measure along the orbits. Its full justification will appear in Proposition \ref{proposition-main} below.

\subsubsection{}

Jacquet's theorem \ref{thmJacquet}, now, admits the following simple reformulation:
\begin{theorem}\label{thmreform}
The diagram
  \[\xymatrix{
 \mathcal D_\hor(M, L) \ar[d]\ar[r]^{\mathcal F}& \mathcal D_\ver(M, L)\ar[d]\\
 \mathcal D_\hor(J^\circ, L) \ar[r]^{\mathcal H_{\Std}} & \mathcal D_\ver(J^\circ, L)
 }\]
commutes, where the horizontal arrows are given by the integrals \eqref{intfoliation} along the leaves of the foliations.
\end{theorem}

A pleasant feature of this reformulation is that it does not require any arbitrary choices, such as the choice of the subset $wA\subset G$ for representing the orbital integrals \eqref{orbital-Kuz} of the Kuznetsov quotient $(N,\psi)\backslash G/(N,\psi)$.

\subsection{The case of $\GL_2$}\label{ssGL2}

We will use the case $n=2$, both to explicitly verify that the Hankel transforms $\mathcal H_{\Std}$ described by Theorems \ref{thmJacquet} and \ref{thmreform} coincide. The proof of the general case will only be sketched in the next subsection; it uses an inductive argument where every step is almost identical to the case of $\GL_2$.  The explicit verification that Theorem \ref{thmreform} amounts to Jacquet's formula \eqref{HStd} in higher rank is also similar to the case of $\GL_2$, and will be left to the reader; the validity of both theorems proves that, indirectly.

\subsubsection{Verification of \eqref{HStd} for $n=2$} 

Taking the description of Theorem \ref{thmreform} for the operator $\mathcal H_\Std$, let us see that it is given by \eqref{HStd}. 

First of all, let us calculate the set $M^\circ_{-f,f}$, and the (free) $N\times N$-action on it. It is not hard to see that this set consists of all elements of the form 
\begin{multline}\label{NNaction}
\left( w \begin{pmatrix} a_1 & \\  & a_2 \end{pmatrix}, \begin{pmatrix} b_1 + a_1^{-2}b_2^{-1} & a_1^{-1}\\ a_1^{-1} & b_2\end{pmatrix}w \right)\cdot  \left(\begin{pmatrix} 1 & y \\  & 1\end{pmatrix} \begin{pmatrix} 1 & x\\  & 1 \end{pmatrix} \right)=\\
 \left( w \begin{pmatrix} a_1 & a_1 x \\ - a_1 y & a_2 - a_1 xy \end{pmatrix}, \begin{pmatrix} b_1 +a_1^{-2}b_2^{-1} -  a_1^{-1}(x-y) - b_2 xy& a_1^{-1}-x b_2 \\ a_1^{-1}+y b_2 & b_2 \end{pmatrix}w \right),
\end{multline}
with $a_1,a_2,b_1,b_2\ne 0$. The symplectic form \eqref{sympM} on $M$ restricts to the form 
\begin{equation}\label{sympJ}\omega_J = db_1\wedge da_1 + db_2 \wedge da_2 + a_1^{-2}b_2^{-2} db_2\wedge da_1
\end{equation}
on $M_{-f,f}^\circ$, which descends to a symplectic form on $J$. Note that 
\[\omega_J\wedge \omega_J = - db_1 \wedge db_2 \wedge da_1\wedge da_2.\]

Next, we 
describe the elements of $\mathcal D_\hor(J^\circ, L)$: Denoting, as before, by $\mathscr G_\ver$ the ``vertical'' foliation $J^\circ\to A_1$, and by $\widetilde{\mathscr G}_\ver$ the corresponding foliation $M^\circ_{-f,f}\to J^\circ\to A$ of Lemma \ref{lemmalagrangian}, we have the following description, whose proof will be left to the reader.

\begin{lemma}\label{lemma-parallel}
The $\mathscr G_\ver$-parallel sections of $L$ can be described as sections $(\Phi,\Phi^*)$ of $L$ over $M^\circ_{-f,f}$, such that $\Phi$ depends only on the projection to $\End(V)$, and $\Phi(m\cdot (n_1,n_2)) = \psi(n_1^{-1} n_2) \Phi(m)$, for all $m$. (This completely determines $\Phi^*$, by the defining property \eqref{phiphistar} of the line bundle $L$.)  Thus, the elements of $\mathcal D_\hor(J^\circ, L)$ are half-densities on $A_1$ valued in this line bundle. Similarly, elements of $\mathcal D_\ver(J^\circ, L)$ are half-densities on $A_2$ valued in the line bundle whose sections are given by pairs $(\Phi,\Phi^*)$ on $M_{-f,f}^\circ$ such that $\Phi^*$ depends only on the projection to $\End(V^*)$, and $\Phi^*(m\cdot (n_1,n_2)) = \psi(n_1^{-1} n_2) \Phi^*(m)$. (Again, this completely determines $\Phi$.)
\end{lemma}

The elements of $\mathcal D^-_{L(\Std,\frac{1}{2})}(\mathfrak Y)$, $\mathcal D^-_{L(\Std^\vee,\frac{1}{2})}(\mathfrak Y)$ of Theorem \ref{thmJacquet} are obtained from the elements of $\mathcal D_\hor(J^\circ, L)$, $\mathcal D_\ver(J^\circ, L)$ of Theorem \ref{thmreform} by evaluating $\Phi$, resp.\ $\Phi^*$, at the elements of $M^\circ_{-f,f}$ living over representatives $wA\subset G\hookrightarrow \End(V)$ and  $wA\subset G\hookrightarrow \End(V^*)$, respectively. Such representatives are given by the pairs 
\begin{equation}\label{representatives1} A\ni (a_1, a_2)  \mapsto  \left(w \begin{pmatrix} a_1 \\ & a_2 \end{pmatrix},  \begin{pmatrix} b_1 + a_1^{-2} b_1^{-1} &  a_1^{-1} \\ a_1^{-1} & b_2 \end{pmatrix}w\right) \in M^\circ_{-f,f}
\end{equation}
(with arbitrary $b_1, b_2 \in F^\times$), and 	
\begin{equation}\label{representatives2} A\ni (b_1, b_2)  \mapsto  \left( w\begin{pmatrix} a_1 & b_2^{-1} \\ b_2^{-1} & a_2 + a_1^{-1} b_2^{-2} \end{pmatrix} , \begin{pmatrix} b_1 &  \\  & b_2 \end{pmatrix}w\right) \in M^\circ_{-f,f}
\end{equation} (with arbitrary $a_1, a_2 \in F^\times $).

Now we compute the Hankel tranform of such an element $\varphi\in \mathcal D_\hor(J^\circ, L)$, according to Theorem \ref{thmreform}. 

Writing $(\underline{a}, \underline{b}) = (a_1,a_2, b_1,b_2)$ and $\varphi (\underline{a}, \underline{b}) = (\Phi,\Phi^*) |d\underline{b}|^\frac{1}{2}$, we have $\mathcal H_{\Std}\varphi = (\Psi,\Psi^*) |d\underline{a}|^\frac{1}{2}$, with 
\begin{multline*} \Psi^*\left( w\begin{pmatrix} * & b_2^{-1} \\ b_2^{-1} & * \end{pmatrix} , \begin{pmatrix} b_1 &  \\  & b_2 \end{pmatrix}w\right) = \\
\int \Phi^* \left( w\begin{pmatrix} a_1 & b_2^{-1} \\ b_2^{-1} & a_2 + a_1^{-1} b_2^{-2} \end{pmatrix} , \begin{pmatrix} b_1 &  \\  & b_2 \end{pmatrix}w\right) |da_1\, da_2| = \\
\int \Phi \left( w\begin{pmatrix} a_1 & b_2^{-1} \\ b_2^{-1} & a_2 + a_1^{-1} b_2^{-2} \end{pmatrix} , \begin{pmatrix} b_1 &  \\  & b_2 \end{pmatrix}w\right) \psi(-a_1 b_1 - a_2 b_2 - a_1^{-1}b_2^{-1})  |da_1\, da_2| = \\
\int \Phi \left( w\begin{pmatrix} a_1 &  \\  & a_2 \end{pmatrix} , \begin{pmatrix} *  & a_1^{-1} \\ a_1^{-1} & * \end{pmatrix}w\right) \psi(-a_1 b_1 - a_2 b_2 + a_1^{-1}b_2^{-1} )  |da_1\, da_2|.
\end{multline*}

This is \eqref{HStd-expl}, 
taking into account that the embedding $\GL_n \hookrightarrow \Mat_n^* \xrightarrow\sim \Mat_n$ takes $\begin{pmatrix} b_1 &  \\  & b_2 \end{pmatrix}w$ to $w\begin{pmatrix} b_1^{-1} &  \\  & b_2^{-1} \end{pmatrix}$.

\subsubsection{Proof for $n=2$}

The proof of Theorem \ref{thmreform} for $\GL_2$ will occupy the rest of \S~\ref{ssGL2}.
Let us rephrase what needs to be proven: We have Lagrangian foliations $\mathscr F_\ver$, $\mathscr F_\hor$ on $M$, and also two Lagrangian foliations $\widetilde{\mathscr G}_\ver$, $\widetilde{\mathscr G}_\hor$ on its subspace $M^\circ_{-f,f}$. Our goal is to show that, starting from a ``horizontal'' half-density $\varphi$ (i.e., parallel along $\mathscr F_\ver$), we obtain the \emph{same} half-density that is parallel along $\widetilde{\mathscr G}_\hor$, either by integrating $\varphi$ 
directly\footnote{It is immediate to see that the composition \rotatebox[origin=c]{180}{$\Lsh$} in the diagram of Theorem \ref{thmreform} is the same as directly integrating over the leaves of $\widetilde{\mathscr G}_\hor$.} over the leaves of $\widetilde{\mathscr G}_\hor$, or by first integrating it over the leaves of $\mathscr F_\hor$ (= Fourier transform), and then integrating it over the leaves of $\widetilde{\mathscr G}_\hor$. 

Such relations would be immediate, by the theory of the Weil representation (up to a certain factor that has to do with the 8-fold metaplectic cover), if the foliations denoted by $\widetilde{\mathscr G}$ were \emph{linear}, by which we mean each fiber to be open dense in an affine subspace of $M$. The fact that it remains true for the nonlinear foliations of the theorem is remarkable, but the proof will take advantage of such statements (``Weil's formula'') in the linear case. The main observation behind the proof is that \emph{there is another, linear Lagrangian foliation} $\widetilde{\mathscr G}$ of $M^\circ_{-f,f}$, that will serve as an intermediary between  $\widetilde{\mathscr G}_\ver$ and $\widetilde{\mathscr G}_\hor$.

Let us recall once more that the fibers of $M^\circ_{-f,f}\to J^\circ\simeq A_1\times A_2$ are $N\times N$-orbits. 
We will write $(\ua,\ub)$ for an element of $A_1\times A_2$ (the underline because we will sometimes explain things in coordinates, and think of $\ua$ as the pair $(a_1,a_2)$, and similarly for $\ub$), and $\mathcal O_{\ua,\ub}$ for its fiber. Using $\approx$ for two sets whose intersection is open dense in both, we have, in coordinates, 
\[ \bigcup_{a_2, b_1} \mathcal O_{\ua,\ub} \approx M_{a_1, b_2},\]
where $M_{a_1,b_2}$ denotes the affine subspace 
\begin{equation}\label{Mab} \left\{\left( w \begin{pmatrix} a_1 & a_1 x \\ - a_1 y & a_2 \end{pmatrix}, \begin{pmatrix} b_1 & a_1^{-1}-x b_2 \\ a_1^{-1}+y b_2 & b_2\end{pmatrix}w \right),\,\, (a_2,b_1, x, y)\in F^4\right\}.\end{equation}

Coordinate-independently, we have made a choice of unipotent subgroup $N$, hence of the line $L=V^N\subset V$, and these affine subspaces are the fibers of the composition
\[ M^\circ_{-f,f}\hookrightarrow \End(V)\times \End(V^*) \to \Hom(L, V/L) \times \Hom(L^\perp, V^*/L^\perp),\] 
over pairs of nonzero morphisms in codomain. In what follows, we also use $(a_1, b_2)$ to denote such a pair in  $\Hom(L, V/L) \times \Hom(L^\perp, V^*/L^\perp)$. These spaces $M_{a_1, b_2}$ (or rather, their open intersections $M_{a_1, b_2}^\circ$ with $M^\circ_{-f,f}$), as the pair $(a_1,b_2)$ varies, form a Lagrangian foliation $\tilde{\mathscr G}$ of $M^\circ_{-f,f}$.  
We will write $\mathscr G$ for the corresponding foliation of the quotient $J^\circ = M^\circ_{-f,f}/(N\times N)$. 

\subsubsection{Integrals over Lagrangians}\label{sssLagrangians}

We now come to one of the finest, albeit elementary, parts of the argument, which will also answer the question that we posed in Remark \ref{remarkhalfdensities}: why is the chosen way to define pushforward of half-densities the ``correct'' one? We start by asking the question: \emph{If $\mathscr F$ is a Lagrangian foliation on $M$, what does it mean to average an element of $\mathcal D_{\mathscr F}(M,L)$ over a Lagrangian subspace $\mathscr L$?}

In the context where $\mathscr L$ is part of another Lagrangian foliation, the answer is given by formula \eqref{intfoliation} and Remark \ref{remarkint}. \emph{It is important to observe that these do not make sense with knowledge of a single Lagrangian subspace $\mathscr L$, without the foliation.} Indeed, the foliation allows the identification of \emph{all} conormal spaces to points of $\mathscr L$ with the cotangent fiber of their image in $M/\mathscr F$ (here, $\mathscr F_1=\mathscr F$), giving a canonical, up to scalar, trivialization of the determinant of the conormal bundle over any leaf of $\mathscr F$. 

However, in our context, we need to integrate over Lagrangians which are not part of a foliation of the entire space $M$, but just of its subspace $M^\circ_{-f,f}$. Let us call such a foliation of a (non-dense) subspace a ``partial'' foliation. In this context, as we have seen, \emph{the integrals only make sense after we pick a Haar half-density on each $N\times N$-orbit}, in order to first push forward to the symplectic variety $M^\circ_{-f,f}/(N\times N) = J^\circ$, where our Lagrangians do form a foliation. The choice of Haar half-densities	 on the $N\times N$-orbits affects the answer; understanding which choice is right is the question we posed in Remark \ref{remarkhalfdensities}. That remark gave an interpretation to our choice, but it remains to be seen why this is the correct choice.

For the proof of Theorem \ref{thmreform}, it is particularly important to compare our integrals over the fibers $M_{a_1,b_2}$ of the partial foliation $\widetilde{\mathscr G}$ with the corresponding integrals for the linear foliations containing them. Note that these affine subspaces are not parallel to each other -- they belong to different linear foliations. But if we fix a pair $(a_1,b_2)$, we can consider the linear foliation $\mathscr F = \mathscr F_{a_1,b_2}$ of parallel Lagrangians to it. We then have two versions of the ``integral over $M_{a_1,b_2}$'' of an element of $\mathcal D_\hor(M, L)$ (or of $\mathcal D_\ver(M, L)$ -- but we present the former):
\begin{enumerate}[(a)]
 \item The integral corresponding to the foliation $\mathscr F$. Since $\dim M_{a_1,b_2}/\mathscr F_{\ver} = 3$, this is given by \eqref{intfoliation}, with the exponent of $|\omega|$ being $\frac{3}{2}$. It produces a section of the line bundle multiplied by a half-density on the $4$-dimen\-sional quotient space $M/\mathscr F$. 
 \item The integral arising as the pushforward to $\mathcal D_\hor(J^\circ,L)$ (``integration over $N\times N$-orbits,'' depending on our choices of Haar half-densities on those), followed by the integral \eqref{intfoliation} along the fibers of the foliation $\mathscr G$. Here, since $\mathscr G_x/\mathscr G_\ver$ has dimension $1$, the exponent of $|\omega|$ is $\frac{1}{2}$. It produces a section of the line bundle multiplied by a half-density on the quotient $J^\circ/\mathscr G$. 
\end{enumerate}

How can we compare the outcomes of the two integrals? It clearly does not make sense to say that their restrictions on $M_{a_1,b_2}$ are ``equal,'' since the meaning of these restrictions is different: Both have a factor that can be intepreted as a parallel section of $L$ along this Lagrangian, but in the first case this is multiplied by a ``transverse'' half-density for the quotient $M/\mathscr F$ (which is $4$-dimensional), while the latter is multiplied by a ``transverse'' half-density for the quotient $J^\circ/\mathscr G$ (which is $2$-dimensional). These half-densities, restricted to $M_{a_1,b_2}$,  live  in the line bundle $|\det\mathcal N^*|^{\frac{1}{2}}$, where $\mathcal N^*$ denotes the conormal bundle in each case. 

The quotient of the normal bundle of $M_{a_1,b_2}^\circ$ in $M$ by its normal bundle in $M^\circ_{-f,f}$ is identified with the tangent space of $(f,f)$ in $\mathfrak n^* \times \mathfrak n^*$ via the moment map, hence the cokernel of the map of normal bundles can be identified with $\mathfrak n^* \times \mathfrak n^*$, hence: 
\begin{lemma}
	We have a canonical isomorphism of line bundles over $M_{a_1,b_2}^\circ$,
	\begin{equation}\label{factordeterminants}
	\det \mathcal N^*_{M/\mathscr F} = \det \mathcal N^*_{J^\circ/\mathscr G} \otimes (\det \mathfrak n)^{\otimes^2},
\end{equation}
	where the indices specify in which manifold the conormal bundles are taken. 
\end{lemma}

Let $\mathscr \mathcal D_{M/\mathscr F}$ be the line bundle of half-densities on $M/\mathscr F$, pulled back to $M^\circ_{a_1,b_2}$. Let $\mathscr \mathcal D_{J^\circ/\mathscr G}$ be the line bundle of half-densities on $J^\circ/\mathscr G$, pulled back to the same leaf. Then, \eqref{factordeterminants} implies that we have a canonical isomorphism of line bundles:
\begin{equation}\label{factordensities}
	\mathscr \mathcal D_{M/\mathscr F} = \mathscr \mathcal D_{J^\circ/\mathscr G} \otimes |\det \mathfrak n|.
\end{equation}

It follows that \emph{it makes sense to compare the restriction of an element of $\mathcal D_{\mathscr F}(M,L)$ to $M_{a_1,b_2}$ with the restriction of an  element of $\mathcal D_{\mathscr G}(J^\circ, L)$ multiplied by a section of $|\det(\mathfrak n \times \mathfrak n)|^\frac{1}{2} = |\det \mathfrak n|$.}

Now, the canonical, up to scalar, trivializations of $\det \mathcal N^*$ (for either of the foliations) along $M^\circ_{a_1,b_2}$ can be thought of as flat connections with trivial monodromy on these line bundles. The line bundle with fiber $(\det \mathfrak n)^{\otimes^2}$ also has such a connection, of course (whose parallel sections are the constant ones). The answer to the question posed in Remark \ref{remarkhalfdensities} lies in the following.

\begin{proposition}\label{proposition-main}
	The canonical isomorphism \eqref{factordeterminants} of line bundles on $M_{a,b}^\circ$ is an isomorphism of local systems, i.e., it preserves the canonical lines of sections. 
\end{proposition}

\begin{proof}
	After this abstract discussion, the reader will probably appreciate an explicit calculation. In coordinates $w\begin{pmatrix} A & B \\ C & D \end{pmatrix} \times \begin{pmatrix} A' & B' \\ C' & D' \end{pmatrix}w$ for $M$, coordinates for the linear quotient $M/\mathscr F$ are given by $(A, D', D' B + A B', D' C + A C')$, hence $\det \mathcal N^*_{M/\mathscr F}$ is trivialized, on $M_{a_1,b_2}^\circ$, by the section 
	\begin{multline*}\eta_{\mathscr F} := dA \wedge dD' \wedge d(D' B + A B') \wedge d(D' C + A C') \\ = dA \wedge dD' \wedge d(b_2 B + a_1 B') \wedge d(b_2 C + a_1 C').\end{multline*}
	Coordinates for $J^\circ/\mathscr G$ are given by $(A,D')$ (recall that now we are restricting to the subset $M_{-f,f}$), hence $\det \mathcal N^*_{J^\circ/\mathscr G}$ is trivialized, up to scalar, by the section 
	\[ \eta_{\mathscr G} := dA \wedge dD'.\]

We could express $\eta_{\mathscr F}$ as the wedge of $\eta_{\mathscr G}$ with differentials factoring through the moment map to $\mathfrak n^*\times \mathfrak n^*$, but equivalently we can also use the symplectic duality to translate these sections to sections of the determinant of the cotangent bundle of $M^\circ_{a_1,b_2}$: Using the symplectic form \eqref{sympM}, the dual of $\eta_{\mathscr F}$ is (up to scalar) the form  $dD \wedge dA' \wedge d(b_2 B - a_1 B') \wedge d(b_2 C - a_1 C')$ on $M^\circ_{a_1,b_2}$, which in coordinates \eqref{Mab} translates (up to scalar) to the form 
\[ \eta_{\mathscr F}^*:= db_1 \wedge da_2 \wedge dx \wedge dy.\] 
Using the symplectic form \eqref{sympJ} on $J^\circ$, the dual of $\eta_{\mathscr G}$ is (up to sign) the form 
\[ \eta_{\mathscr G}^*:= db_1 \wedge da_2.\] 

Hence, we see that $\eta_{\mathscr F}^*$ is the product of $\eta_{\mathscr G}^*$ by a fixed Haar volume form on $N\times N$, which proves the claim.
\end{proof}

\subsubsection{Completion of the proof of Theorem \ref{thmreform} for $n=2$}

We keep fixing a pair $(a_1,b_2)$, and denoting by $\mathscr F$ the corresponding ``linear'' foliation of $M$.

\begin{proposition}\label{prop-linear}
 The integrals \eqref{intfoliation} along the leaves of the following Lagrangian foliations give rise to a commutative diagram  
 \[\xymatrix{
 \mathcal D_\hor(M, L) \ar[dr]\ar[rr]^{\mathcal F}&& \mathcal D_\ver(M, L)\ar[dl]\\
 & \mathcal D_{\mathscr F}(M, L).
 }\] 
\end{proposition}

\begin{proof}
In the theory of the Weil representation, such diagrams of integrals along ``linear'' foliations  commute up to an $8$-th root of unity, see \cite[\S~3.5]{WWLi}. For the case at hand, this root of unity is trivial; this is ``Weil's formula,'' whose (very simple) proof is recalled in \cite[Proposition 2]{Jacquet}, and we leave the verification to the reader.
\end{proof}

We can now complete the proof of Theorem \ref{thmreform} for the case of $\GL_2$. Using a fixed Haar measure on $N\times N$ to define the pushforwards \eqref{intM}, we have a diagram 
\begin{equation}\label{intdiagram} \xymatrix{ 
 \mathcal D_\hor(M, L) \ar[d] \ar[dr]\ar[rr]^{\mathcal F}&& \mathcal D_\ver(M, L)\ar[dl] \ar[d]\\
\mathcal D_\hor(J^\circ, L)\ar[r]^\alpha & \mathcal D_{\mathscr G}(J^\circ, L) & \mathcal D_\ver(J^\circ, L) \ar[l]_\beta,
}
\end{equation}
where all the arrows are given by integrals along the Lagrangian leaves of each foliation. In particular, the maps to $\mathcal D_{\mathscr G}(J^\circ, L)$ are given by the integrals over the leaves of the foliation $\widetilde{\mathscr G}$ on $M^\circ_{-f,f}$. We discussed in \S~\ref{sssLagrangians} how to make sense of these integrals. Proposition \ref{proposition-main} implies that the evaluations of those on each leaf $M_{a_1,b_2}$ are \emph{equal}, up to a fixed half-density on $|\mathfrak n^*\times\mathfrak n^*|$ which can be taken to be the one corresponding to the chosen Haar half-density on $N\times N$, to the ``linear'' integral along $M_{a_1,b_2}$ corresponding to the linear foliation $\mathscr F = \mathscr F_{a_1,b_2}$. Proposition \ref{prop-linear}, now, implies that the upper triangle commutes. The two lower triangles also commute, by construction; i.e., the diagram \eqref{intdiagram} is commutative.

Finally, the inverse of the arrow labeled by $\beta$ is also the ``integral \eqref{intfoliation} over Lagrangians.'' This can be reduced to a usual, linear Fourier transform in dual $\underline{a}$ and $\underline{b}$ variables, using representatives as in \eqref{representatives1}, \eqref{representatives2}, and the symplectic form \eqref{sympJ}. Therefore, the composition $\beta^{-1} \circ \alpha$ is the map $\mathscr H_\Std$ of Theorem \ref{thmreform}.

\subsection{Sketch of the proof in the general case}

If one understands how the argument of Jacquet translates to the argument we presented in \S~\ref{ssGL2} for the case of $n=2$, it is straightforward to adapt it for general $n$. The only additional idea needed, already present in \cite{Jacquet}, is to apply induction in $n$. I will give a vague and impressionistic summary of the argument; most of the details remain to be filled in by the interested reader, who will also need to consult Jacquet's paper. This summary is not meant as a stand-alone account of the argument; its goal is simply to convince the reader that a translation to the setting of Theorem \ref{thmreform} is possible, even straightforward, given the argument for $n=2$.

The inductive step needed corresponds to the following factorization of the Kuznetsov orbital integrals \eqref{orbital-Kuz}
\begin{equation}\label{twosteps} O_a (\Phi) = \int_{U_n' \times U_n} \int_{N_{n-1}'\times N_{n-1}} \Phi(u_1 n_1 w a n_2 u_2) \psi^{-1}(u_1n_1n_2u_2) d(n_1,n_2) d(u_1,u_2),
\end{equation}
where, if $N$ is the unipotent subgroup of the Borel stabilizing a flag $V_1\subset V_2\subset\dots \subset V_n=V$ (which, in coordinates, to take to be the standard flag of $F^n$, with $N$ upper triangular),  $U_n$ is the unipotent radical of the stabilizer of $V_1$, and $N_{n-1}$ is the corresponding group for the $(n-1)$-dimensional space $V/V_1$, identified with a subgroup of $N$ by choosing a splitting of the quotient (which, in coordinates, we will do using the standard basis of $F^n$. The groups $U_n'$ and $N_{n-1}'$ are defined similarly, in terms of the dual filtration on the dual space. 

The inner integral of \eqref{twosteps} is then the Kuznetsov orbital integral for a function $\Phi_1$ on $\GL_{n-1}$, defined as 
\[\Phi_1 (g) = \int_{U_n' \times U_n} \Phi\left(u_1 w_n \begin{pmatrix} a_1 \\ & g \end{pmatrix}  u_2\right) \psi^{-1}(u_1 u_2)  d(u_1,u_2),\]
where $w_n$ is the permutation matrix $\begin{pmatrix} & I_{n-1} \\ 1 \end{pmatrix}$. An inductive application of the theorem, then, relates the Kuznetsov orbital integrals of $\Phi_1$ (appropriately normalized -- i.e., we need to work with half-densities again) to those of its Fourier transform, and applies the same argument, using Weil's formula (a direct generalization of Proposition \ref{prop-linear}) to relate the Fourier transform of $\Phi_1$ to the Fourier transform of $\Phi$. 

It is clear that this argument directly translates to prove our reformulation \ref{thmreform} of Jacquet's theorem. I will only make a few comments on the inductive step: Starting with the same foliations $\mathscr F_\ver: M\to \End(V)$ and $\mathscr F_\hor: M\to \End(V)$ as before, we can interpret the definition of $\Phi_1$ as an integral of an element of $\mathcal D_\hor(M,L)$ over appropriate Lagrangians corresponding to a foliation $\mathscr G_1$ of the symplectic reduction $U_n'\bbbslash_f M \ssslash_f U_n$ (or rather, of its  ``open Bruhat cell''). These Lagrangians are parametrized by pairs $(a_1, g)$ as above, which, coordinate-independently, can be thought of as invertible elements $a_1\in \Hom(V_1, V_n/V_{n-1})$ and $g\in \Hom(V/V_1, V_{n-1})$. Moreover, fixing $a_1$, these Lagrangians map to a similar to $\mathscr F_\ver$ linear foliation of (an open dense subset of) $T^* \Hom(V/V_1, V_{n-1})$, in a way that allows for an inductive application of the theorem to compute the Kuznetsov orbital integrals of the original function in terms of the Kuznetsov orbital integrals of the Fourier transform of $\Phi_1$. Finally, the relation between the Fourier transform of $\Phi_1$ and the Fourier transform of $\Phi$ (and its proof) can be thought of as the statement that a diagram of the form
\begin{equation}\label{intdiagram2} \xymatrix{ 
 \mathcal D_\hor(M, L) \ar[d] \ar[dr]\ar[rr]^{\mathcal F}&& \mathcal D_\ver(M, L)\ar[dl] \ar[d]\\
\mathcal D_{\mathscr G_1}(U_n'\bbbslash_f M \ssslash_f U_n, L)\ar[r]^\alpha & \mathcal D_{\mathscr G}(U_n'\bbbslash_f M \ssslash_f U_n, L) & \mathcal D_{\mathscr G_2}(U_n'\bbbslash_f M \ssslash_f U_n, L) \ar[l]_\beta
}
\end{equation}
commutes. This diagram is the analog of \eqref{intdiagram}, with $\mathscr G_2$ a foliation parametrized by pairs $(b_1,g')$, with $b_1, g'$ in the duals $\Hom(V_n/V_{n-1},V_1)$,  $\Hom(V_{n-1},V/V_1)$ of $\Hom(V_1,V_n/V_{n-1})$, $\Hom(V/V_1, V_{n-1})$, respectively. The argument for the proof of this remains essentially the same, with an intermediate ``linear'' foliation $\mathscr G$, parametrized by pairs $(a_1, g')$ (notation as before), where one can apply Weil's formula.

\bibliographystyle{alphaurl}
\bibliography{biblio}

\end{document}